\documentclass[a4paper,11pt]{article}

\usepackage{amsmath}
\usepackage{theorem}
\usepackage{amsfonts}
\usepackage{stmaryrd}
\usepackage{amssymb}

\usepackage[latin1]{inputenc}
\usepackage[all]{xy}

\usepackage{graphics}

\theoremstyle{plain}
\newtheorem{theo}{Theorem}[section]
\newtheorem{cor}[theo]{Corollary}
\newtheorem{prop}[theo]{Proposition}
\newtheorem{lem}[theo]{Lemma}

\newtheorem{fact}[theo]{Fact}

\newtheorem{deuf}[theo]{Definition}

\newtheorem{ps}[theo]{Remark}

\newenvironment{eg} {\noindent {\em \bf{Example}} \\ }{ \hfill \newline }
\newenvironment{preuve} {\noindent {\em \bf{Proof}} \\ }{\hfill $\square$ \\ }
\newenvironment{preuvede}{\noindent {\em \bf{Proof}}}{\hfill $\square$ \newline}

\newcommand{\dC}{\ensuremath{\mathbb{C}}}
\newcommand{\dD}{\ensuremath{\mathbb{D}}}

\newcommand{\dF}{\ensuremath{\mathbb{F}}}

\newcommand{\dH}{\ensuremath{\mathbb{H}}}

\newcommand{\dN}{\ensuremath{\mathbb{N}}}

\newcommand{\dR}{\mathbb{R}}
\newcommand{\dS}{\ensuremath{\mathbb{S}}}

\newcommand{\dZ}{\ensuremath{\mathbb{Z}}}

\newcommand{\masum}[2]{\displaystyle{\sum_{#1}^{#2}}}

\input{epsf}

\def \psl{PSL(2,\dR)}

\def \Homeop{\text{Homeo}^+}

\def \Tr{\text{\rm Tr}}

\def \vertinj{\cap}

\def \Id{\text{\rm Id}}

\def \so{SO_\dR(n,1)}
\def \sop{SO_\dR^+(n,1)}

\title{Connected components of the compactification of representation spaces of surface groups}

\author{Maxime Wolff}

\begin{document}

\begin{abstract}
  The Thurston compactification of Teichm\"uller spaces has been generalized
  to many different representation spaces by J. Morgan, P. Shalen, M. Bestvina,
  F. Paulin, A. Parreau and others. In the simplest case of representations of
  fundamental groups of closed hyperbolic surfaces in $\psl$, we prove that
  this compactification is very degenerated: the nice behaviour of the
  Thurston compactification of the Teichm\"uller space contrasts with
  wild phenomena happening on the boundary of the other connected components
  of these representation spaces. We prove that it is more natural to consider
  a refinement of this compactification, which remembers the orientation of the
  hyperbolic plane. The ideal points of this compactification are {\em fat}
  $\dR$-trees, \textsl{i.e.}, $\dR$-trees equipped with a planar structure.
\end{abstract}

\maketitle

\tableofcontents

\section{Introduction}
Let $\Gamma$ be a group with a given finite generating set $S$. In all this
text, we consider the space $R_\Gamma(n)=Hom(\Gamma,Isom^+(\dH^n))$ of
actions of $\Gamma$ on the real hyperbolic space of dimension $n$ by
isometries preserving the orientation. The set $R_\Gamma(n)$ naturally
embeds in $\left(Isom^+(\dH^n)\right)^S$, giving it a Hausdorff, locally
compact topology. The Lie groups $Isom^+(\dH^n)$ and $Isom(\dH^n)$ act on
$R_\Gamma(n)$ by conjugation, and we will consider the quotients
$R_\Gamma(n)/Isom(\dH^n)$ and $R_\Gamma(n)/Isom^+(\dH^n)$, equipped with
their quotient topologies.

We will mainly focus on the case when $\Gamma=\pi_1\Sigma_g$ is the
fundamental group of a closed surface of genus $g$, with a given standard
presentation (\textsl{i.e.}, a {\em marking}); we then denote
$R_g(n)=R_{\pi_1\Sigma_g}(n)$. Also, we are mainly interested in the case
$n=2$ (thus $Isom^+(\dH^2)=\psl$), and we denote $R_g=R_g(2)$.

The space $R_g$ is a real algebraic variety (see \cite{CullerShalen}),
and is smooth, of dimension $6g-3$, outside of the set of abelian
representations (see \cite{Goldman84}). Outside of the set of (classes of)
{\em elementary} representations (\textsl{i.e.}, having a global fixed point
in $\dH^2\cup\partial\dH^2$), the quotient $R_g/\psl$, equipped with the
quotient topology, is again smooth of dimension $6g-6$; one of its connected
components is naturally identified with the Teichm\"uller space of the
surface (see \cite{GoldmanThese,Goldman88}).

\smallskip

A function $e\colon R_g\rightarrow\dZ$, called the
{\em Euler class}, plays a key role in understanding the spaces
$R_g$ and $R_g/\psl$. There are several ways to define it; for instance,
a representation $\rho\in R_g$ defines a circle bundle on $\Sigma_g$, which
has a $\dZ$-valued characteristic class, the Euler class. We will give
complete recalls on this Euler class in section \ref{SectionEulerClass}, and
refer to
\cite{Milnor,Wood,Ghys,Matsumoto,Goldman88,Calegari} for a deep understanding
of this class.

In 1988, W. Goldman proved \cite{Goldman88} that $R_g$ possesses
$4g-3$ connected components, which are the preimages $e^{-1}(k)$, for
$2-2g\leq k\leq 2g-2$ (the fact that the absolute value of the Euler class
is bounded by $2g-2$ was previously known as the Milnor-Wood inequality,
see \cite{Milnor,Wood}). He proved, moreover, that the components
$e^{-1}(2-2g)$ and $e^{-1}(2g-2)$ consist precisely in {\em all} the discrete
and faithful representations of $\pi_1\Sigma_g$ into $\psl$. He also proved
that the space $R_g(3)=Hom(\pi_1\Sigma_g,PSL(2,\dC))$ has two connected components.
The elements of $R_g(2)$ of even Euler class, on one hand, and the ones of
odd Euler class on the other hand, fall into these two different components
of $R_g(3)$.

The Euler class is still well-defined in the quotient space $R_g/\psl$
(we denote $\mathbf{e}\colon R_g/\psl\rightarrow\dZ$), whence the space
$R_g/\psl$ also admits $4g-3$ connected components, similarly. However, only
the absolute value of the Euler class is defined on $R_g/Isom(\dH^2)$, which
has $2g-1$ connected components. In $R_g/\psl$, the connected component
$\mathbf{e}^{-1}(2-2g)$ is naturally identified with the Teichm\"uller
space of $\Sigma_g$, and $\mathbf{e}^{-1}(2g-2)$ is the space of marked
hyperbolic metrics of constant curvature $-1$ on the surface $\Sigma_g$
equipped with the opposite orientation.
These consist in all the classes of discrete and faithful representations,
given by the holonomies of such structures.

\medskip

In 1976, W. Thurston \cite{Thurston} introduced a natural compactification
of Teichm\"uller spaces; his construction was intensively studied and detailed
in \cite{FLP}. Recall that the Teichm\"uller space of the surface,
$\mathbf{e}^{-1}(2-2g)$, is the space of marked hyperbolic metrics on
$\Sigma_g$, of constant curvature $-1$. A natural system of coordinates, on
the Teichm\"uller space, consists in taking the lengths of closed geodesics on
the surface. In that way, one proves that the Teichm\"uller space is
homeomorphic to $\dR^{6g-6}$ (indeed, the lengths of $6g-6$ curves suffice
to parameter hyperbolic metrics). W. Thurston's compactification of the
Teichm\"uller space consists in embedding this space into a projective
space (thus, considering lengths only up to a scalar), in which the
Teichm\"uller space has relatively compact image. A very important feature of
this compactification is that the boundary added to the Teichm\"uller space
is homeomorphic to a sphere of dimension $6g-7$, in such a way that the
Thurston compactification of the Teichm\"uller space is homeomorphic to a
closed ball of dimension $6g-6$.

This compactification has been extended to other connected components of
$R_g/\psl$, and to other representation spaces in successive works.
In 1984, J. Morgan and P. Shalen
\cite{MorganShalen}, using techniques of algebraic geometry,
defined a compactification of the real algebraic variety
$X_{\Gamma,SL(2,\dR)}$ of characters of representations
of $\Gamma$ in $SL(2,\dR)$. We denote
$\overline{X_{\Gamma,SL(2,\dR)}}^{MS}$ this compactification, and
$\overline{X_{g,SL(2,\dR)}}^{MS}$ in the case $\Gamma=\pi_1\Sigma_g$.
Since all the representations of even Euler class in $\psl$
(and, in particular, elements of the Teichm\"uller space)
lift to $SL(2,\dR)$
(see \textsl{e.g.} \cite{Goldman88}), this defines a
compactification of
Teichm\"uller spaces, and it coincides with the one of W. Thurston's.
In \cite{Morgan}, J. Morgan generalized this construction to the group
$SO(n,1)$, for $n\geq 2$.
In 1988 \cite{Bestvina88,Paulin88}, M. Bestvina and F. Paulin, independently,
gave a much more geometric viewpoint to this compactification,
for representations in
$Isom^+(\dH^n)$, $n\geq 2$. F. Paulin proved that the quotient topology, on
the space $R_\Gamma(n)/Isom(\dH^n)$, coincides with the equivariant Gromov
topology (which we will discuss in detail in the section \ref{SectionPaulin}).
Equipped with this topology, the space
$m_\Gamma^{fd}(n)$ of (conjugacy classes of) faithful and discrete
representations has a natural compactification, which recovers the
compactification of \cite{Morgan,MorganShalen}.
Note, finally, that A. Parreau extended \cite{AnneComp} this compactification
to representation spaces in higher rank groups.

\bigskip

This compactification being well-defined, it is a natural question to ask if
it respects the topology and geometry of $R_g/\psl$, as does the
Thurston compactification, when restricted to the Teichm\"uller space alone.
We shall prove that it is not the case, and that this compactification, as
defined in all the works mentioned above, leads to a very degenerate space.

The works of J. Morgan and P. Shalen, of M. Bestvina and of F. Paulin all yield
the same compactification, and we shall follow the approach of F. Paulin,
in which it is easier to add a notion of orientation.
In order to use F. Paulin's construction, we first define {\em explicitly}
the biggest Hausdorff quotients of $R_\Gamma(n)/Isom(\dH^n)$ and
$R_\Gamma(n)/Isom^+(\dH^n)$, that we denote $m_\Gamma^u(n)$ and
$m_\Gamma^o(n)$, respectively (indeed, these spaces are not
Hausdorff in general, so that M. Bestvina and F. Paulin's constructions cannot
be extended literally to these spaces). Again, we write
$m_\Gamma^u=m_\Gamma^u(2)$, $m_g^u(n)=m_{\pi_1\Sigma_g}(n)$ and so on.
Note that, in the space $m_g^u$, only the absolute value of the Euler class
is still defined, and $m_g^u$ possesses $2g-1$ connected components.
We denote by $\overline{m_g^u(n)}$ the compactification of $m_g^u(n)$ as it
is constructed by F. Paulin.

Also, in all this text, for all $k\in\{2-2g,\ldots,2g-2\}$,
we will denote by $m_{g,k}^o$ the subset of $m_g^o$ consisting of classes of
representations of Euler class $k$, and for $k\in\{0,\ldots,2g-2\}$,
$m_{g,k}^u$ will denote the connected component of $m_g^u$  consisting of
classes of representations whose Euler class, in absolute value, equals $k$.

\begin{theo}\label{PaBoIntro}
  Let $g\geq 4$ and $k\in\{0,\ldots,2g-3\}$.
  Then, in $\overline{m_g^u}$, the boundary of the Teichm\"uller space,
  $\partial m_{g,2g-2}^u$, is contained in
  $\partial m_{g,k}^u$ as a
  closed, nowhere dense subset.
\end{theo}

In particular,
\begin{cor}\label{introcc}
  For all $g\geq 4$, $\overline{m_g^u(2)}$ is connected.
\end{cor}
We also prove that for $g=2$ and $g=3$, the space $\overline{m_g^u(2)}$
possesses at most two connected components. Actually, theorem \ref{PaBoIntro}
should hold for all $g\geq 2$, but in the case where $k=1$ the proofs presented
here request that $g\geq 4$.

In particular, the two connected components of
$Hom(\pi_1\Sigma_g,PSL(2,\dC))/PSL(2,\dC)$ meet at their boundary in this
compactification, as soon as $g\geq 3$:
\begin{cor}
  For all $g\geq 3$, the space $\overline{m_g^u(3)}$ is connected.
\end{cor}

If we consider only representations of even Euler class (\textsl{i.e.}, that
lift to $SL(2,\dR)$) then the result holds for all $g\geq 2$:
\begin{cor}\label{introcarac}
  For all $g\geq 2$,
  the space $\overline{X_{g,SL(2,\dR)}}^{MS}$ is connected.
\end{cor}

\smallskip

Theorem \ref{PaBoIntro} not only implies that the space $\overline{m_g^u}$
is connected, but it also implies that this space is particularly degenerated.
Since the connected components $m_{g,k}^u$ are of dimension $6g-6$
(see \cite{Goldman84}, section 1), one should expect that the boundary
$\partial m_{g,k}^u$ has dimension $6g-7$. However, the boundary of the
Teichm\"uller space itself has dimension $6g-7$ (see \textsl{e.g.}
\cite{FLP}, expos\'e 1, theorem 1). Therefore, it follows from theorem
\ref{PaBoIntro} that the compactifications $\overline{m_{g,k}^u}$ of
the ``exotic'' connected components (\textsl{i.e.}, the connected components
$m_{g,k}^u$ such that $|k|\leq 2g-3$)
have no PL structures, or cell complex structures or so on, compatible with
the compactification. This contrasts very strongly with the behaviour of
the Thurston compactification of the Teichm\"uller spaces.

\smallskip

The proof of theorem \ref{PaBoIntro} uses the following fact,
interesting for itself:
\begin{prop}\label{oneendintro}
  For all $k$ such that $|k|\leq 2g-2$, the connected components $m_{g,k}^o$
  and $m_{g,|k|}^u$ are one-ended (\textsl{i.e.}, all the connected components
  of $m_g^o$ and $m_g^u$ are one-ended).
\end{prop}
Actually, for all $k\neq 0$, this follows from a theorem of N. Hitchin
(see \cite{Hitchin}, proposition 10.2), which says that the connected
component $m_{g,k}^u$ is homeomorphic to a complex vector bundle over the
$(2g-2-|k|)$-th symmetric product of the surface $\Sigma_g$. However,
the proof of proposition \ref{oneendintro} is far more simple and extends
to the case $k=0$.

\medskip

The degeneration of the compactification of
$m_g^u$ is due (at least) to the fact that the equivariant Gromov topology
forgets the orientation of the space
$\dH^n$. This information, in the case $n=2$, is the information carrying
the Euler class, and which separates the space
$m_g^u$ into its connected components.
By restauring this orientation we will define a new compactification of
these representation spaces, cancelling (partly, at least) this
degeneration.

We define a notion of convergence in the sense of Gromov for
{\em oriented} spaces, which preserves the orientation.
This enables us to define a new compactification of
$m_\Gamma^o$,
in which the ideal points, added at the boundary, are
{\em fat} $\dR$-trees, that is, $\dR$-trees together with an ``orientation''.
These fat $\dR$-trees form a set
$\mathcal{T}^o(\Gamma)$, which compactifies the space $m_\Gamma^o$:

\begin{theo}\label{introordonnable}
  The map $m_\Gamma^o\rightarrow m_\Gamma^o\cup\mathcal{T}^o(\Gamma)$ induces
  a natural compactification of
  $m_\Gamma^o$. Moreover, the natural map
  $\pi:\overline{m_\Gamma^o}\rightarrow\overline{m_\Gamma^u}$,
  which consists in forgetting the orientation, is onto, and some of its fibres
  have the same cardinality as $\dR$.
\end{theo}
By {\em natural}, we mean, following F. Paulin \cite{Paulin04}, that the
action of the mapping class group extends continuously to the compact space
$\overline{m_g^o}$.

One can define an Euler class in a quite general context, as we will see in
the section \ref{SectionEulerClass}.
In particular, the actions of $\pi_1\Sigma_g$ on fat $\dR$-trees admit an
Euler class, and we shall prove the following:

\begin{theo}\label{EulerContinuIntro}
  The Euler class $\mathbf{e}:\overline{m_g^o}\rightarrow\dZ$ is a continuous
  function. In particular, the
  compactification $\overline{m_g^o}$ possesses as many connected components
  as the space $m_g^o$.
\end{theo}

In theorem \ref{introordonnable}, the fact that some fibres of the map $\pi$
have the same cardinality as $\dR$ will depend on the following remarkable
property
of rigidity of non discrete representations of $\pi_1\Sigma_g$ in $\psl$.
Let $o$ denote the natural cyclic order on $\partial\dH^2=\dS^1$,
\textsl{i.e.}, for all $a,b,c\in\dS^1$, put $o(a,b,c)=0$ if two of these
three points coincide, $1$ if $a,b,c$ are met in the anticlockwise order on
the circle, and $-1$ otherwise.
\begin{prop}\label{introDiscretesRigides}
  Let $\rho_1,\rho_2\in R_g$ be two non discrete, non elementary
  representations, and let $a_1,a_2\in\partial\dH^2$. Suppose that for all
  $\gamma_1,\gamma_2,\gamma_3\in\pi_1\Sigma_g$,
  \[ o(\rho_1(\gamma_1)\cdot a_1,\rho_1(\gamma_2)\cdot a_1,
  \rho_1(\gamma_3)\cdot a_1)=o(\rho_2(\gamma_1)\cdot a_2,
  \rho_2(\gamma_2)\cdot a_2,\rho_2(\gamma_3)\cdot a_2). \]
  Then there exists $g\in\psl$ such that $\rho_1=g\rho_2 g^{-1}$.
\end{prop}

In other words, two representations are conjugate if and only if they have
orbits which are isomorphic as cyclically ordered sets. It follows that two
non elementary, non discrete representations are conjugate by an element of
$\psl$ if and only if they are conjugate by an element of $Homeo^+(\dS^1)$,
as representations in $Homeo^+(\dS^1)$. Since all the orbits of these
representations are dense in $\dS^1$, it follows
from proposition 2-2 of E. Ghys \cite{Ghys87} that two non elementary, non
discrete representations are conjugate in $\psl$ if and only if they are
topologically semi-conjugate, which is equivalent to saying that they
define the same element in the bounded cohomology group
$H_b^2(\pi_1\Sigma_g,\dZ)$ (see \cite{Ghys87}).

This contrasts drastically with the case of representations of maximal
Euler class. A theorem of S. Matsumoto \cite{Matsumoto} states that
if two representations of $\pi_1\Sigma_g$ in $Homeo^+(\dS^1)$ both have
Euler class $2g-2$ then they are conjugate. In other words, in the space
$Hom(\pi_1\Sigma_g,Homeo^+(\dS^1))/Homeo^+(\dS^1)$, the subspace
$\mathbf{e}^{-1}(2g-2)$ is a single point, whereas $\mathbf{e}^{-1}(k)$, for
$|k|\leq 2g-3$, contains the injective image of a space of dimension $6g-6$
(indeed, it is proved in \cite{noninjrep}, proposition 4.4, see also
\cite{TheseMaxime}, proposition 4.1.15) that the space of non elementary, non
discrete representations is an open and dense subspace of each connected
component of $Hom(\pi_1\Sigma_g,\psl)/\psl$ of non extremal Euler class.

\bigskip

By theorem \ref{introordonnable}, it is a necessary condition, for an action of
$\pi_1\Sigma_g$ on an $\dR$-tree to be in $\overline{m_g^u}$,
to preserve some orientation. In particular, we can construct explicit actions
on $\dR$-trees which are not the limit of actions of surface groups on $\dH^2$.
These explicit actions can even be obtained as limits of actions of surface
groups on
$\dH^3$, so that we get the following proposition:

\begin{prop}\label{degSkoraintro}
  Let $g\geq 3$.
  There exist minimal actions of $\pi_1\Sigma_g$ on $\dR$-trees by isometries,
  which are in $\partial\overline{m_{g}^u(3)}$
  but not in $\partial\overline{m_g^u(2)}$.
\end{prop}

This contrasts severely with the case of discrete and faithful representations.
Indeed, R. Skora \cite{Skora} proved that a minimal action of
$\pi_1\Sigma_g$ on an $\dR$-tree has small arc stabilizers
(\textsl{i.e.}, the stabilizer of any pair of distinct points of the tree
is virtually abelian)
if and only if it is the limit of discrete and faithful representations of
$\pi_1\Sigma_g$ in $\psl$ (or equivalently, it is a point in the boundary of
the Teichm\"uller space), and it is well-known
(see \textsl{e.g.} \cite{MorganShalen,Paulin88,Bestvina88}) that limits of
discrete and faithful representations in
$Isom^+(\dH^n)$ enjoy that property. As a corollary of R. Skora's result, we
thus have
$\partial\overline{m_{g}^{fd}(n)}=\partial\overline{m_g^{fd}(2)}$.
When we consider representations which may not be faithful and discrete,
proposition \ref{degSkoraintro} states that this equality does not hold any
more.

\medskip

Note that C. McMullen has developed, independently, a theory of oriented
$\dR$-trees, under the name of {\em ribbon $\dR$-trees}, in \cite{McMu},
in order to
compactify the set of proper holomorphic maps from the unit disk
$\Delta\subset\dC$ into itself.

\bigskip

This text is organized as follows. Section 2 gathers every background material
concerning the spaces that we wish to compactify.
Section 2.1 is devoted to the explicit construction of $m_\Gamma^u(n)$ and
$m_\Gamma^o(n)$, the biggest Hausdorff quotients of
$R_\Gamma(n)/Isom(\dH^n)$ and $R_\Gamma(n)/Isom^+(\dH^n)$. In section 2.2,
we give some recalls on F. Paulin's point of view on the compactification
of $m_\Gamma^u(n)$, while adapting it slightly so that it indeed defines a
compactification of the whole space $m_\Gamma^u(n)$. In section 2.3, we recall
the construction of the Euler class, and establish a technical lemma which
will enable us to prove that the Euler class extends continuously to the
boundary (in the oriented compactification) of $m_g^o$ (theorem
\ref{EulerContinuIntro}). Section 2.4 is devoted to the proof of proposition
\ref{introDiscretesRigides}.
Section 2.5 recalls an argument of Z. Sela in the
context of limit groups, which plays a key role in the proof of theorem
\ref{PaBoIntro}. Section 3 contains the core of our results on the
compactification. It begins with the proof of proposition \ref{oneendintro},
then turns to the degenerations of $\overline{m_g^u}$, and finally we
construct the oriented compactification.

\medskip

{\bf Acknowledgements.}
\small

This article, essentially, gathers the main results of my PhD thesis, directed
by Louis
Funar; he inspired this work from the beginning and suggested several ideas
developed here.
I am very grateful to Vincent Guirardel, Gilbert Levitt and Fr\'ed\'eric Paulin
who corrected some mistakes in earlier versions of this text, and whose
numerous remarks have deeply improved the presentation of this paper.
I would also like to thank Thierry Barbot, G\'erard Besson,
L\'ea Blanc-Centi, Damien Gaboriau, Sylvain Gallot,
Anne Parreau, Vlad Sergiescu,
for inspiring discussions and encouragements. I also wish to thank Bill
Goldman for sending me a copy of \cite{GoldmanThese}, and all the members of
the LATP in Marseille, where I finished the redaction of this paper.

\normalsize

\section{Preliminaries}\label{SectionPreliminaires}

\subsection{The representation spaces}

Let us first remark that the quotient spaces
$m_\Gamma'(n)=R_\Gamma(n)/Isom^+(\dH^n)$ and $R_\Gamma(n)/Isom(\dH^n)$ are
not Hausdorff in general. Indeed,
in $Isom^+(\dH^2)=PSL(2,\dR)$, the matrix
$\left(\begin{array}{cc}1 & 1 \\ 0 & 1\end{array}\right)$ is conjugate to
$\left(\begin{array}{cc}1 & \frac{1}{t^2} \\ 0 & 1\end{array}\right)$,
for all $t\in \dR^*$, hence its conjugacy class cannot be separated from the
one of the identity. As soon as there exists a morphism of
$\Gamma$ onto $\dZ$, we can therefore construct abelian representations of
$\Gamma$ in $PSL(2,\dR)$, and more generally in $Isom^+(\dH^n)$ (which
contains isomorphic copies of $PSL(2,\dR)$), which are not separated from
the trivial representation in these quotient spaces.

However, in order to use F. Paulin's construction of the compactification of
these representation spaces, we need to work with Hausdorff spaces.

\medskip

Most of the arguments presented in this section can be expressed in a more
algebraic manner, and can be generalized to every semi-simple Lie groups
(see \cite{AnneComp}). The construction we propose here uses only tools of
hyperbolic geometry.

\bigskip

In order to set up the notations, we first make some recalls on the real
hyperbolic space $\dH^n$ and its group of orientation-preserving isometries.
We will be using the hyperboloid model, and we refer to
\cite{BenedetiPetronio}, Chapter A, for a complete overview).

Equip the space $\dR^{n+1}$ with the quadratic form
$$q(x_0,x_1,\ldots,x_n)=x_0x_1+\frac{1}{2}\left(x_2^2+\cdots+x_n^2\right).$$

The set $\left\{\underline{x}\in\dR^{n+1}|q(\underline{x})=-1\right\}$
(where $\underline{x}$ stands for $(x_0,\ldots,x_n)\in\dR^{n+1}$) has two
connected components, and the set $\dH^n$ is defined as
\[ \dH^n=\left\{\underline{x}\in\dR^{n+1}|q(\underline{x})=-1,x_1-x_0>0\right\}.\]
The form
$\langle\underline{x},\underline{y}\rangle=x_0y_1+x_1y_0+x_2y_2+\cdots+x_ny_n$
defines a scalar product on the tangent space at each point of $\dH^n$, hence
a Riemannian metric on $\dH^n$. The geodesics of this space are its
intersections with the (vectorial) planes of $\dR^{n+1}$.

Denote by $\so$ the subgroup of $SL(n+1,\dR)$ consisting of elements preserving
$q$.
Denote by $\sop$ the index $2$ subgroup of $\so$ formed by elements preserving
$\dH^n$ (the other elements exchange the two connected components of
$q^{-1}(-1)$). Then $\sop$ is the group of orientation-preserving
isometries of $\dH^n$.

The image of $\dH^n$ on the hyperplane
$\{x_0-x_1=0\}$,
under the stereographic projection of center
$\left(\frac{1}{\sqrt{2}},-\frac{1}{\sqrt{2}},0,\ldots,0\right)$ is
an open disk of centre $0$ and radius $1$; we denote it
$\dD^n$. This yields the usual compactification of
$\dH^n$, and the projection of any geodesic of
$\dH^n$ on $\dD^n$ is a geodesic in $\dD^n$ for the Poincar\'e metric.

In $\sop$, the stabilizer of the point
$\left(\frac{1}{\sqrt{2}},-\frac{1}{\sqrt{2}},0,\ldots,0\right)\in
\partial\dD^n$
is the subgroup of $\sop$ formed by matrices of the form
$$\left(\begin{array}{ccc}\lambda & 0 & (0) \\ r & \frac{1}{\lambda} & Y\\
Z & (0) & A \end{array}\right),$$
with $\lambda>0$, $A\in SO(n-1)$, $||Y||^2=-\frac{2r}{\lambda}$ and
$A^t Y=-\frac{1}{\lambda}Z$.
The subgroup of elements which also fix the point
$\left(-\frac{1}{\sqrt{2}},\frac{1}{\sqrt{2}},0,\ldots,0\right)
\in\partial\dD^n$
(and in particular, which preserve globally the geodesic
$\dH^n\cap\{x_2=\cdots=x_n=0\}$) is formed by matrices of the form
$$\left(\begin{array}{ccc}\lambda & 0 & (0) \\ 0 & \frac{1}{\lambda} & (0) \\
(0) & (0) & A \end{array}\right),$$
with $\lambda>0$ et $A\in SO(n-1)$.
By analogy with the classical case $n=2$, and with the model of the upper
half plane, we denote by $0$ the point
$\left(\frac{1}{\sqrt{2}},-\frac{1}{\sqrt{2}},0,\ldots,0\right)\in
\partial\dD^n$,
and the point
$\left(-\frac{1}{\sqrt{2}},\frac{1}{\sqrt{2}},0,\ldots,0\right)
\in\partial\dD^n$ is denoted $\infty$.

Elements of $\sop$ having a fixed point in $\dH^n$
are called {\em elliptic}, elements with a unique fixed point in
$\partial\dH^n$ are termed {\em parabolic}, and elements
$\varphi\in\sop$ such that $\inf_{x\in\dH^n}d(x,\varphi(x))>0$ (this lower
bound is then achieved)
are called {\em loxodromic}. Loxodromic isometries fix two points in
$\partial\dH^n$, and are conjugate to a matrix of the form
$$\left(\begin{array}{ccc}\lambda & 0 & (0) \\ 0 & \frac{1}{\lambda} & (0) \\
(0) & (0) & A \end{array}\right).$$
When $A=\Id$, we also say that $\varphi$ is {\em hyperbolic}.

In some texts (which do not include this paper), in the case $n=2$, the
identity is not considered as an elliptic element.

\smallskip

The lower bound
\[d(u)=\inf_{x\in\dH^n}\,d(x,u\cdot x)\]
is achieved if and only if $u$ is non parabolic. In that case
we denote by
$\min(u)$ the set $\{x\in\dH^n\,|\,d(x,u\cdot x)=d(u)\}$,
and if $r>0$,
the set $\{x\in\dH^n\,|\,d(x,u\cdot x)<r+d(u)\}$
is denoted $\min_r(u)$.

The following lemma will be necessary to prove the local compactness of
the quotients $m_\Gamma^o(n)$ and $m_\Gamma^u(n)$.

\begin{lem}\label{argh}
  Let $u\in Isom^+(\dH^n)$ be a non parabolic element. Then $\min_r(u)$ is at
  bounded distance of $\min(u)$. More explicitly, for all $r>0$, there exists
  $k>0$ such that for all $x\in\min_r(u)$ we have $d(x,\min(u))<k$.
\end{lem}

\begin{preuve}
  First suppose that $u$ is loxodromic. In the following figure in the
  hyperbolic plane, the cosine law I gives the formula
  
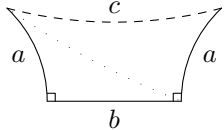
\[\epsfbox{Figures/FiguresPreliminaires.9}\]
  \[\cosh c=1+\cosh^2 a(\cosh b-1).\]
  In particular, for all $x\in\dH^n$ we have
  \[\cosh\left(d(x,u\cdot x)\right)\geq1+\left(\cosh(d(u))-1\right)\cdot
  \cosh^2\left(d(x,\min(u))\right),\]
  hence the lemma, in that case.

  Now suppose that $u$ fixes a unique point $x_0$ in $\dH^n$. Then the angle
  between $[x_0,x]$ and $[x_0,ux]$, as $x$ describes the sphere of radius $1$
  around $x_0$, achieves a minimum
  $\theta_0\neq 0$. For all $x\in\dH^n\smallsetminus\{x_0\}$, denote by
  $\theta(x)$ the angle, at $x_0$, in the triangle $\Delta(x,x_0,u\cdot x)$.
  The cosine law I in this triangle then gives
  \[\begin{array}{rcl}\cosh(d(x,u\cdot x)) & = &
  1+\sinh^2(d(x,\min(u)))(1-\cos\theta) \\ & \geq &
  1+\sinh^2(d(x,\min(u)))(1-\cos\theta_0), \end{array} \]
  which proves the lemma in this case.

  Suppose finally that $u$ is elliptic, and fixes pointwise (at least) a
  geodesic in $\dH^n$. Choose a point $x_0\in\min(u)$.
  For all $x\in\dH^n$, denote $\pi(x)$ the projection of $x$ on $\min(u)$.
  Up to conjugating by an isometry which commutes with $u$, we may suppose that
  $\pi(x)$ is the point $x_0$.
  Then, as $x$
  describes the compact set $\{x\in\dH^n|\pi(x)=x_0,\,d(x,x_0)=1\}$,
  the angle between the segments
  $[x_0,x]$ and $[x_0,\pi(x)]$ achieves a minimum $\theta_0\neq 0$.
  One finds once again, as in the preceding case, for all $x\in\dH^n$,
  the formula
  \[\cosh(d(x,u\cdot x))\geq 1+\sinh^2(d(x,\min(u)))(1-\cos\theta_0)\]
  which completes the proof.
\end{preuve}

Also, we shall use the following classical property:
\begin{prop}\label{bornecompact}
  For every $x\in\dH^n$, and every $d\in\dR$, the set
  \[\{\gamma\in Isom(\dH^n)\,|\,d(x,\gamma x)\leq d\}\]
  is compact.
\end{prop}

Now we can define the space $m_\Gamma^o(n)$:

\begin{deuf}
  We denote by $m_\Gamma^o(n)$ the subspace of $m_\Gamma'(n)$ formed by
  classes of representations which have either $0$, or at least $2$
  global fixed points in $\partial\dH^n$.
\end{deuf}

We have an application
$i\colon m_\Gamma^o(n)\hookrightarrow m_\Gamma'(n)$. We can also define an
application $\pi\colon m_\Gamma'(n)\twoheadrightarrow m_\Gamma^o(n)$ as
follows.
If $c\in m_\Gamma^o(n)$, put $\pi(c)=c$
(in particular, $\pi$ is onto).
If $c=[\rho]\in m_\Gamma'(n)\smallsetminus m_\Gamma^o(n)$, then $\rho$
possesses a unique fixed point $r_1\in\partial\dH^n$. Choose another point
$r_2\in\partial\dH^n\smallsetminus\{r_1\}$ arbitrarily, and denote by
$g_k\in\sop$ the hyperbolic isometry of axis $(r_1,r_2)$, and attractive point
$r_1$, with translation distance $k$.

\begin{lem}\label{lemmg1}
  The sequence $\left(g_k^{-1}\rho g_k\right)_{k\in\dN}$ converges to a
  representation $\rho_\infty\in R_\Gamma(n)$
  such that $[\rho_\infty]\in m_\Gamma^o(n)$,
  and such that $[\rho_\infty]$ does not depend on the choice of $\rho$ in the
  conjugacy class $c$ either of the choice of $r_2$.
\end{lem}

We can therefore set $\pi(c)=[\rho_\infty]$.

\begin{preuve}
  Choose a representant $\rho$ of the conjugacy class $c$ so that
  $\rho$ fixes $0\in\partial\dH^n$ (in other words, conjugate $\rho$ by an
  isometry sending $r_1$ to $0$).
  Take $r_2=\infty$. Then for all
  $\gamma\in\Gamma$, $\rho(\gamma)$ is of the form
  $$\rho(\gamma)=\left(\begin{array}{ccc}\lambda(\gamma) & 0 & (0) \\
  r(\gamma) & \frac{1}{\lambda(\gamma)} & Y(\gamma) \\
  Z(\gamma) & (0) & A(\gamma) \end{array}\right).$$
  In this basis, $g_k$ has the form
  $$g_k=\left(\begin{array}{ccc}t_k & & (0) \\  & \frac{1}{t_k} &  \\
  (0) & & {\rm I}_{n-1} \end{array}\right),$$
  where $t_k\rightarrow 0$ as $k\rightarrow+\infty$. Then $g_k^{-1}\rho g_k$
  converges to the representation $\rho_\infty$ such that for all
  $\gamma\in\Gamma$,
  $$\rho_\infty(\gamma)=\left(\begin{array}{ccc}\lambda(\gamma) & & (0) \\
  & \frac{1}{\lambda(\gamma)} & \\ (0) & & A(\gamma) \end{array}\right),$$
  which indeed fixes the points $0$ and $\infty$ at the boundary.

  Choosing another representant $\rho'$ of $c$ fixing $0$ would simply be
  the same as
  considering a conjugate $\rho'=h^{-1}\rho h$, where $h\in\sop$ has the form
  $$h=\left(\begin{array}{ccc}\lambda & 0 & (0) \\
  r & \frac{1}{\lambda} & Y \\ Z & (0) & A \end{array}\right),$$
  and conjugation by $h$ does not touch the elements $\lambda(\gamma)$ and
  $A(\gamma)$, which determine the representation $\rho_\infty$. Finally, the
  choice of another $r_2$
  amounts to conjugate $\rho$ by an orientation-preserving isometry of $\dH^n$
  fixing $0$; we have just dealt with this case.
\end{preuve}

It follows that we can equip the set $m_\Gamma^o(n)$ both with the induced
topology, and of the quotient topology determined by the surjection
$\pi\colon m_\Gamma'(n)\twoheadrightarrow m_\Gamma^o(n)$.
The object of this section is to prove the following.

\begin{theo}\label{mg}
  The induced topology and the quotient topology coincide on
  $m_\Gamma^o(n)$ (in particular, $\pi$ is continuous). Moreover, the space
  $m_\Gamma^o(n)$ is Hausdorff, and locally compact.
\end{theo}
In particular, this topology is also the quotient topology corresponding to
the surjection $R_\Gamma(n)\twoheadrightarrow m_\Gamma^o(n)$.

It follows that $m_\Gamma^o(n)$ is the biggest Hausdorff quotient of
$m_\gamma'(n)$, in the following sense:
\begin{deuf}
  Let $X$ be a topological space. A quotient space
  $\pi\colon X\rightarrow X_s$ is called {\em the biggest Hausdorff quotient
  of $X$}
  if $X_s$ is Hausdorff and if for every continuous mapping
  $f\colon X\rightarrow Y$ to a Hausdorff space $Y$, there exists a unique
  function
  $\overline{f}\colon X_s\rightarrow Y$ such that $f=\pi\circ\overline{f}$.
\end{deuf}

In particular, such a quotient, if it exists, is unique, up to
{\em canonical} homeomorphism. Clearly, every Hausdorff space is itself
its biggest Hausdorff quotient.

\bigskip

The representations which have no fixed points in $\partial\dH^n$ are called
{\em non parabolic}, and in the sequel of this section, everything that
concerns only non parabolic representations is already written in
\cite{AnneComp}, in the much more general context of actions on
$CAT(0)$ spaces.

We first prove that non parabolic representations form an open subset of
$m_\Gamma'(n)$:

\begin{lem}[compare with \cite{AnneComp}, proposition 2.6]\label{lemmg4}
  Let $\rho\in R_\Gamma(n)$ be non parabolic. Then $[\rho]$
  possesses, in $m_\Gamma'(n)$,
  a neighbourhood consisting of non parabolic representations.
\end{lem}

\begin{preuve}
  The space $R_\Gamma(n)\subset \left(M_{n+1}(\dR)\right)^S$ is a metric space,
  and the map $R_\Gamma(n)\rightarrow m_\Gamma'(n)$ is open, hence the space
  $m_\Gamma'(n)$ is locally numerable (every point possesses a numerable
  fundamental system of neighbourhoods); it follows that we can use sequential
  criteria in this space.

  Consider $\left(\rho_k\right)_k\in m_\Gamma'(n)^\dN$ such that for all
  $k$, $\rho_k$ has a fixed point $r_k\in\partial\dH^n$; and suppose that
  $\rho_k\rightarrow\rho$: let us prove that $\rho$ has a fixed point in
  $\partial\dH^n$. Up to extract a subsequence, $r_k$ converges to a point
  $r\in\partial\dH^n$. Then there exist $h_k\in\sop$ such that
  $h_k(r_k)=r$ and such that $h_k\rightarrow\Id$. Then $h_k\rho_k h_k^{-1}$
  fixes $r$
  globally, and converges to $\rho$ hence $\rho$ fixes $r$ globally.
  It follows that the set of representations which have at least one fixed
  point in $\partial\dH^n$ is a closed subset of $m_\Gamma'(n)$.
\end{preuve}

The first step towards the continuity of $\pi$ is the following.

\begin{lem}\label{lemmg2}
  Let $\left(\rho_k\right)_k$ be a sequence of representations, each having a
  unique fixed point in $\partial\dH^n$, converging to a representation
  $\rho$ which has at least two fixed points in
  $\partial\dH^n$. Then
  $\pi([\rho_k])$ converges to $[\rho]$, in the space $m_\Gamma'$.
\end{lem}

\begin{preuve}
  Let us first prove that up to considering a subsequence, the sequence
  $\pi([\rho_k])$ converges to $[\rho]$.
  The fixed point $r_k$ of $\rho_k$ stays in the compact space
  $\partial\dH^n$, hence there is a subsequence $r_{\varphi(k)}$ of fixed
  points of $\rho_{\varphi(k)}$ which converges to a point
  $r\in\partial\dH^n$, which is therefore a fixed point of $\rho$.
  For all $k$, there exists an orientation-preserving isometry
  $h_{\varphi(k)}$ of $\dH^n$ such that
  $h_{\varphi(k)}(r_{\varphi(k)})=r$, satisfying
  $h_{\varphi(k)}\rightarrow\Id$. Then
  $h_{\varphi(k)}^{-1}\rho_{\varphi(k)} h_{\varphi(k)}$ converges to $\rho$,
  and fixes $r$. In other words,
  we may suppose that $\rho_{\varphi(k)}$ fixes $r$, and up to conjugation we
  can further suppose that $r=0$. Now for all $\gamma\in\Gamma$,
  $\rho_{\varphi(k)}(\gamma)$ and $\rho(\gamma)$ are under the following form:
  \[\rho_{\varphi(k)}(\gamma)=\left(\begin{array}{ccc}
  \lambda_{\varphi(k)}(\gamma) & 0 & (0) \\
  r_{\varphi(k)}(\gamma) & \frac{1}{\lambda_{\varphi(k)}(\gamma)} &
  Y_{\varphi(k)}(\gamma) \\ Z_{\varphi(k)}(\gamma) & (0) &
  A_{\varphi(k)}(\gamma) \end{array}\right),\hspace{0.8cm}
  \rho(\gamma)=\left(\begin{array}{ccc}\lambda(\gamma) & 0 & (0) \\
  0 & \frac{1}{\lambda(\gamma)} & (0) \\
  (0) & (0) & A(\gamma) \end{array}\right),\]
  and by construction, a representant (denote it $\pi\rho_{\varphi(k)}$)
  of the conjugacy class $\pi([\rho_{\varphi(k)}])$ has the form:
  \[\pi\rho_{\varphi(k)}(\gamma)=\left(\begin{array}{ccc}
  \lambda_{\varphi(k)}(\gamma) & 0 & (0) \\
  0 & \frac{1}{\lambda_{\varphi(k)}(\gamma)} & (0) \\
  (0) & (0) & A_{\varphi(k)}(\gamma) \end{array}\right).\]
  Moreover, $\rho_{\varphi(k)}(\gamma)\rightarrow\rho(\gamma)$ for all
  $\gamma\in\Gamma$,
  hence $\lambda_{\varphi(k)}(\gamma)\rightarrow\lambda(\gamma)$ and
  $A_{\varphi(k)}(\gamma)\rightarrow A(\gamma)$, so that
  $\pi\rho_{\varphi(k)}\rightarrow\rho$.

  The proof we have just given works for every subsequence of the sequence
  $(\rho_k)_k$. In particular, every subsequence of
  $\left(\pi([\rho_k])\right)_k$ possesses a subsequence converging to
  $[\rho]$. This implies that $\pi([\rho_k])$ converges to $[\rho]$.
\end{preuve}

We now use an argument of M. Besvina (\cite{Bestvina88}, proposition 1.2).

\begin{lem}[compare with \cite{AnneComp}, proposition 2.5]\label{lemmg3}
  Let $[\rho]\in m_\Gamma^o(n)$. Then the minimum
  \[\min_{x\in\dH^n}\max_{s\in S}d(x,\rho(s)\cdot x)\]
  is achieved.
\end{lem}

\begin{preuve}
  Consider a minimizing sequence $\left(x_n\right)_{k\in\dN}$ for this number.
  If $x_k$ leaves every compact subset of $\dH^n$, then up to considering a
  subsequence, $x_n$ converges to a boundary point
  $r\in\partial\dH^n$. In that case, $r$ is a global fixed point of
  $\rho$, and since $[\rho]\in m_\Gamma^o(n)$, there exists at least one other.
  Hence $\rho$ fixes (globally) a geodesic line of $\dH^n$, and acts by
  translations on this geodesic, and then every point of this geodesic achieves
  the minimum.
  On the other hand, if $(x_k)_k$ is bounded, then it has a subsequence
  converging to some point $x_\infty\in\dH^n$, which realizes this minimum.
\end{preuve}

In the sequel, for every $\rho$ such that $[\rho]\in m_\Gamma^o(n)$, we put
\[d(\rho)=\min_{x\in\dH^n}\max_{s\in S}d(x,\rho(s)\cdot x),\]
\[\min(\rho)=\left\{x\in\dH^n\left|\max_{s\in S} d(x,\rho(s)\cdot x)=d(\rho)
\right.\right\}\]
and
$$\text{\rm min}_\varepsilon(\rho)=\left\{x\in\dH^n\left|\max_{s\in S}
d(x,\rho(s)\cdot x)<d(\rho)+\varepsilon
\right.\right\}.$$

\begin{lem}
  Let $\rho\in m_\Gamma^o(n)$.
  If $\rho$ fixes at least one point at the boundary, then $\min(\rho)$
  is the convex hull of the fixed points of $\rho$ in $\partial\dH^n$; this
  is a totally geodesic subspace of $\dH^n$.
  Otherwise, $\min(\rho)$ is compact. In every case,
  $\min_\varepsilon(\rho)$ is at bounded distance to $\min(\rho)$; that is,
  for all $\varepsilon>0$, there exists $k>0$ such that for all
  $x\in\min_\varepsilon(\rho)$, we have $d(x,\min_\varepsilon(\rho))<k$.
\end{lem}

\begin{preuve}
  Suppose that $\rho$ has no fixed points in $\partial\dH^n$.
  If there existed $\varepsilon>0$ such that $\min_\varepsilon(\rho)$ is not
  bounded, there would exist a sequence $(x_n)_{n\in\dN}$ of elements of
  $\min_\varepsilon(\rho)$, converging to a point
  $x_\infty\in\partial\dH^n$, and then $x_\infty$ would be a fixed point of
  $\rho$. Hence, for all $\varepsilon>0$, $\min_\varepsilon(\rho)$ is bounded;
  moreover $\min(\rho)$ is a closed subset of $\dH^n$, hence compact;
  this finishes the proof in this case.

  Now suppose that $\rho$ has at least two distinct fixed points $x_1,x_2$
  in $\partial\dH^n$. If $d(\rho)\neq 0$, then there exists $s\in S$
  such that $\rho(s)$ is loxodromic, of axis $(x_1,x_2)$, and then
  $\min(\rho)=\min(\rho(s))=(x_1,x_2)$; so $\rho$ does not have any other fixed
  points in $\dH^n$, and for all $\varepsilon>0$, $\min_\varepsilon(\rho)$ lies
  at bounded distance of $\min(\rho)$, by lemma \ref{argh}.
  If $d(\rho)=0$, then $\rho$ fixes pointwise at least the line $(x_1,x_2)$.
  Then
  \[\min(\rho)=\left\{\left.x\in\dH^n\right|\forall s\in S,\ \rho(s)\cdot x
  =x\right\}=\bigcap_{s\in S}\min(\rho(s))\]
  is an intersection of totally geodesic subspaces of $\dH^n$,
  hence it is a totally geodesic subspace of $\dH^n$; and it is the convex hull
  of the fixed points of $\rho$ in $\partial\dH^n$.
  For all $\varepsilon>0$, we then have
  $\min_r(\rho)=\displaystyle{\bigcap_{s\in S}}\min_r(\rho(s))$.
  By lemma \ref{argh}, for all $s\in S$, the set
  $\min_\varepsilon(\rho(s))$ is at bounded distance of the subspace
  $\min(\rho)$. One then checks by induction on $Card(S)$ that
  $\min_\varepsilon(\rho)$ is at bounded distance to $\min(\rho)$.
\end{preuve}

This function
$d\colon m_\Gamma^o(n)\rightarrow\dR_+$
is quite natural:

\begin{prop}\label{dcontinu}
  The function $d'\colon R_\Gamma(n)\rightarrow\dR_+$ defined by
  $d'(\rho)=\displaystyle{\inf_{x\in\dH^n}\max_{s\in S}}\,d(x,\rho(s)x)$
  is continuous.
\end{prop}

Of course, this function $d'$ is constant on conjugation classes,
hence passes to the quotient $m_\Gamma'(n)$; the restriction of this function
to the space $m_\Gamma^o(n)$ is the function $d$.

\begin{preuve}
  Suppose $\rho_k\rightarrow\rho$.
  By construction,
  \[\{\rho|d'(\rho)<a\}=
  \left\{\rho\left|\inf_{x\in\partial\dH^n}\max_{s\in S}
  d(x,\rho(s)x)<a\right.\right\}=
  \left\{\rho\left|\exists x\in\dH^n,\ \forall s\in S,\ 
  d(x,\rho(s)x)<a\right.\right\}
  \]
  is open (that is, the function $d'$ is upper semi-continuous).
  For all $\varepsilon$, we thus have $d'(\rho_k)< d'(\rho)+\varepsilon$,
  for $k$ large enough. In particular $d'$ is continuous at every
  $\rho$ such that $d'(\rho)=0$.

  If $\rho$ has at least one global fixed point $x\in\partial\dH^n$,
  then $d'(\rho)$ is the maximum of translation distances of the
  $\rho(s)$, $s\in S$, such that $\rho(s)$ is loxodromic. Let $s_0$
  be the element of $S$ maximizing this translation distance. Then,
  if $\varepsilon>0$,
  for all $k$ large enough, $\rho_k(s_0)$ is loxodromic, and its translation
  distance is close to that of $\rho(s_0)$, so that
  $d'(\rho_k)>d'(\rho)-\varepsilon$.

  Now suppose that $\rho$ is non parabolic.
  Let $\varepsilon>0$.
  Since the topology of $Isom^+(\dH^n)$ coincides with the compact-open
  topology and since the set
  \[F_{3\varepsilon}=\left\{x\in\dH^n\left|
  \max_{s\in S} d(x,\rho(s)x)\leq 3\varepsilon\right.\right\}\]
  is compact, for all $k$ large enough, for all $x\in F$ and all $s\in S$,
  we thus have
  $|d(x,\rho_k(s)x)-d(x,\rho(s)x)|<\varepsilon$.
  Hence, the convex function
  $x\mapsto\displaystyle{\max_{s\in S}}\,d(x,\rho_k(s)x)$
  has a local minimum
  in the open set $\min_{3\varepsilon}(\rho)$.
  It follows that this minimum is global, whence
  $|d(\rho)-d(\rho_k)|\leq 4\varepsilon$ for $k$ large enough.
\end{preuve}

We turn to a second step towards the proof of theorem \ref{mg}.

\begin{prop}\label{propmg1}
  The space $m_\Gamma^o(n)$, equipped with the induced topology, is Hausdorff.
\end{prop}

\begin{preuve}
  Once again, the space $m_\Gamma'(n)$ being locally second countable,
  we can use sequences in this space.

  Let $[\rho_1]$, $[\rho_2]\in m_\Gamma^o(n)$.
  Suppose that $[\rho_1]$ and $[\rho_2]$ cannot be separated by open sets.
  This means that there exists a sequence
  $\left([\rho_k]\right)_k\in\left(m_\Gamma'(n)\right)^\dN$ converging both to
  $[\rho_1]$ and $[\rho_2]$, in other words, there exist $g_k$,
  $h_k\in\sop$ such that $g_k\rho_k g_k^{-1}\rightarrow \rho_1$ and
  $h_k\rho_k h_k^{-1}\rightarrow\rho_2$. Up to conjugating $\rho_k$ by
  $h_k$, we may suppose that $h_k=1$.

  First suppose that $\rho_2$ is non parabolic.
  Let $x\in\min(\rho_1)$; fix $\varepsilon>0$.
  Then for all $\gamma\in\Gamma$,
  $d(g_k\rho_k(\gamma)g_k^{-1} x,x)=d(\rho_k(\gamma)g_k^{-1}x,g_k^{-1}x)$,
  hence for all $k$ large enough,
  $g_k^{-1}x\in\text{\rm min}_\varepsilon(\rho_2)$,
  since $d(\rho_1)=d(\rho_2)$.
  Since $\text{\rm min}_\varepsilon(\rho_2)$ is bounded, $g_k$ stays in a
  compact set (by proposition \ref{bornecompact})
  and hence has a converging subsequence. Thus, up to taking subsequences, we
  have $g_k\rightarrow g_\infty\in\sop$, and $\rho_1$ and $\rho_2$ are
  conjugate.
  Of course this argument still works after exchanging the roles of
  $\rho_1$ and $\rho_2$.

  Now suppose that $\rho_2$ has at least two distinct fixed points in
  $\partial\dH^n$. Fix $x\in\min(\rho_1)$. Then, as before, for all
  $\varepsilon>0$ and all $k$ large enough we have
  $g_k^{-1}(x)\in\min_\varepsilon(\rho_2)$. If $d(x,g_k^{-1}x)$ is
  bounded then we can finish as in the preceding case. Let us then suppose that
  up to considering a subsequence,
  the sequence $\left(g_k^{-1}x\right)_k$ converges to a point
  $r_1\in\partial\dH^n$. Then $r_1$ is a fixed point of
  $\rho_2$. Choose another fixed point $r_2\in\partial\dH^n$ of $\rho_2$.
  Denote by $r_k$ the second boundary point of the axis $(r_2,g_k^{-1}x)$.
  Then the sequence $(r_k)_k$ in $\partial\dH^n$ converges to $r_1$.
  Then there exists a sequence $(u_k)_{k\in\dN}$ of elements of
  $Isom^+(\dH^n)$,
  converging to $id_{\dH^n}$ and such that $u_k$ fixes $r_2$, and sends
  $r_k$ to $r_1$. Since $r_k$ converges to $r_1$, which is distinct of $r_2$,
  the distance of $x$ to the axis $(r_2,r_k)$ is bounded: there exists $B>0$
  such that for all $k$, $d(x,(r_2,r_k))<B$. Denote by $p_k$ the projection of
  $x$ on the axis $(r_2,r_k)$, and denote by $\varphi_k$ the hyperbolic element
  of axis $(r_1,r_2)$ which sends $u_kg_k^{-1}x$ to $u_kp_k$. Then
  $u_k^{-1}\varphi_k u_k$ sends $g_k^{-1}x$ to $p_k$, so that
  $d\left(x,g_k\cdot\left(u_k^{-1}\varphi_k u_k\right)^{-1}x\right)<B$.
  Hence, by proposition \ref{bornecompact}, up to extract it, the sequence
  $g_k(u_k^{-1}\varphi_k u_k)^{-1}$ converges to an element
  $\varphi\in Isom^+(\dH^n)$. We have
  $g_k\rho_k g_k^{-1}\rightarrow\rho_1$ thus
  $\left(u_k^{-1}\varphi_ku_k\right)\rho_k\left(u_k^{-1}\varphi_ku_k\right)^{-1}
  \rightarrow\varphi^{-1}\rho_1\varphi$. Denote $\rho_k'=u_k\rho_ku_k^{-1}$.
  Since $\rho_k$ converges to $\rho_2$ and since $u_k$ converges to
  $id_{\dH^n}$,
  we have again $\rho_k'\rightarrow\rho_2$. Similarly, we also have
  $\varphi_k\rho_k'\varphi_k^{-1}\rightarrow\varphi^{-1}\rho_1\varphi$.
  Up to conjugate everything, we may suppose that the axis $(r_1,r_2)$ is
  the axis $(0,\infty)$. Since $\rho_2$ and $\varphi_k$ preserve that axis,
  for every $\gamma\in\Gamma$ we can write the elements
  $\rho_2(\gamma)$, $\rho_k'(\gamma)$ and $\varphi_k$ under the form
  \[\rho_2(\gamma)=\left(\begin{array}{ccc}\lambda(\gamma) & 0 & (0) \\
  0 & \frac{1}{\lambda(\gamma)} & (0) \\
  (0) & (0) & A(\gamma) \end{array}\right),\ 
  \rho_k'(\gamma)=\left(\begin{array}{ccc}a_k(\gamma) & b_k(\gamma) &
  X_k(\gamma) \\ c_k(\gamma) & d_k(\gamma) & Y_k(\gamma) \\
  Z_k(\gamma) & W_k(\gamma) & A_k(\gamma) \end{array}\right)\]
  and
  $\varphi_k=\left(\begin{array}{ccc}t_k &  & (0) \\ & \frac{1}{t_k} &  \\
  (0) &  & I_{n-1} \end{array}\right)$,
  where $\displaystyle{\lim_{k\rightarrow+\infty}}\,t_k=+\infty$,
  up to conjugate in order to exchange the points $0$ and $\infty$ in
  $\partial\dH^n$.
  Then
  \[\varphi_k\rho_k'(\gamma)\varphi_k^{-1}=\left(\begin{array}{ccc}a_k(\gamma)
  & \frac{1}{t_k^2}b_k(\gamma) & \frac{1}{t_k}X_k(\gamma) \\
  t_k^2 c_k(\gamma) & d_k(\gamma) & t_k Y_k(\gamma) \\ t_kZ_k(\gamma) &
  \frac{1}{t_k}W_k(\gamma) & A_k(\gamma) \end{array}\right)\]
  so that
  \[\varphi\rho_1(\gamma)\varphi^{-1}=\left(\begin{array}{ccc}\lambda(\gamma) &
  0 & (0) \\
  r(\gamma) & \frac{1}{\lambda(\gamma)} & Y(\gamma) \\
  Z(\gamma) & (0) & A(\gamma) \end{array}\right).\]
  Since $[\rho_1]\in m_\Gamma$, the representation $\varphi\rho_1\varphi^{-1}$
  fixes another point of $\partial\dH^n$ than $0$; denote it $r_3$. There
  exists an isometry $\phi\in Isom^+(\dH^n)$ fixing $0$ and sending $r_3$
  to $\infty$, and now $\phi\varphi\rho_1\varphi^{-1}\phi^{-1}$ and
  $\rho_2$ are conjugate.
\end{preuve}

Now we prove another step towards theorem \ref{mg}:
\begin{prop}
  On the set $m_\Gamma^o(n)$, the induced topology defined by
  $m_\Gamma^o(n)\hookrightarrow m_\Gamma'(n)$
  and the quotient topology defined by
  $\pi\colon m_\Gamma'(n)\twoheadrightarrow m_\Gamma^o(n)$
  coincide (in particular, $\pi$ is continuous).
\end{prop}

\begin{preuve}
  By definition of $\pi$, and by proposition \ref{propmg1},
  for all $x\in m_\Gamma'(n)$, $\pi(x)$ is the unique element of
  $m_\Gamma^o(n)$ such that every neighbourhood of $\pi(x)$ contains $x$. Let
  $U$ be an open subset of $m_\Gamma^o(n)$ for the quotient topology. Then
  $\pi^{-1}(U)$ is open in $m_\Gamma'(n)$, hence the set
  $m_\Gamma^o(n)\cap\pi^{-1}(U)$ is open for the induced topology.
  But $m_\Gamma^o(n)\cap\pi^{-1}(U)=U$, since $\pi$ is the identity on
  $m_\Gamma^o(n)$. Consequently, in order to prove the proposition it suffices
  to prove that $\pi$ is continuous.

  Recall that $m_\Gamma'(n)$ is locally second countable. In order to prove the
  continuity of $\pi\colon m_\Gamma'(n)\rightarrow m_\Gamma^o(n)$, we can use
  a sequential criterium.

  Suppose that $\rho_k,\rho_\infty\in R_\Gamma(n)$ and
  $\rho_k\rightarrow\rho_\infty$. We want to prove that
  $\pi([\rho(k)])$ converges to $\pi([\rho_\infty])$.
  The elements $\rho_k$ have either zero, or at least one fixed point in
  $\partial\dH^n$. Up to consider two distinct subsequences, we may suppose
  that this situation does not depend on $k$.
  If $\rho_k$ is non parabolic, for all $k$,
  then $\pi([\rho_k])=[\rho_k]\rightarrow[\rho_\infty]$, but ever neighbourhood
  of $\pi([\rho_\infty])$ contains $[\rho_\infty]$, hence by definition the
  sequence $\pi([\rho_k])$ also converges to $\pi([\rho_\infty])$.
  Suppose finally that $\rho_k$ has at least one fixed point in
  $\partial\dH^n$, for all $k\in\dN$. By lemma
  \ref{lemmg4}, the representation $\rho_\infty$ has at least one fixed point
  in $\partial\dH^n$, and as before, $[\rho_k]$ also converges
  to $\pi([\rho_\infty])$.
  It then follows from lemma
  \ref{lemmg2} that $\pi([\rho_k])$ converges to $\pi([\rho_\infty])$.
\end{preuve}

Now that we proved that these two topologies coincide, we will equip the set
$m_\Gamma^o(n)$ with this topology in the sequel, without having to precise
which topology we consider. We shall finish the proof of theorem \ref{mg} now:

\begin{cor}\label{dcontinupropre}
  The function $d\colon m_\Gamma^o(n)\rightarrow\dR_+$ is continuous, and
  proper. In particular, the space $m_\Gamma^o(n)$ is locally compact.
\end{cor}

\begin{preuve}
  The map $R_\Gamma(n)\rightarrow m_\Gamma'(n)$ is open, and the continuous
  function $d'\colon R_\Gamma(n)\rightarrow\dR_+$ is constant on the conjugacy
  classes, hence it defines a continuous function of $m_\Gamma'(n)$ in $\dR_+$.
  Hence, $d$ is continuous, by considering the induced topology on
  $m_\Gamma^o(n)$.

  Now, let $A>0$; let us prove that $d^{-1}([0,A])$ is compact. Fix
  $x_0\in\dH^2$, and denote by $R_\Gamma^A\subset R_\Gamma(n)$ the set of
  representations $\rho$ satisfying
  \[ \max_{s\in S}\, d(x_0,\rho(s)x_0)\leq A\, . \]
  By proposition \ref{bornecompact}, this set $R_\Gamma^A$ is compact.
  Denote here by $p$ the map
  $p\colon R_\Gamma(n)\rightarrow m_\Gamma^o(n)$: it is continuous,  and takes
  values in a Hausdorff space, hence $p(R_\Gamma^A)$ is compact; it therefore
  suffices to check that $p(R_\Gamma^A)=d^{-1}([0,A])$.

  Let $[\rho]\in d^{-1}([0,A])$. Then, by lemma \ref{lemmg3},
  there exists a point $x\in\dH^n$ such that
  $\displaystyle{\max_{s\in S}}(x,\rho(s)x)=d(\rho)$, and up to conjugating
  $\rho$ we have $x=x_0$. We then have $\rho\in R_\Gamma^A$, so that
  $[\rho]\in p(R_\Gamma^A)$.
  Now let $\rho\in R_\Gamma^A$. By definition of $d(\rho)$, we have
  $d(\rho)\leq A$. And the map $d$ passes to the quotient
  $p\colon R_\Gamma(n)\rightarrow m_\Gamma^o(n)$, hence
  $d([\rho])=d(\rho)\leq A$, and $p(\rho)\in d^{-1}([0,A])$. This proves that
  $d\colon m_\Gamma^o(n)\rightarrow \dR_+$ is proper.

  Since the space $\dR_+$ is locally compact, so is $m_\Gamma^o(n)$.
\end{preuve}

\bigskip

The group $Isom(\dH^n)$ of isometries of $\dH^n$ which may not be
orientation-preserving, acts on $m_\Gamma^o(n)$ by conjugation.
We denote by $m_\Gamma^u(n)$ the quotient of $m_\Gamma^o(n)$ by this action
($u$ standing for ``unoriented'').

In the case $n=2$, we shall prove here that this quotient identifies to the
space $X_\Gamma(2)$ formed by characters of representations.
For all $\rho\in Hom(\Gamma,\psl)$, denote by
$\chi(\rho)\colon\Gamma\rightarrow\dR_+$ its character. Then:

\begin{prop}\label{mgcaracteres}
  Let $[\rho_1]$, $[\rho_2]\in m_\Gamma^o(2)$. Then
  $\chi(\rho_1)=\chi(\rho_2)$ if and only if there exists an isometry $u$
  of $\dH^2$ such that $\rho_1=u\rho_2 u^{-1}$.
\end{prop}

\begin{preuve}
  Of course, the character is a conjugation invariant in
  $Hom(\Gamma,\psl)$; hence there is only one direction to prove.
  The proof given here inspires of the proof of proposition 1.5.2 of M. Culler
  and P. Shalen \cite{CullerShalen}.

  It can be checked very easily (see \textsl{e.g.} \cite{TheseMaxime},
  proposition 1.1.16) that a representation $\rho$ is elementary
  (\textsl{i.e.}, has a global fixed point in $\dH^2\cup\partial\dH^2$) if and
  only if $\Tr(\rho([\gamma_1,\gamma_2]))$ for all
  $\gamma_1,\gamma_2\in\Gamma$; in particular, elementary representation are
  characterized by their characters.

  Suppose first that $\rho_1$ and $\rho_2$ are not elementary, and suppose
  further that $\rho_1$ is not ``dihedral'' (following \textsl{e.g.}
  \cite{Paulin89}, we say that a representation is {\em dihedral} if it is not
  elementary but fixes globally a pair $\{x,y\}\subset\partial\dH^2$).
  Then there exists $\gamma_0\in\Gamma$ such that $\rho_1(\gamma_0)$ is
  hyperbolic; $\rho_2(\gamma_0)$ is hyperbolic too, and up to conjugating
  $\rho_1$ and $\rho_2$ by an element of $\psl$ we have
  $\rho_1(\gamma_0)=\rho_2(\gamma_0)=\left(\begin{array}{cc}\lambda & 0 \\
  0 & \frac{1}{\lambda}\end{array}\right)$,
  with $\lambda>1$.
  Now let $\gamma_1\in\Gamma$ such that $\rho_1(\gamma_1)$ has no fixed points
  in common with $\rho_1(\gamma_0)$, and is not an elliptic element whose fixed
  point is on the axis of $\rho_1(\gamma)$
  (this is always possible, when $\rho_1$ is not dihedral), and denote
  $\rho_1(\gamma_1)=\left(\begin{array}{cc}a & b \\
  c & d\end{array}\right)$ and
  $\rho_2(\gamma_1)=\left(\begin{array}{cc}a' & b' \\
  c' & d'\end{array}\right)$.
  Then for all $n\in\dZ$,
  $\left|a\lambda^n+\frac{d}{\lambda^n}\right|=\left|a'\lambda^n+\frac{d'}{\lambda^n}\right|$,
  hence, up to change the signs of $a'$, $b'$, $c'$ and $d'$ we have
  $a=a'$ and $d=d'$. Up to conjugating by diagonal matrices we can also suppose
  that $|b|=|b'|=1$, since $\gamma_0$ and $\gamma_1$ do not share any fixed
  points. Up to conjugate by the reflection of axis
  $(0,\infty)$, we can further suppose that $b=b'=1$, and then $c=c'$.

  Now consider any $\gamma\in\Gamma$ and denote
  $\rho_1(\gamma)=\left(\begin{array}{cc}x & y \\ z & t\end{array}\right)$
  and
  $\rho_2(\gamma)=\left(\begin{array}{cc}x' & y' \\ z' & t'\end{array}\right)$.
  Then, for all $\gamma_1$, up to change the signs of $x'$, $y'$,
  $z'$ and $t'$ we have $x=x'$ and $t=t'$. Then, the equality
  $\Tr(\rho_1(\gamma_0^n\gamma_1\gamma))=\Tr(\rho_2(\gamma_0^n\gamma_1\gamma))$
  yields
  \[\left|\lambda^n(ax+bz)+\frac{1}{\lambda^n}(cy+dt)\right|=
  \left|\lambda^n(ax+bz')+\dfrac{1}{\lambda^n}(cy'+dt)\right|\]
  for all $n\in\dN$, where
  $z=z'$ and $y=y'$ or $z=-z'+\frac{2ax}{b}$ and $y=-y'+\frac{2dt}{c}$.
  We can check that this last case implies that
  $\rho_1(\gamma_1)$ is an elliptic element whose fixed point is on the axis
  of $\rho_1(\gamma_0)$, and we have supposed that this is not the case.
  This finishes the proof, in the ``generic'' case.
  Of course, $\rho_1$ and $\rho_2$ play symmetric roles, hence this also covers
  the case where $\rho_2$ is non dihedral.

  If $\rho_1$ and $\rho_2$ are dihedral, then as before,
  we can find $\gamma_0$, $\gamma_1\in\Gamma$
  such that
  $\rho_1(\gamma_0)=\rho_2(\gamma_0)=\left(\begin{array}{cc}\lambda & 0 \\
  0 & \frac{1}{\lambda}\end{array}\right)$
  with $\lambda>1$, and
  $\rho_1(\gamma_1)=\rho_2(\gamma_1)=\left(\begin{array}{cc}0 & 1 \\
  -1 & 0 \end{array}\right)$.
  For all $\gamma\in\Gamma$, $\rho_i(\gamma)$ is then entirely determined by
  the traces of $\rho_i(\gamma_0^n\gamma)$ and
  $\rho_i(\gamma_0^n\gamma_1\gamma)$.

  We still have to deal with the case when $\rho_1$ and $\rho_2$ are
  elementary. Two distinct situations may occur. If $\rho_1$ possesses at least
  one global fixed point in $\partial\dH^2$, then it possesses two, hence all
  the elements $\rho_1(\gamma)$ are hyperbolic and share the same axis;
  $\rho_1$ then identifies to an action of $\Gamma$ on an axis
  $\dR$ by translations, which is determined by its character.
  If $\rho_1$ possesses a (unique) global fixed point in $\dH^2$, then so does
  $\rho_2$. The representations $\rho_1$ and $\rho_2$ factor through morphisms
  $\varphi_i\colon\Gamma\rightarrow\dR$, such that
  $|\cos\varphi_1|=|\cos\varphi_2|$. We can check again that
  $\rho_1$ and $\rho_2$ are conjugate by an isometry.
\end{preuve}

The two quotients $m_\Gamma^u(2)$ and $X_\Gamma(2)$ of the space
$R_\Gamma(2)$ are therefore identical. Hence, $m_\Gamma^o(2)$ is just an
oriented version of the space of characters.

\subsection{Recalls on the Bestvina-Paulin compactification of $m_g^u$}\label{SectionPaulin}

\subsubsection{Equivariant Gromov topology}\label{sectionTopGrom}

We are now going to recall F. Paulin's construction of the compactification
of representation spaces, in order to adapt it to the whole space
$m_\Gamma^u(n)$. W refer to \cite{Paulin04} for an efficient exposition of
this construction.

Here we will be interested in metric spaces $(X,d)$, equipped with actions
of $\Gamma$ by isometries, \textsl{i.e.}, morphisms
$\rho\colon\Gamma\rightarrow Isom(X,d)$; we say that $(\rho,X,d)$ is equivalent
to $(\rho',X',d')$ if there exists a $\Gamma$-equivariant isometry between
$X$ and $X'$.
Let $\mathcal{E}$ be a set of classes of equivariant isometry of actions of
$\Gamma$ on metric spaces.
If $(\rho,X,d)\in\mathcal{E}$ (in order to avoid too heavy notations, we denote
again by $(\rho,X,d)$ its class of equivariant isometry),
if $\varepsilon>0$,
if $\{x_1,\ldots,x_p\}$ is a finite collection in $X$ (denote it $K$) and if
$P$ is a finite subset of $\Gamma$, we define
$U_{K,\varepsilon,P}(\rho,X,d)$
as the subset of $\mathcal{E}$ consisting of those
$(\rho',X',d')$
such that there exist
$\{x_1',\ldots,x_p'\}\in X'$, such that for all $g,h\in P$, and all
$i,j\in\{1,\ldots,p\}$, we have
\[\left|d(\rho(g)\cdot x_i,\rho(h)\cdot x_j)-
d'(\rho'(g)\cdot x_i',\rho'(h)\cdot x_j')\right|
<\varepsilon.\]
In that case we also say that $\{x_1',\ldots,x_p'\}$ is a $P$ equivariant
{\em $\varepsilon$-approximation} of the collection $K$.
The sets $U_{K,\varepsilon,P}(\rho,X,d)$ form a basis of open sets of a
topology, called the {\em equivariant Gromov topology}
(see \cite{Paulin88,Guirardel}).

By definition, every representation $\rho\in R_\Gamma(n)$ defines an action
of $\Gamma$ by isometries on the metric space $(\dH^n,d_{\dH^n})$
(where $d_{\dH^n}$ is the usual distance on $\dH^n$), and every conjugation
by an isometry of $\dH^n$ defines an equivariant isometry. Hence, every
element $[\rho]\in m_\Gamma^u(n)$ defines a unique equivariant isometry class
of actions $(\rho,\dH^n,d_{\dH^n})$; we can therefore consider the set
$m_\Gamma^u(n)$ as a set of (equivariant isometry classes of) actions of
$\Gamma$ on $(\dH^n,d_{\dH^n})$ and we can equip this set with the equivariant
Gromov topology. Then:

\begin{prop}[F. Paulin \cite{Paulin88}, proposition 6.2]\label{propPaulin}
  On the set $m_\Gamma^u(n)$,
  the usual topology and the equivariant Gromov topology coincide.
\end{prop}

Since we will have to adapt this proposition to the oriented case when $n=2$,
we recall here the proof given in \cite{Paulin88}, in that case. Note, anyway,
that the general case $n\geq 2$ is proved similarly.

\begin{preuvede} in the case $n=2$.

  The usual topology is of course finer than the equivariant Gromov topology,
  since the distances considered are continuous for the usual topology on
  $m_\Gamma^u$.

  Conversely, fix $\varepsilon>0$, and let $[\rho_k]$ be a sequence converging
  to $[\rho_\infty]$ for the equivariant Gromov topology.
  Let us prove that up to conjugate these representations
  $\rho_k\rightarrow\rho_\infty$ (since $\dH^2$ is separable, the space
  the space $m_\Gamma^u(n)$ is locally second countable; hence we can indeed
  use sequences in that space).
  Consider three points $x_1,x_2,x_3\in\dH^2$ which form a non degenerate
  triangle. Then, for all $\varepsilon'>0$, and for $k$ large enough, there
  exists a collection $\{x_1^k,x_2^k,x_3^k\}\in\dH^2$ such that for every
  $i,j\in\{1,2,3\}$ and every $s_1,s_2\in S$,
  \[|d(\rho_\infty(s_1)x_i,\rho_\infty(s_2)x_j)-
  d(\rho_k(s_1)x_i^k,\rho_k(s_2)x_j^k)|<\varepsilon.\]
  Now let $S=\{s_1,\ldots,s_n\}$ and
  $y_1=\rho_\infty(s_1)x_1,y_2=\rho_\infty(s_1)x_2,\ldots,
  y_{3n}=\rho_\infty(s_n)x_3$, and similarly define
  $y_1^k,\ldots,y_{3n}^k\in\dH^2$.
  The following fact will enable us to conclude:
  \begin{fact}\label{lemmememetop}
    For all $\varepsilon>0$, there exists $\varepsilon'>0$ such that for all
    $x_1,\ldots,x_n,x_1',\ldots,x_n'$ in $\dH^2$, if
    for all $i,j$, $|d(x_i,x_j)-d(x_i',x_j')|<\varepsilon'$
    then there exists an isometry $\varphi$ of $\dH^2$ such that
    $d(\varphi(x_i),x_i')<\varepsilon$.
  \end{fact}
  This is left as an exercise (see \textsl{e.g.} \cite{TheseMaxime},
  proposition 1.1.8), and follows from the fact that the sine and cosine laws,
  in the hyperbolic plane, can be used to recover continuously a triangle
  from its three lengths.
  In particular, up to conjugate $\rho_k$ by an isometry of $\dH^2$, we thus
  have $d(y_i,y_i')<\varepsilon$, and for every $s\in S$, it follows that
  $d(x_i,\rho_k(s)\cdot\rho_\infty^{-1}(s) x_i)<\varepsilon$, for the three
  non aligned points $x_1,x_2,x_3$, so that $\rho_k\rightarrow\rho_\infty$
  in the usual topology.
\end{preuvede}

For all $\rho\in m_\Gamma^u(n)$, denote
\[\ell(\rho)=\max\left(1,d(\rho)\right)\]
and equip the set $\dH^n$ with the distance
$\dfrac{d_{\dH^n}}{\ell(\rho)}$.
From now on, every element $[\rho]\in m_\Gamma^u(n)$ will be associated to
$(\rho,\dH^n,\frac{d_{\dH^n}}{\ell(\rho)})$ instead of
$(\rho,\dH^n,d_{\dH^n})$ (this is another realization of $m_\Gamma^u(n)$
as a set of (classes of) actions of $\Gamma$ on $\dH^n$, and
the equivariant Gromov topology is still the same, by proposition
\ref{propPaulin} and corollary \ref{dcontinupropre}).

As such, the equivariant Gromov topology does not separate any action on
a space from the restricted actions on invariant subspaces. In particular,
if we consider the whole space $m_\Gamma^u(n)$ as well as actions of $\Gamma$
on $\dR$-trees, including actions on lines, this will yield a non-Hausdorff
space, since some of the elements of $m_\Gamma^u(n)$ have an invariant line.
In order to get rid of this little degeneracy, we are going to modify
slightly the definition of the equivariant Gromov topology, so that elementary
representations in $m_\Gamma^u$ will be separated from the corresponding
actions on lines, when considered as actions on $\dR$-trees.
If $\mathcal{E}$ is a set of (classes of) actions of $\Gamma$ by isometries
on spaces which are {\em hyperbolic in the sense of Gromov}, put
$U_{K,\varepsilon,P}'(\rho,X,d)$
to be the subset of $\mathcal{E}$ consisting of those
$(\rho',X',d')$
such that there exist
$\{x_1',\ldots,x_p'\}\in X'$, such that for all $g,h\in P$, and all
$i,j\in\{1,\ldots,p\}$, we have
\[\left|d(\rho(g)\cdot x_i,\rho(h)\cdot x_j)-
d'(\rho'(g)\cdot x_i',\rho'(h)\cdot x_j')\right|<\varepsilon
\text{ and }|\delta(X)-\delta(X')|<\varepsilon.\]

As we shall see very soon (proposition \ref{deltaproche}), this extra
condition changes the equivariant Gromov topology only at the neighbourhood
of elementary representations.
It is comparable to the fact of adding $2$ to characters, as it is done by
J. Morgan and P. Shalen in
\cite{MorganShalen}, in order not to bother with the neighbourhood of the
trivial representation.

\subsubsection{The space $\overline{{m}_\Gamma^u(n)}$}\label{sectionCompactif}

Let us begin this section by a few recalls. Let $X$ be a topological space.
A {\em compactification of $X$} is a couple
$(\overline{X},i)$ such that $\overline{X}$ is a compact Hausdorff space,
$i\colon X\hookrightarrow\overline{X}$ is a homeomorphism on its image, and
such that $i(X)$ is open and dense in $\overline{X}$.

Note that we request $\overline{X}$ to be Hausdorff. In particular, a space
$X$ needs to be locally compact in order to admit a compactification.
If $X$ is a locally compact space and $(\overline{X},i)$ is a
compactification of $X$, the compact set $\overline{X}\smallsetminus X$ is
called the {\em boundary} of $X$, and is denoted $\partial\overline{X}$, or
$\partial X$ if no confusion is possible between different compactifications
of $X$. The points in the boundary are called
the {\em ideal points} of the compactification.

Note that, if $X$ is locally compact and if $(\overline{X},i)$
is a compactification of $X$, then the open subsets of $\overline{X}$
containing $\partial\overline{X}$ are precisely the complements, in
$\overline{X}$, of the compact subsets of $X$ (indeed, if
$K\subset X$ is compact then $i(K)$ is compact, hence
$\overline{X}\smallsetminus i(K)$ is open and contains
$\partial\overline{X}$, and if $U$ is an open subset of $\overline{X}$
containing
$\partial\overline{X}$ then $\overline{X}\smallsetminus U$ is compact,
contained in $i(X)$, so that it is the image of a compact set, since
$i$ is a homeomorphism on its image).

Finally, if $X$ is a locally compact space, if $Y$ is Hausdorff and
$f\colon X\rightarrow Y$ is continuous and has a relatively compact image,
then we can define a compactification of $X$ as follows (see
\cite{MorganShalen}, p. 415). Denote by $\widehat{X}=X\cup\{\infty\}$ the
Alexandrov compactification of $X$ (that is, the one-point compactification),
we define $i\colon X\rightarrow\widehat{X}\times Y$ by $i(x)=(x,f(x))$.
Denote $\overline{X}$ the closure of $i(X)$ in $\widehat{X}\times Y$.
Then $(\overline{X},i)$ is a compactification of $X$; we say that it is
the {\em compactification defined by $f$}.

\bigskip

We say that an action of $\Gamma$ on an $\dR$-tree $T$ is
{\em minimal} if $T$ has no proper invariant sub-tree.
The equivalence classes of $\dR$-trees equipped with minimal actions of
$\Gamma$ by isometries, up to equivariant isometry, form a set, and we denote
$\mathcal{T}'(\Gamma)$ the subset formed by trees not reduced to a point.
In order to exhibit this set, one can prove that the $\dR$-tree $T$ and the
action of $\Gamma$ are entirely determined by the set
$\{d(p,\gamma\cdot p)|\gamma\in\Gamma\}$ (see \cite{AlperinMoss,MorganShalen}).

\smallskip

We have modified the definition of the equivariant Gromov topology so that
we would be considering Hausdorff spaces, and for this reason we are also
going to restrict the set of $\dR$-trees we consider.
If $(\rho,T)$ possesses an end which is globally fixed by $\rho$
(we then say that this action is {\em reducible}, see \textsl{e.g.}
\cite{Paulin89}),
then $(\rho,T,d)$ is not separated, in the equivariant Gromov topology,
from the action on a line which has the same translation lengths.
Therefore, we shall restrict ourselves to the subset
$\mathcal{T}(\Gamma)\subset\mathcal{T}'(\Gamma)$ consisting of actions on lines
and actions on $\dR$ trees without fixed ends, such that
$\displaystyle{\min_{x_0\in T}
\max_{\gamma\in S}}\, d(x_0,\gamma\cdot x_0)=1$.

In the sequel, we equip the set $m_\Gamma^u(n)\cup\mathcal{T}(\Gamma)$
with the modified equivariant Gromov topology, and we denote by
$\overline{m_\Gamma^u(n)}$ the closure of the set
$m_\Gamma^u(n)$ in the space $m_\Gamma^u(n)\cup\mathcal{T}(\Gamma)$.

The following proposition
is simply a detailed version of an argument of
F. Paulin, in \cite{Paulin88}, saying that $\dR$-trees and hyperbolic
structures are not undistinguishable by the equivariant Gromov topology, so
that the modified equivariant Gromov topology changes it only at the boundary
of degenerate representations.

We say that a finite collection $K$ is {\em big} if it contains four points
$A,B_1,B_2,B_3$ such that $d(A,B_i)=1$, $d(B_1,B_3)=2$ and
$d(B_1,B_2)=d(B_2,B_3)$. In an $\dR$-tree, this means that their convex hull
is a tripod of centre $A$, and this implies that
$d(B_1,B_2)=2$. In $\dH^n$, this means that the segments
$[B_1,B_3]$ and $[B_2,A]$ meet orthogonally at $A$.

\begin{prop}\label{deltaproche}
  Let $(\rho,X,d)\in m_\Gamma^u(n)\cup\mathcal{T}(\Gamma)$
  and let $K\subset X$ be a big finite collection.
  Let $\varepsilon>0$. Then for all $\varepsilon'>0$
  small enough, for all
  $(\rho',X',d')\in U_{K,\varepsilon',\{1\}}(\rho,X,d)$,
  we have $|\delta(X')-\delta(X)|<\varepsilon$.
\end{prop}

Here, by $\varepsilon'$ small enough, we mean:
$\varepsilon'<\upsilon(\delta(X),\varepsilon)$, where
$\upsilon\colon\dR_+\times\dR_+^*\rightarrow\dR_+^*$ is some (universal)
function.

\begin{preuve}
  Let $\varepsilon'>0$, and let $K'=\{A',B_1',B_2',B_3'\}\subset X'$ be an
  $\varepsilon'$-approximation of $K$.
  Put $B_1''=B_1'$. We have
  $|d(B_1'',B_3')-2|<\varepsilon'$ so we can choose a point $B_3''\in X'$ such
  that $d(B_3',B_3'')<\varepsilon'$ and $d(B_1'',B_3'')=2$. Denote by $r(t)$
  the geodesic segment, in $X'$, such that $r(0)=B_1''$ and $r(2)=B_3''$, and
  put $A''=r(1)$. Then the $CAT(0)$ property on $X'$ implies that
  $d(A',A'')<\sqrt{2\varepsilon'+\varepsilon'^2}$. Finally, we have
  $|d(A',B_2')-1|<\varepsilon'$ so there exists $B_2''\in X$ such that
  $d(B_2',B_2'')<\varepsilon'+\sqrt{2\varepsilon'+\varepsilon'^2}$ and
  $d(A'',B_2'')=1$. We also have the inequalities
  $|d(B_1'',B_2'')-1|<2\varepsilon'+\sqrt{2\varepsilon'+\varepsilon'^2}$,
  and
  $|d(B_3'',B_2'')-d(B_1'',B_2'')|<
  5\varepsilon'+2\sqrt{2\varepsilon'+\varepsilon'^2}$.

  First suppose that $\delta(X)\neq 0$, and denote
  $x=\dfrac{\delta_{\dH^2}}{\delta(X)}$.
  The cosine law I, in $(\dH^2,d_{\dH^2})$, implies that
  $\cosh(x d(B_1,B_2))=\cosh^2(x)$.
  For all $\varepsilon'$ small enough, we have
  $d(A'',B_i'')=1$, $d(B_1'',B_3'')=2$, $d(B_1'',B_2'')<2$ and
  $d(B_3'',B_2'')<2$, so that $X'$ cannot be an $\dR$-tree.
  Put also
  $x'=\dfrac{\delta_{\dH^2}}{\delta(X')}$.
  Then the cosine law I in $X'$ implies that
  $\cosh(x'd(B_1'',B_2''))+\cosh(x'd(B_2'',B_3''))=2\cosh^2(x')$.
  It follows from the study of the function
  $F\colon [0,2]\times\dR_+\rightarrow\dR$
  defined by $F(b,x)=\cosh^2(x)-\cosh(xb)$, that by taking
  $\varepsilon'$ small enough we can force $x$ and $x'$ to be arbitrarily
  close.
  
  Now suppose that $\delta(X)=0$. If $\delta(X')=0$ then there is nothing to
  do. Otherwise, we have again
  $\cosh(x'd(B_1'',B_2''))+\cosh(x'd(B_2'',B_3''))=2\cosh^2(x')$,
  where, for $\varepsilon'$ small enough, the distances
  $d(B_1'',B_2'')$ and $d(B_2'',B_3'')$ can be taken arbitrarily close to $2$,
  which implies that $x'$ can be forced to be arbitrarily large.
\end{preuve}

In particular, the modified equivariant Gromov topology and the equivariant
Gromov topology coincide in $m_\Gamma^u(n)$.

\smallskip

Now, every arguments of M. Bestvina and F. Paulin
(see \cite{Bestvina88,Paulin88,Paulin04}) work, and we have the following.

\begin{theo}[M. Bestvina, F. Paulin]\label{theoBestPau}
  The space $\overline{m_\Gamma^u(n)}$, equipped with the function
  $m_\Gamma^u(n)\hookrightarrow\overline{m_\Gamma^u(n)}$, is a natural
  compactification of $m_\Gamma^u(n)$.
\end{theo}

By ``natural'', we mean that the action of $Out(\Gamma)$ on $m_\Gamma^u(n)$
extends continuously to an action of $Out(\Gamma)$
on $\overline{m_\Gamma^u(n)}$.

We refer to
\cite{Paulin04,Paulin89} and \cite{KapovichLeeb} for a complete proof of this
result.

\subsubsection{Other compactifications}

We now give a (very) short recall on the compactification
of $X_{\Gamma,SL(2,\dR)}$ by J. Morgan and P. Shalen.
The numerable collection
$(f_\gamma)_{\gamma\in\Gamma}$,
$f_\gamma\colon\chi\mapsto \chi(\gamma)$ generates the coordinate ring of
$X_{\Gamma,SL(2,\dR)}$. Denote $P\dR^\Gamma$ the quotient of
$[0,+\infty)^\Gamma\smallsetminus\{0\}$ by positive multiplication,
and let $\theta\colon X_{\Gamma,SL(2,\dR)}\rightarrow P\dR^\Gamma$
defined by
$\theta(x)=\left[\log(|f_\gamma(x)|+2)\right]_{\gamma\in\Gamma}.$
J. Morgan and P. Shalen proved (see \cite{MorganShalen}, proposition I.3.1)
that the image of $X_{\Gamma,SL(2,\dR)}$ under $\theta$ is relatively compact,
so that $\theta$ defines a compactification of $X_{\Gamma,SL(2,\dR)}$.

We now restrict to the group $\Gamma=\pi_1\Sigma_g$ with $g\geq 2$,
and we denote by $m_g^u$ the space $m_\Gamma^u(2)$. Then the absolute value
of the Euler class is defined on $m_g^u$, and we denote
$m_{g,\text{even}}^u$ the subspace of $m_g^u$ consisting of representations
of even Euler class
(recall, indeed, that a representation
$\rho\colon\pi_1\Sigma_g\rightarrow\psl$ lifts to $SL(2,\dR)$ if and only if
its Euler class is even; we will see that in section \ref{Euleretpsl}).
Then the map $\theta$ factors through
$\theta'\colon m_{g,\text{even}}^u\rightarrow P\dR^\Gamma$, which defines
a compactification $\overline{m_{g,\text{even}}^u}^{MS}$ of
$m_{g,\text{even}}^u$. Similarly, the functions
$m_{g,\text{even}}^u\hookrightarrow m_g^u\hookrightarrow\overline{m_g^u}$
define a compactification $\overline{m_{g,\text{even}}^u}$ of that space.
The ideal points of the compactifications
$\overline{m_{g,\text{even}}^u}$ and
$\overline{m_{g,\text{even}}^u}^{MS}$ are actions of $\pi_1\Sigma_g$
on $\dR$-trees. These actions on $\dR$-trees are irreducible, \textsl{i.e.},
without global fixed points, or are actions on a line. The topology on
$\partial\overline{m_{g,\text{even}}^u}^{MS}$ is the
{\em axis topology}; it is the coarsest topology such that the functions
$\ell_T(\gamma)=\displaystyle{\inf_{x\in T}}\, d(x,\gamma x)$ are continuous.
By the main theorem of \cite{Paulin89}, the spaces
$\partial\overline{m_{g,\text{even}}^u}$ and
$\partial\overline{m_{g,\text{even}}^u}^{MS}$ are homeomorphic, and hence, as
F. Paulin explains it in \cite{Paulin88}, the spaces
$\overline{m_{g,\text{even}}^u}$ and
$\overline{m_{g,\text{even}}^u}^{MS}$ are homeomorphic (we refer to
\cite{MorganShalen,Paulin88,Paulin89} for details).
In particular, the corollary \ref{introcarac} concerns simply the space
$\overline{m_{g,{\text{\rm even}}}}$ and it is under that form that we shall
prove it in section \ref{ChapitreCompact}.

\subsection{Euler class}\label{SectionEulerClass}


We shall now introduce the Euler class, for actions of groups on cyclically
ordered sets. For the sake of completeness we shall give a complete
construction, but equivalent definitions can be found in
\cite{Ghys87,Thurston3}. We refer to \cite{Calegari} for a detailed overview.

\subsubsection{Cyclically ordered sets}

\begin{deuf}\label{defordrecyclique}
  Let $X$ be a set. A {\em (total) cyclic order} on $X$ is a function
  $o:X^3\rightarrow \{-1,0,1\}$ such that:
  \begin{itemize}
    \item[(i)] $o(x,y,z)=0$ if and only if $\text{card}\{x,y,z\}\leq 2$;
    \item[(ii)] For all $x$, $y$ and $z$, $o(x,y,z)=o(y,z,x)=-o(x,z,y)$;
    \item[(iii)] For all $x$, $y$, $z$ and $t$,
    if $o(x,y,z)=1$ and $o(x,z,t)=1$ then $o(x,y,t)=1$.
  \end{itemize}
\end{deuf}

\begin{ps}
  If $o(x,y,z)=1$ and $o(x,z,t)=1$ then we also have $o(x,z,t)=1$
  and $o(y,z,t)=1$. Indeed, $o(z,x,y)=o(z,t,x)=1$ so, by condition
  (iii) of the definition, $o(z,t,y)=1$, that is,
  $o(y,z,t)=1$. Similarly, $o(x,y,t)=1$ so $o(t,x,y)=1$, which, together with
  $o(t,y,z)=1$, yields $o(t,x,z)=1$, \textsl{i.e.} $o(x,z,t)=1$. In other
  words, the transitivity relation (iii) implies all the other
  ``natural'' transitivity relations. In particular, for instance, on a set
  of $4$ elements, there are
  as many total cyclic orders as injections of that set in the oriented circle,
  up to homeomorphism, that is, $6$.
\end{ps}

\begin{ps}
  In all this text, we use only
  {\em total} cyclic orders ({\em every} triple defines an order).
  Consequently, we will sometimes forget the word ``total'' when we refer to
  a cyclic order.
\end{ps}

Now fix a set $X$ equipped with a cyclic order $o$ and a base point $x_0\in X$.

\begin{deuf}
  We set $y<_{x_0}z$ if $o(x_0,y,z)=1$, and $y\leq_{x_0}z$ if
  $y<_{x_0}z$ or $y=z$.
\end{deuf}

The following proposition follows directly from the properties of $o$.

\begin{prop}
  The relation $\leq_{x_0}$ is a total order on
  $X\smallsetminus\{x_0\}$.
\end{prop}

\begin{ps}
  Reciprocally, if $\leq$ is a total order on
  $X\smallsetminus\{x_0\}$, then there exists a unique cyclic total order
  on $X$ which verifies $o(x,y,z)=1$ for all $x,y,z\neq x_0$ such that
  $x<y<z$, and verifying $o(x_0,y,z)=1$ as soon as $y<z$ (by enumerating the
  different possible cases we can check that the axioms of definition
  \ref{defordrecyclique} are satisfied).
  For all $x_0\in X$, these two constructions realize a bijection, and its
  inverse, between the set of total cyclic orders on $X$ and the set of total
  orders on $X\smallsetminus\{x_0\}$.
  In particular, on any set $X$ there exists at least one total cyclic order,
  by the axiom of choice.
\end{ps}

\begin{deuf}
  On the set $\dZ\times X$ we put:
  \begin{itemize}
    \item $(m,x)<_{x_0}(n,y)$ when $m<n$, for any $x,y\in X$,
    \item $(k,y)<_{x_0}(k,z)$ when $y<_{x_0}z$,
    \item $(k,x_0)<_{x_0}(k,y)$ for all $y\in X\smallsetminus\{x_0\}$,
  \end{itemize}
  and we put $(m,y)\leq_{x_0}(n,z)$ if $(m,y)<_{x_0}(n,z)$ or $(m,y)=(n,z)$.

  In particular, when restricted to the set
  $\dZ\times(X\smallsetminus\{x_0\})$ it is the lexicographic order.
\end{deuf}

We check easily the following.

\begin{prop}
  The relation $\leq_{x_0}$ is a total order on the set
  $\dZ\times X$.
\end{prop}

\begin{eg}
  If $X=\dS^1$ and $x_0\in\dS^1$, then $X\smallsetminus\{x_0\}$ is an
  interval,
  \[\epsfbox{Figures/FiguresOrdre.1}\]
  and $\dZ\times X$ identifies naturally to $\dR$, and this identification
  depends canonically on $x_0$.
\end{eg}

Another example is the deck of cards. It is equipped with a cyclic order, which
is preserved as we make a ``cut''. The choice of the
``cut'' consists in choosing a card $x_0$, and determines a total order in the
deck. Here our definition of the order on $\dZ\times X$ consists in choosing a
cut, and then to put $\dZ$ copies of the deck the ones above the others.

\medskip

\subsubsection{Applications and lifts}

We consider here a set $X$ equipped with a total cyclic order $o$.
We shall see that the bijections of $X$ preserving $o$ can be lifted to
$\dZ\times X$ in some natural way (as soon as we take a base point $x_0$),
similarly as the orientation-preserving homeomorphisms of $\dS^1$
can be lifted naturally to its {\em cover} $\dR$.

\begin{deuf}
  Let $f:X\rightarrow X$ be any map.
  We call an {\em arbitrarily lift} of $f$ any function
  $\tilde{f}$ such that the following diagram commutes:
  \[\xymatrix{\dZ\times X\ar[r]^{\tilde{f}}\ar[d]^{pr_2} &
  \dZ\times X\ar[d]^{pr_2} \\ X\ar[r]^f & X\ .}\]
  In that case we say that $\tilde{f}$ {\em projects on} $f$.
\end{deuf}

We define an application $h\colon \dZ\times X\rightarrow \dZ\times X$
by $h(n,x)=(n+1,x)$.

\begin{prop}\label{uniqueh}
  Let $f:\dZ\times X\rightarrow \dZ\times X$ be a bijection preserving the
  order $\leq_{x_0}$ and which projects to $id_X$. Then there exists an integer
  $n$ such that $f=h^n$.
\end{prop}

\begin{preuve}
  Let $f$ be such a function; put $n=pr_1(f(0,x_0))$. Then the function
  $pr_1\circ f(\cdot,x_0):\dZ\rightarrow\dZ$ is a bijection preserving the
  order, hence for all $k\in\dZ$, $f(k,x_0)=(n+k,x_0)$. Now, for all
  $y\in X\smallsetminus\{x_0\}$, we have
  $f(k,x_0)\leq_{x_0}f(k,y)\leq_{x_0}f(k+1,x_0)$,
  that is,
  $(n+k,x_0)\leq_{x_0}f(k,y)\leq_{x_0}(n+k+1,x_0)$, hence $f(k,y)=(n+k,y)$,
  since $f$ projects on $id_X$.
\end{preuve}

\begin{prop}\label{transportordre}
  For all $x_0,x_1\in X$, there exists a unique bijection
  $F_{x_0x_1}$ of $\dZ\times X$
  which projects on $id_X$ such that for all $a,b\in \dZ\times X$,
  $a<_{x_0}b\Leftrightarrow F_{x_0x_1}(a)<_{x_1}F_{x_0x_1}(b)$, and such that
  $F_{x_0x_1}(0,x_0)=(0,x_0)$.
\end{prop}

\begin{preuve}
  If $x_0=x_1$, then by the preceding proposition, the unique possible function
  is $F_{x_0x_1}=id_{\dZ\times X}$. Now suppose that $x_0\neq x_1$.
  The application $F_{x_0x_1}$ must project on
  $id_X$, so that we can only change indices, in the way suggested by the
  following picture.
  
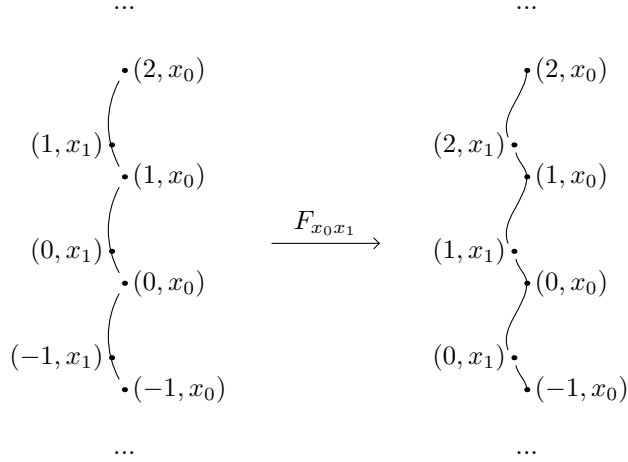
\[\epsfbox{Figures/FiguresOrdre.2}\]
  We put $F_{x_0x_1}(k,x_0)=(k,x_0)$ and $F_{x_0x_1}(k,x_1)=(k+1,x_1)$.
  For all $y\in X\smallsetminus\{x_0,x_1\}$, if $o(x_0,x_1,y)=-1$ we put
  $F_{x_0x_1}(k,y)=(k,y)$ otherwise we put $F_{x_0x_1}(k,y)=(k+1,y)$:
  we then have $a<_{x_0}b\Leftrightarrow F_{x_0x_1}(a)<_{x_1}F_{x_0x_1}(b)$.
  Now, if $f$ is another bijection satisfying the same conditions, then
  $F_{x_0x_1}\circ f^{-1}$ is a bijection of $\dZ\times X$ which preserves
  the order $\leq_{x_0}$, which projects on $id_X$ and which fixes
  $(0,x_0)$ hence, by proposition \ref{uniqueh}, it is the identity.
\end{preuve}

Note that for all distinct points $x_0,x_1\in X$ we have
$F_{x_0x_1}\circ F_{x_1x_0}=h$,
contrarily to what our notation may suggest.

\begin{prop}
  Let $X$ be a set equipped with a (total) cyclic order $o$, let
  $x_0\in X$, and let $f\colon X\rightarrow X$ be a bijection preserving $o$.
  Then $f$ admits at least a lift $\tilde{f}$ preserving the order
  $\leq_{x_0}$. Moreover, if $\tilde{f}$ and $\tilde{f'}$ are two such lifts
  then there exists $n\in\dZ$ such that $\tilde{f'}=h^n\cdot\tilde{f}$.
\end{prop}


\begin{preuve}
  Let $f\colon X\rightarrow X$ preserving $o$. We check easily that the map
  $F:\dZ\times X\rightarrow\dZ\times X$ defined by $F(n,x)=(n,f(x))$ verifies:
  $\forall a,b\in \dZ\times X$,
  $a\leq_{x_0}b\Leftrightarrow F(a)\leq_{f(x_0)}F(b)$.
  We therefore put $\tilde{f}=F_{f(x_0)x_0}\circ f$. Then $\tilde{f}$ is indeed
  a lift of $f$, preserving $\leq_{x_0}$.

  Now, if $\tilde{f'}$ is another lift of $f$ preserving $\leq_{x_0}$,
  then $\tilde{f'}\circ\tilde{f}^{-1}$ is a lift of $id_X$ preserving
  $\leq_{x_0}$, hence, by proposition \ref{uniqueh},
  $\tilde{f'}\tilde{f}^{-1}=h^n$, for some $n\in\dZ$.
\end{preuve}

Now denote by $Ord(X,o)$ the group of bijections of $X$ which preserve
$o$, and by $\widetilde{Ord(X,o,x_0)}$ the group formed by lifts of elements of
$Ord(X,o)$ à $\dZ\times X$ which preserve the order $\leq_{x_0}$.
By proposition \ref{transportordre}, the conjugation by
$F_{x_0x_1}$ realizes a canonical isomorphism between the groups
$\widetilde{Ord(X,o,x_0)}$ and $\widetilde{Ord(X,o,x_1)}$, for all
$x_0,x_1\in X$. This group, considered up to isomorphism, is denoted
$\widetilde{Ord(X,o)}$, and if there is no confusion possible about $o$, these
two groups are denoted $Ord(X)$ and $\widetilde{Ord(X)}$.


\subsubsection{Euler class}

Let $(X,o)$ be a cyclically ordered set, let $\Sigma$ be a (connected) surface,
and let $\rho\colon\pi_1\Sigma\rightarrow Ord(X)$ be a representation: we shall
define an element $e(\rho)\in H^2(\Sigma,\dZ)$, depending (only) on
$\rho$, and which we call the Euler class of the representation.

Since we are considering the representations of $\pi_1\Sigma$ in a non abelian
group, we first need to choose a base point $*\in\Sigma$.
Choose also a base point $x_0\in X$. Let then $C=C_0\cup C_1\cup C_2$;
$C_i=\left\{\sigma_\alpha^i\right\}_\alpha$, be a cellulation of
$\Sigma$, where every cell is equipped with an orientation, and such that
$*\in C_0$.
For each loop $\gamma$ (based at $*$)
in the $1$-skeleton, we can write $\gamma$ as a word
$(\sigma_1^1)^{s_1}\cdots(\sigma_k^1)^{s_k}$, $s_j=\pm 1$,
in the elements of $C_1$. We shall say that a function
$f\colon C_1\rightarrow\widetilde{Ord(X)}$ {\em projects on $\rho$}
if for every loop $\gamma$ based at $*$
in the $1$-skeleton, the element
$\left(f(\sigma_1^1)\right)^{s_1}\circ\cdots\circ
\left(f(\sigma_k^1)\right)^{s_k}$
(keeping the same notations for $\gamma$),
which we call $f(\gamma)$, is a lift of $\rho(\gamma)$.
Since the $1$-skeleton is homotopic to a bouquet of circles, there exists
at least one such function $f$
(indeed, we can choose a maximal tree $T$ containing $*$ in the
$1$-skeleton, put $f(\sigma^1)=id_X$ for every edge $\sigma^1$ of $T$,
and for each edge $\sigma^1$ which is not in $T$, $\sigma^1$ (equipped with
its orientation) defines, jointly with $T$, an element
$\gamma\in\pi_1\Sigma$, and we put $f(\sigma^1)=\rho(\gamma)$).
Now, the boundary of a $2$-cell $\sigma_\alpha^2$ is a loop
$\gamma(\sigma_\alpha^2)$ (not based at $*$) in the $1$-skeleton,
homotopically trivial in $\Sigma$.
We can again write it as a word
$(\sigma_1^1)^{s_1}\cdots(\sigma_k^1)^{s_k}$, $s_j=\pm 1$,
well defined up to cyclic permutation, and, by proposition \ref{uniqueh} and
by the fact that $h$ is central in $\widetilde{Ord(X)}$, there exists
$n(\sigma_\alpha^2)\in\dZ$ such that
$f(\gamma(\sigma_\alpha^2))=h^{n(\sigma_\alpha^2)}$. Finally,
if $c=\Sigma\lambda_i\sigma_i^2$ is a $2$-cycle, we put
$e(\rho)\cdot c=\Sigma\lambda_i n(\sigma_i^2)$.

\begin{theo}\label{eulerbiendefinie}
  For every $2$-cycle $c$, the integer $e(\rho)\cdot c$ depends only on $(X,o)$
  and on $\rho$; this defines an element $e(\rho)\in H^2(\Sigma,\dZ)$, called
  the {\em Euler class of the representation $\rho$}.
\end{theo}

If the surface $\Sigma$ is oriented, then the evaluation of $e(\rho)$ on the
fundamental class is an integer, which we still call (abusively) the Euler
class of the representation, $e(\rho)\in\dZ$.

\begin{preuve}
  Step 1: we first prove that given base points $*\in\Sigma$, $x_0\in X$ and a
  cellulation, $e(\rho)$ does not depend on the choice of the function $f$.
  Let $f_1$ and $f_2$ two such functions. Choose a maximal tree $T$ in the
  $1$-skeleton, containing $*$.
  For every cell $\sigma^1\in T$, put $f_1'(\sigma^1)=f_2'(\sigma^1)=1$.
  Every cell $\sigma^1$ which is not in $T$ determines, jointly with $T$, a
  loop $\gamma(\sigma^1)$ (since $\sigma^1$ is equipped with an orientation)
  in the $1$-skeleton, based at $*$. For $i=1,2$ we put
  $f_i'(\sigma^1)=f_i(\gamma(\sigma^1))$. Since these elements
  $\sigma^1\not\in T$ generate the fundamental group of the $1$-skeleton,
  for every loop $\gamma$ in the $1$-skeleton we thus have
  $f_i(\gamma)=f_i'(\gamma)$, in particular
  $f_i'$ projects on $\rho$ and defines the same Euler class as $f_i$.
  Now, for all $1$-cell $\sigma^1$, there exists an integer $n(\sigma^1)\in\dZ$
  such that $f_1'(\sigma^1)\circ f_2'(\sigma^1)^{-1}=h^{n(\sigma^1)}$: this is
  $0$ if $\sigma^1\in T$, and otherwise this follows from proposition
  \ref{uniqueh}.
  Then for all $2$-cell $\sigma_\alpha^2$, we have
  $f_1(\gamma(\sigma_\alpha^2))-f_2(\gamma(\sigma_\alpha^2))=
  \displaystyle{\sum_{\sigma^1\in\partial\sigma^2}n(\sigma^1)}$
  (indeed, if $\gamma(\sigma_\alpha^2)$ is given by the cyclic word
  $(\sigma_1^1)^{s_1}\cdots(\sigma_k^1)^{s_k}$ then
  \begin{eqnarray*}
  f_1(\gamma(\sigma_\alpha^2)) & = & f_1(\sigma_1^1)^{s_1}\circ\cdots\circ
  f_1(\sigma_k^1)^{s_k}=
  \left(f_2(\sigma_1^1)^{s_1}\circ h^{s_1n(\sigma_1^1)}\right)\cdots
  \left(f_2(\sigma_k^1)^{s_k}\circ h^{s_1n(\sigma_k^1)}\right) \\ & = &
  f_2(\gamma(\sigma_\alpha^2))\circ
  h^{\sum_{\sigma^1\in\partial\sigma_2}n(\sigma^1)}\, ,
  \end{eqnarray*}
  since $h$ is central). Therefore, for every $2$-cycle $c$, we have
  \[ e(\rho)_1\cdot c-e(\rho)_2\cdot c=
  \sum_{\sigma^2\in c\ }\sum_{\,\sigma^1\in\partial\sigma^2}
  n(\sigma^1)\, , \]
  where $e(\rho)_i$ is the Euler class defined by the function
  $f_i$. This sum is zero, since $c$ has no boundary.

  Step 2: we prove that $e(\rho)$ depends neither on the point $*\in\Sigma$,
  nor on the cellulation, nor on the base point $x_0\in X$.
  Independence with respect to $*$ follows from the fact that $h$ is central
  in $\widetilde{Ord(X)}$, so that no global conjugation of a representation
  can change its Euler class. Moreover, $e(\rho)$
  is defined in terms of the cellulation of $\Sigma$, and hence it is
  invariant under any isotopy of the cellulation of the surface.
  Now, let $C=C_0\cup C_1\cup C_2$ be a cellulation of $\Sigma$.
  Let $\sigma^1$ be a $1$-cell, and $x$ a point of $\sigma^1$ which is not in
  $C_0$: we define a new cellulation
  $C'=C_0'\cup C_1'\cup C_2'$, with
  $C_0'=C_0\cup\{x\}$,
  $C_1'=(C_1\smallsetminus\{\sigma^1\})\cup\{\sigma_1^1,\sigma_2^1\}$
  and $C_2'=C_2$ cutting the $1$-cell $\sigma^1$ into two cells along the
  vertex $x$; the new edges
  $\sigma_1^1$ and $\sigma_2^1$ inherit of the orientation of $\sigma^1$.
  Then the Euler class defined by $C'$ is the same as the one defined by
  $C$. Indeed, if $f:C_1\rightarrow\widetilde{Ord(X)}$ is a function
  which projects on $\rho$, define simply
  $f':C_1'\rightarrow\widetilde{Ord(X)}$
  by $f'(\sigma_1^1)=1$, $f'(\sigma_2^1)=f(\sigma^1)$, and
  $f'(\sigma_\alpha^1)=f(\sigma_\alpha^1)$ for all the other $1$-cells
  $\sigma_\alpha^1$.
  Similarly, if $\sigma^2$ is a $2$-cell, we define a new cellulation
  $C''=C_0''\cup C_1''\cup C_2''$
  with $C_0''=C_0$, $C_1''=C_1\cup\{\sigma_n^1\}$ and
  $C_2''=(C_2\smallsetminus\{\sigma^2\})\cup\{\sigma_1^2,\sigma_2^2\}$
  obtained by cutting the $2$-cell $\sigma_2$ into two, along the new edge
  $\sigma_n^1$, whose ends are vertices in $\partial C_2$.
  Let then $f:C_1\rightarrow\widetilde{Ord(X)}$ be a function which
  projects on $\rho$. The loop defined by $\partial\sigma_1^2$
  (maybe with the opposite orientation) can be written as a word
  $\sigma_n^1\gamma$ in the elements of $C_1''$ (where $\gamma$ stands for a
  word in the elements of $\partial\sigma^2$). We thus put
  $f'(\sigma_n^1)=f(\gamma)$, and $f'(\sigma_\alpha^1)=f(\sigma_\alpha^1)$
  for all the other $1$-cells $\sigma_\alpha^1$.
  Then by construction,
  $f'$ projects on $\rho$. And
  $f(\gamma(\sigma^2))=f'(\gamma(\sigma^2))=f'(\gamma(\sigma_1^2))\circ
  f'(\gamma(\sigma_2^2))$, hence $f'$ defines the same Euler class as $f$.
  In particular, the Euler class defined by $C''$ is still the same as the one
  defined by $C$.
  Now, we use the classical fact that by using these two operations
  (consisting in refining $C$), as well as taking homotopies, from any pair
  of cellulations $C^1$ and $C^2$ of a same surface $\Sigma$, we can find a
  common refinement $C^3$ of $C^1$ and $C^2$.
  Finally, let $x_0$, $x_0'$ be two points in $X$.
  Take a cellulation $C$ of $\Sigma$ and a function
  $f:C_1\rightarrow\widetilde{Ord(X)}_{x_0}$
  which projects on
  $\rho$. Define then $f':C_1\rightarrow\widetilde{Ord(X)}_{x_0'}$
  by $f'=F_{x_0x_0'}\circ f\circ (F_{x_0x_0'})^{-1}$. Clearly, $f'$
  projects on $\rho$ and $f'$ defines the same Euler class as
  $f$, since $h$ commutes with $F_{x_0x_0'}$.
\end{preuve}

\begin{ps}\label{eulermult}
  Let $\Sigma$ be an oriented surface, $\rho$ a representation and
  $h\colon\Sigma'\rightarrow\Sigma$ a covering of degree $d$.
  Then $e(\rho\circ h)=d\cdot e(\rho)$.
  Indeed, it suffices to take a cellulation $C$ of $\Sigma$ which lifts to a
  cellulation of $\Sigma'$, to take a function
  $f\colon C_1\rightarrow Ord(X)$ which projects on $\rho$, and to take
  $f'=f\circ h$ in order to define the Euler classes of $\rho$ and
  $\rho\circ h$ in $\Sigma$ and $\Sigma'$.
\end{ps}

\subsubsection{Milnor's algorithm}\label{SectionAlgMil}

Consider a closed, oriented surface of genus $g$,
$\Sigma_g$, and equip it with a ``standard'' cellulation, featuring
$1$ single $2$-cell, a single vertex and $4g$ edges labelled by $a_i$,
$b_i$, $a_i^{-1}$ and $b_i^{-1}$, for $1\leq i\leq g$. It yields a standard
presentation
$\pi_1\Sigma_g=\left\langle a_1,\ldots,b_g|\Pi_i[a_i,b_i]=1\right\rangle$.
Then a representant of the fundamental class
$c\in H_2(\pi_1\Sigma_g,\dZ)$ is the $2$-cycle consisting of the unique
$2$-cell, equipped with its orientation.

Given a representation $\rho\in Hom(\pi_1\Sigma_g,Ord(X,o))$,
take some $x_0\in X$ and choose an arbitrary lift $\widetilde{\rho(x)}$
for all $x\in\{a_1,b_1,\ldots,a_g,b_g\}$. As we said before, we still denote
by $e(\rho)\in\dZ$ the evaluation of $e(\rho)\in H^2(\Sigma_g,\dZ)$
on the fundamental class $c$.
Thus, by construction of the Euler class we have
\begin{equation}\label{AlgMil}
[\widetilde{\rho(a_1)},\widetilde{\rho(b_1)}]\cdots
[\widetilde{\rho(a_g)},\widetilde{\rho(b_g)}]=h^{e(\rho)}.
\end{equation}

Since $h$ is central in $\widetilde{Ord(X,o,x_0)}$ and commutes with
$F_{x_0x_1}$ for all $x_1\in X\smallsetminus\{x_0\}$, the result of
this product of commutators does not depend on $x_0$ neither on the choices
of the lifts $\widetilde{\rho(a_i)}$, $\widetilde{\rho(b_i)}$.


\subsubsection{Representations in $\psl$}\label{Euleretpsl}

The Lie group $\psl$ is homeomorphic to a solid torus, and
$\pi_1(\psl)\cong\dZ$.
For all $A\in \psl$, taking a lift
$\tilde{A}\in\widetilde{\psl}$ amounts to lifting the homeomorphism
$A\in\Homeop(\dS^1)$ to a homeomorphism of the universal cover of the circle,
$\widetilde{\dS^1}=\dR$ (here we are using the natural injection
$\psl\hookrightarrow\Homeop(\dS^1)$).
In other words, we have the following short exact sequence:
\[ \begin{array}{rcccccccl}
0 & \rightarrow & \dZ & \rightarrow &
\widetilde{\psl} & \rightarrow & \psl & \rightarrow & 1. \\
 & & & & \vertinj & & \vertinj & & \\
 & & & & \Homeop(\dR) & & \Homeop(\dS^1) & &
\end{array}\]

In this diagram, the sign of the generator of $\dZ$ is determined by the
choice of the orientation of $\dS^1$.
As in \cite{Goldman88}, we are going to denote by $z\in\widetilde{\psl}$
the image of this generator $1\in\dZ$.

Of course, all the construction of the Euler class applies here and if
$\widetilde{\rho(a_i)}$, $\widetilde{\rho(b_i)}$ are arbitrary lifts of
$\rho(a_i)$, $\rho(b_i)$ to
$\widetilde{\psl}$, then $e(\rho)\in\dZ$ is given by the following formula:
$$[\widetilde{\rho(a_1)},\widetilde{\rho(b_1)}]\cdots
[\widetilde{\rho(a_g)},\widetilde{\rho(b_g)}]=z^{e(\rho)}.$$
In particular, notice that in the intermediate cover $SL(2,\dR)$ of $\psl$,
the image of $z$ is $-\Id$. Hence, a representation
$\rho\in Hom(\pi_1\Sigma_g,\psl)$ lifts to $SL(2,\dR)$ if and only if
$e(\rho)$ is even.

\subsubsection{Another definition}\label{sectionAutreDef}

The construction on the Euler class that we have just given is oriented towards
the Milnor algorithm, which is essential for the technical propositions proved
in the next section, which will be useful in the sequel. This algorithm also
gives much of the intuition of how the Euler class behaves and how it can be
computed. But the ``true'' definition of the Euler class is in terms of the
(bounded) group cohomology of $\Gamma$, and we give it here. We refer to
\cite{MacLane,Ghys,Thurston} for more details, and to \cite{Calegari} for a
detailed overview.

Let $\Gamma$ be a group. One forms the {\em homogeneous complex} of $\Gamma$
as follows. We denote by $C_n(\Gamma)$ the free abelian group generated by the
equivalence classes of the $(n+1)$-tuples of $\Gamma$, under the diagonal
action of $\Gamma$: for all $\gamma\in\Gamma$,
$(\gamma_0,\ldots,\gamma_n)\sim (\gamma\gamma_0,\ldots,\gamma\gamma_n)$.
The boundary operator is defined as
\[\partial(\gamma_0,\ldots,\gamma_n)=
\sum_{i=0}^n (-1)^i(\gamma_0,\ldots,\hat{\gamma_i},\ldots,\gamma_n),\]
where $(\gamma_0,\ldots,\hat{\gamma_i},\ldots,\gamma_n)$ is the $n$-tuple
formed by forgetting the $(i+1)$-th term $\gamma_i$. If $A$ is a ring, the
homology of the dual complex
$Hom(C_*(\Gamma),A)$ is denoted $H^*(\gamma;A)$; it is the
{\em cohomology of the group $\Gamma$}. This cohomology of $\Gamma$ identifies
to that of $K(\Gamma,1)$.

Now let $(X,o)$ be a cyclically ordered set, and let
$\rho\colon\Gamma\rightarrow Ord(X,o)$ be a representation. We then define a
$2$-cocycle $c$ by letting
$c(\gamma_1,\gamma_2,\gamma_3)=o(\rho(\gamma_1),\rho(\gamma_2),\rho(\gamma_3))$.

If $\Gamma=\pi_1\Sigma_g$ is the fundamental group of the (compact, closed,
orientable) surface of genus $g$, then we have $[c]=2 e(\rho)$.
The following remark follows:
\begin{ps}\label{pscohom}
  Let $\Gamma$ be a group such that $H^2(\Gamma,\dZ)=0$, and let
  $\rho_1\colon\pi_1\Sigma_g\rightarrow\Gamma$,
  $\rho_2\colon\Gamma\rightarrow Ord(X,o)$.
  Then $e(\rho_1\circ\rho_2)=0$.
\end{ps}

\subsubsection{Finite sets suffice}

Denote by $\dF_{2g}$ the free group on the set
$\{a_1,b_1,\ldots,a_g,b_g\}$, and
$w=[a_1,b_1]\cdots [a_g,b_g]$. The images under the canonical surjection
$\dF_{2g}\rightarrow\pi_1\Sigma_g$ of the subwords of $w$ form a set
$P$, and in all the sequel of this text we denote by
$P_{ref}$ the set $P\cup P^{-1}$.
A major interest of Milnor's algorithm is that we need only finitely many
informations, concerning action of the finite set
$P_{ref}$ on the ordered set
$X$, in order to be able to compute the Euler class of a representation.

This is the key idea which will prove, in chapter
\ref{ChapitreCompact}, that the Euler class extends {\em continuously}
to the boundary of $m_g^o$.


More precisely, the idea is the following:

\begin{prop}\label{nbrefini1}
  Let $(X,o)$ and $(X',o')$ be two cyclically ordered sets, equipped with
  ``base points'' $x_0$, $x_0'$. Let
  $\rho:\pi_1\Sigma_g\rightarrow Ord(X,o)$ and
  $\rho':\pi_1\Sigma_g\rightarrow Ord(X',o')$ be two representations.
  Suppose that for all
  $g_1,g_2,g_3\in P_{ref}$,
  $$o(g_1 x_0,g_2 x_0,g_3 x_0)=o'(g_1 x_0,g_2 x_0,g_3 x_0).$$

  Then $e(\rho_1)=e(\rho_2)$.
\end{prop}

In fact we shall need a slightly more subtle statement, since we want the
Euler class to be stable under small degenerations. We shall leave the proof
of proposition \ref{nbrefini1} as an (easy) exercise, and instead we will
prove the following:

\begin{prop}\label{nbrefini2}
  Let $(X,o)$ and $(X',o')$ be two cyclically ordered sets, equipped with
  ``base points'' $x_0$, $x_0'$. Let
  $\rho:\pi_1\Sigma_g\rightarrow Ord(X,o)$ and
  $\rho':\pi_1\Sigma_g\rightarrow Ord(X',o')$ be two representations.
  Let $y_0\in X$ and $y_0'\in X'$.
  Suppose also that $x_0\not\in P_{ref}\cdot y_0$, that
  $card(P_{ref}\cdot y_0)\geq 2$, and that for all
  $g_1,g_2,g_3\in P_{ref}$,
  $$o(g_1 x_0,g_2 x_0,g_3 y_0)=1\Rightarrow o'(g_1 x_0',g_2 x_0',g_3 y_0')=1$$
  and
  $$o(g_1 x_0,g_2 y_0,g_3 y_0)=1\Rightarrow o'(g_1 x_0',g_2 y_0',g_3 y_0')=1.$$
  Then $e(\rho_1)=e(\rho_2)$.
\end{prop}

Everything relies on the two following elementary lemmas:

\begin{lem}\label{lemtechniqueordre1}
  Let $(X,o)$ be a cyclically ordered set,
  and $f\in Ord(X)$. Suppose we have a base point $x_0\in X$
  and an element $y\in X\smallsetminus\{x_0\}$ such that
  $f(y)\neq x_0$. Denote by $\tilde{f}$ the lift of $f$ to $\dZ\times X$
  verifying $\tilde{f}(0,x_0)=(0,f(x_0))$, and denote by $n$ the integer
  such that $\tilde{f}(0,y)=(n,f(y))$. Then $n$ depends only on
  $o(x_0,f(x_0),f(y))$. More precisely, $n=max(0,-o(x_0,f(x_0),f(y)))$.
\end{lem}

\begin{lem}\label{lemtechniqueordre2}
  Let $(X,o)$ be a cyclically ordered set, $f\in Ord(X)$ and $x_0\in X$.
  Let $x_1,x_2\in X$
  be such that $o(x_0,x_1,x_2)=o(f(x_0),x_1,x_2)=1$.
  Then there exists a lift $\tilde{f}$ of $f$ to $\dZ\times X$ such that
  $(-1,x_2)<_{x_0}\tilde{f}(0,x_0)<_{x_0}(0,x_1)$.
  Moreover, if $y\in X$ is such that
  $o(x_1,x_2,y)\leq 0$ and $o(x_1,x_2,f(y))\leq 0$, then this lift
  $\tilde{f}$ verifies
  $\tilde{f}(0,y)=(0,f(y))$.
\end{lem}

\begin{preuvede} of lemma \ref{lemtechniqueordre1}.

  We have $f(y)\neq f(x_0)$ hence $y\neq x_0$
  and hence $(0,x_0)<_{x_0}(0,y)<_{x_0}(1,x_0)$. The function $\tilde{f}$
  is increasing, so that
  $(0,f(x_0))<_{x_0}(n,f(y))<_{x_0}(1,f(x_0))$.

  If $o(x_0,f(x_0),f(y))=0$ then $x_0=f(x_0)$ and hence $n=0$.

  If $o(x_0,f(x_0),f(y))=1$
  then $(0,f(x_0))<_{x_0}(0,f(y))<_{x_0}(1,f(x_0))$ and then $n=0$.

  If $o(x_0,f(x_0),f(y))=-1$
  then $(-1,f(x_0))<_{x_0}(0,f(y))<_{x_0}(0,f(x_0))$ and in that case $n=1$.
\end{preuvede}

\begin{preuvede} of lemma \ref{lemtechniqueordre2}.

  There are three cases to consider here.
  \begin{itemize}
    \item If $f(x_0)=x_0$, then
      $(-1,x_2)<_{x_0}(0,f(x_0))<_{x_0}(0,x_1)$, so we take
      $\tilde{f}$ such that $\tilde{f}(0,x_0)=(0,f(x_0))$. We then have
      $(0,x_0)<_{x_0}(0,y)<_{x_0}(1,x_0)$, hence, by applying
      $\tilde{f}$ (which is strictly increasing):
      $(0,x_0)<_{x_0}\tilde{f}(0,y)<_{x_0}(1,x_0)$, so that
      $\tilde{f}(0,y)=(0,f(y))$.
    \item If $o(x_0,x_1,f(x_0))=-1$ then
      $(-1,x_2)<_{x_0}(0,x_0)<_{x_0}(0,f(x_0))<_{x_0}(0,x_1)$
      so we take again $\tilde{f}$ such that $\tilde{f}(0,x_0)=(0,f(x_0))$.
      By lemma \ref{lemtechniqueordre1}, it suffices to prove that
      $o(x_0,f(x_0),f(y))\geq 0$ in order to have $n=0$ and thus
      $\tilde{f}(0,y)=(0,f(y))$. But if $o(x_0,f(y),f(x_0))=1$, since
      $o(x_0,f(x_0),x_1)=1$ we get
      $o(x_0,f(y),x_1)=1$ so $o(x_1,x_0,f(y))=1$ which, together with
      $o(x_1,x_2,x_0)=1$, gives $o(x_1,x_2,f(y))=1$,
      a contradiction.
    \item If $o(x_0,x_1,f(x_0))=1$ then we have $o(x_0,x_2,f(x_0))=1$
      (indeed, this follows from the equalities
      $o(f(x_0),x_0,x_1)=o(f(x_0),x_1,x_2)=1$), hence
      \[(-1,x_2)<_{x_0}(-1,f(x_0))<_{x_0}(0,x_0)<_{x_0}(0,x_1)\]
      and we take $\tilde{f}$ such that
      $\tilde{f}(0,x_0)=(-1,f(x_0))$. In particular,
      $\tilde{f}=\tilde{f}'\circ h^{-1}$,
      where $\tilde{f}'(0,x_0)=(0,f(x_0))$. And we have $o(x_0,f(x_0),f(y))=-1$
      (indeed, if $f(y)=x_1$ or $x_2$ we already have this equality, and
      otherwise $o(x_2,x_1,f(y))=1$, which, together with the equality
      $o(x_2,x_0,x_1)=1$, gives
      $o(x_2,x_0,f(y))=1$, \textsl{i.e.} $o(x_0,f(y),x_2)=1$, which, together
      with $o(x_0,x_2,f(x_0))=1$, yields $o(x_0,f(y),f(x_0))=1$) therefore
      lemma \ref{lemtechniqueordre1} applied to $\tilde{f}'$ implies that
      $\tilde{f}'(0,y)=(1,f(y))$, whence $\tilde{f}(0,y)=(0,f(y))$ once again.
  \end{itemize}
\end{preuvede}

We shall now prove proposition \ref{nbrefini2}.
We stated proposition \ref{nbrefini1} only in order to set up the idea
of proposition \ref{nbrefini2}, and we shall not use it. The proof is
similar, though a little easier, than the proof of proposition \ref{nbrefini2}.

\begin{preuvede} of proposition \ref{nbrefini2}.

  Denote
  $y_i=[\rho(a_{i+1}),\rho(b_{i+1})]\cdots[\rho(a_g),\rho(b_g)]\cdot y_0$.
  Then in particular $y_g=y_0$.
  The integers $m_i$ such that
  $[\rho(a_i),\rho(b_i)](0,y_i)=(m_i,y_{i-1})$
  do not depend on the choices of the lifts
  $\widetilde{\rho(a_i)}$, $\widetilde{\rho(b_i)}$, and
  $e(\rho)=\masum{i=1}{g}m_i$, by Milnor's algorithm.
  We also use the same notations in $X'$.

  The finite set $P_{ref}\cdot y_0\subset X\smallsetminus\{x_0\}$,
  equipped with the order $<_{x_0}$, contains a smallest element $x_1$ and
  a biggest element $x_2$. Similarly we define $x_1'$ and $x_2'$ in $X'$.
  Let then $\gamma_1$ be an element of
  $P_{ref}$ such that $\gamma_1y_0'=x_1'$. Then for all
  $\gamma\in P_{ref}$, $o'(x_0',\gamma_1y_0',\gamma y')\leq 0$, hence
  $o(x_0,\gamma_1y_0,\gamma y_0)\leq 0$, so that $\gamma_1y_0$ is minimal among
  $P_{ref}\cdot y_0$ in $X\smallsetminus\{x_0\}$ for the order $<_{x_0}$,
  that is, $\gamma_1y_0=x_1$. Similarly, $x_2$ and $x_2'$ correspond to
  (at least) one same element $\gamma_2\in P_{ref}$.
  Moreover, since $Card(P_{ref}\cdot y_0)\geq 2$, we have $x_1<_{x_0}x_2$,
  hence $o(x_0,x_1,x_2)=1$.

  For every element $\gamma\in\{a_1,b_1,\ldots,a_g,b_g\}$ we define
  $\widetilde{\rho(\gamma)}$ and $\widetilde{\rho'(\gamma)}$ as follows.
  If $\rho(\gamma)\cdot x_0\neq x_0$, we choose $\widetilde{\rho(\gamma)}$
  such that $\widetilde{\rho(\gamma)}(0,x_0)=(0,\rho(\gamma)\cdot x_0)$;
  and we choose $\widetilde{\rho'(\gamma)}$ such that
  $\widetilde{\rho'(\gamma)}(0,x_0')=(0,\rho'(\gamma)\cdot x_0')$.
  Otherwise, if $\rho(\gamma)\cdot x_0=x_0$ then we have
  $o(x_0,x_1,x_2)=o(\rho(\gamma)\cdot x_0,x_1,x_2)=1$ hence,
  by lemma \ref{lemtechniqueordre2}, $\rho(\gamma)$ possesses a lift
  $\widetilde{\rho(\gamma)}$ such that
  $(-1,x_2)<_{x_0}\widetilde{\rho(\gamma)}(0,x_0)<_{x_0}(0,x_1)$,
  and in that case again we have
  $o'(x_0',x_1',x_2')=o'(\rho'(\gamma)\cdot x_0',x_1',x_2')=1$
  (indeed, $o(x_0,x_1,x_2)=o(x_0,\gamma_1y_0,\gamma_2y_0)=1$ so that
  $o'(x_0',x_1',x_2')=1$ and, similarly, $o(\gamma x_0,x_1,x_2)=1$ hence
  $o'(\rho'(\gamma)x_0',x_1',x_2')=1$, since
  $\gamma,\gamma_1,\gamma_2\in P_{ref}$)
  hence, still by applying lemma \ref{lemtechniqueordre2},
  we can define
  $\widetilde{\rho'(\gamma)}$ in such a way that
  $(-1,x_2')<_{x_0'}\widetilde{\rho'(\gamma)}(0,x_0')<_{x_0'}(0,x_1')$.

  Now denote by $n_{i_1},n_{i_2},n_{i_3},n_{i_4}$ the integers such that
  \[\widetilde{\rho(b_i)}^{-1}(0,y_i)=(n_{i_4},\rho(b_i)^{-1}\cdot y_i),\]
  \[\widetilde{\rho(a_i)}^{-1}\cdot(0,\rho(b_i)^{-1}\cdot y_i)=
  (n_{i_3},\rho(a_i)^{-1}\rho(b_i)^{-1}\cdot y_i),\]
  \[\ldots,\]
  and similarly we define integers $n_{i_1}'$, \ldots, $n_{i_4}'$.

  Let us first check that $n_{i_4}=n_{i_4}'$. Denote
  $\gamma=[a_{i+1},b_{i+1}]\cdots[a_g,b_g]$ (that is the element of
  $\pi_1\Sigma_g$ defined by the subword following $b_{i-1}$ in $w$).
  \begin{itemize}
    \item
      If $\rho(b_i)\cdot x_0\neq x_0$, then
      $o(x_0,\rho(b_i^{-1})\cdot x_0,\rho(b_i^{-1}\gamma)\cdot y_0)\neq 0$
      (indeed, these three points are all distinct since $\gamma$
      and $b_i^{-1}\gamma$ are in $P_{ref}$),
      and hence
      \[o'(x_0',\rho'(b_i^{-1})\cdot x_0',\rho'(b_i^{-1})\cdot y_i')=
      o(x_0,\rho(b_i^{-1})\cdot x_0,\rho(b_i^{-1})\cdot y_i)\]
      so by lemma
      \ref{lemtechniqueordre1} (applied to $f=\rho(b_i^{-1})$, $y=y_i$
      and to $f=\rho'(b_i^{-1})$ and $y=y_i'$)
      we have $n_{i_4}=n_{i_4}'$.
    \item
      If $\rho(b_i)\cdot x_0=x_0$, then by lemma
      \ref{lemtechniqueordre2},
      this time we have $n_{i_4}=n_{i_4}'=0$.
  \end{itemize}
  Similarly we get $n_{i_3}-n_{i_4}=n_{i_3}'-n_{i_4}'$, \ldots,
  $n_{i_1}-n_{i_2}=n_{i_1}'-n_{i_2}'$,
  so that $m_i=m_i'$, and hence $e(\rho)=e(\rho')$.
\end{preuvede}

\subsection{Rigidity of non discrete representations}\label{SectionOrdreCaracterise}

In all this section, $o$ will denote the natural cyclic order on
$\partial\dH^2=\dS^1$. The object of this section is to prove the following
statement:
\begin{theo}\label{ordrecaracterise}
  Let $\rho_1,\rho_2\in R_g$ be two non discrete, non elementary
  representations, and let $a_1,a_2\in\partial\dH^2$. Suppose that for all
  $\gamma_1,\gamma_2,\gamma_3\in\pi_1\Sigma_g$,
  \[ o(\rho_1(\gamma_1)\cdot a_1,\rho_1(\gamma_2)\cdot a_1,
  \rho_1(\gamma_3)\cdot a_1)=o(\rho_2(\gamma_1)\cdot a_2,
  \rho_2(\gamma_2)\cdot a_2,\rho_2(\gamma_3)\cdot a_2). \]
  Then there exists $g\in\psl$ such that $\rho_1=g\rho_2 g^{-1}$.
\end{theo}

The proof of theorem \ref{ordrecaracterise} depends on the following lemma:
\begin{lem}\label{ordrecaractech}
  For any $A\in\psl$ and any $a\in\partial\dH^2$:
  \begin{itemize}
    \item[(1)] $o(A^{n_1}a,A^{n_2}a,A^{n_3}a)=0$ for all
      $(n_1,n_2,n_3)\in\dZ^3$ if and only if $a$ is a fixed point of $A^2$,
    \item[(2)] if $o(A^{n_1}a,A^{n_2}a,A^{n_3}a)=1$ for all
      $(n_1,n_2,n_3)\in\dZ^3$
      such that $n_1<n_2<n_3$ then $A$ is parabolic or hyperbolic,
    \item[(3)] similarly, if $o(A^{n_1}a,A^{n_2}a,A^{n_3}a)=1$ for all
      $(n_1,n_2,n_3)\in\dZ^3$
      such that $n_1<n_2<n_3$ then $A$ is parabolic or hyperbolic,
    \item[(4)] otherwise, $A$ is a non trivial elliptic element of $\psl$, and
      the data of $o(A^{n_1}a,A^{n_2}a,A^{n_3}a)$ determines $A$ up to
      conjugation.
  \end{itemize}
\end{lem}

\begin{preuve}
  If $A$ is parabolic or hyperbolic, then one of the cases {\it (1)}, {\it (2)}
  or
  {\it (3)} arises depending on whether $a$ is a fixed point of $A$, and on the
  direction in which $A$ moves the region of $\dS^1$ minus the fixed point(s)
  of $A$, containing $a$. Conversely, it is clear that the cases {\it (2)} and
  {\it (3)} can happen only if $A$ is parabolic or hyperbolic.

  Now suppose that $A$ is elliptic. Then, up to conjugating $A$, we have
  $A=\left(\begin{array}{cc}\cos\theta & -\sin\theta \\ \sin\theta &
  \cos\theta \end{array}\right)$, for some $\theta\in[0,\pi)$.
  In particular, we note that in this case, the function
  $(n_1,n_2,n_3)\mapsto o(A^{n_1}a,A^{n_2},A^{n_3})$ does not depend on $a$
  any more, since $A$ acts as a rotation of angle $2\theta$ on the circle
  $\dS^1\simeq\dR/2\pi\dZ$. Now remark that $o(a,Aa,A^2a)$ equals $1$ if
  $2\theta\in(0,\pi)$, it equals $-1$ if $2\theta\in(\pi,2\pi)$ and it equals
  $0$ otherwise. Similarly, if
  $\frac{\theta}{\pi}=\frac{\alpha_1}{2}+\frac{\alpha_2}{4}+\cdots+
  \frac{\alpha_n}{2^n}+\cdots$ is the dyadic expression of $\frac{\theta}{\pi}$
  (we do not consider expressions with $\alpha_n=1$ for all $n$ large enough,
  hence this expression is unique), then:
  \begin{itemize}
    \item if $o(a,A^{2^{n-1}}a,A^{2^n}a)=1$ then $\alpha_n=1$ and the
      expression does not stop at $\alpha_n$ (\textsl{i.e.}, one of the
      $\alpha_{n+1}$, \ldots, is not zero),
    \item if $o(a,A^{2^{n-1}}a,A^{2^n}a)=-1$ then $\alpha_n=0$ and the
      expression continues after $\alpha_n$,
    \item if $o(a,A^{2^{n-1}}a,A^{2^n}a)=0$ then the angle of the rotation
      $A^{2^{n-1}}$ is $0$ or $\pi$, \textsl{i.e.} $\alpha_k=0$ for all
      $k\geq n+1$.
  \end{itemize}
  In particular, the data of $o(A^{n_1}a,A^{n_2}a,A^{n_3}a)$ for all
  $(n_1,n_2,n_3)\in\dZ^3$ determines the angle $2\theta\in[0,2\pi)$, unless
  if $2\theta=0$ or $\pi$ (since in both cases, this number is always $0$).
\end{preuve}

\begin{ps}
  As we said in section \ref{sectionAutreDef}, the cocycle we are computing
  here is twice the Euler class. In particular, lemma \ref{ordrecaractech}
  is related to the fact that $H_b^2(\dZ,\dZ)=\dR/\dZ$ (see \cite{Ghys87},
  proposition 3-1). Since the element $\frac{1}{2}$ has order $2$ in
  $\dR/\dZ$, and since we are computing twice the Euler class, lemma
  \ref{ordrecaractech} does not distinguish the rotation of half the circle
  from the identity.
\end{ps}

\begin{preuvede} of theorem \ref{ordrecaracterise}.

  It is well-known (see \textsl{e.g.} \cite{Katok}, exercise 2.14 or
  \cite{TheseMaxime}, proposition 1.1.14 for a proof) that every non discrete,
  non
  elementary subgroup of $\psl$ contains an elliptic element of infinite order.
  Let $\gamma\in\pi_1\Sigma_g$ be such that $\rho_1(\gamma)$ is elliptic of
  infinite order. Conjugate $\rho_1$ so that
  $\rho_1(\gamma)=\left(\begin{array}{cc}\cos\theta & -\sin\theta \\ \sin\theta
  & \cos\theta \end{array}\right)$, for some irrational $\theta\in[0,\pi)$.
  Now lemma \ref{ordrecaractech} implies that $\rho_2(\gamma)$
  is conjugate to $\rho_1(\gamma)$.
  Now let $\gamma'\in\pi_1\Sigma_g$. Write
  $\rho_1(\gamma)=\left(\begin{array}{cc}a & b \\ c & d \end{array}\right)$.
  Then
  \begin{equation}\label{equordcarac}
    \Tr(\rho_1(\gamma^n\gamma'))=
    \left|(a+d)\cos(n\theta)+(b-c)\sin(n\theta)\right|
  \end{equation}
  is the absolute value of the
  scalar product of the Euclidean plane vectors
  $\left(\begin{array}{c}a+d \\ b-c\end{array}\right)$ and
  $\left(\begin{array}{c}\cos(n\theta) \\ \sin(n\theta)\end{array}\right)$.
  The latter can be chosen to be arbitrarily close to any direction, since
  $\theta$ is irrational, and the former has norm greater or equal to $4$
  (since $ad-bc=1$), thus
  $\Tr(\rho_1(\gamma^n\gamma'))\in[\frac{1}{2},\frac{3}{2}]$ for infinitely
  many values of $n$; hence $\rho_1(\gamma^n\gamma')$ is elliptic, of order
  different from $2$, for infinitely many values of $n$. For these values
  it follows from lemma \ref{ordrecaractech} that
  $\Tr(\rho_1(\gamma^n\gamma'))=\Tr(\rho_2(\gamma^n\gamma'))$, so that, from
  the formula (\ref{equordcarac}) we derive that
  $\Tr(\rho_1(\gamma'))=\Tr(\rho_2(\gamma'))$. It now follows from
  proposition \ref{mgcaracteres} that $\rho_1$ and $\rho_2$ are conjugate
  by an element of $Isom(\dH^2)$, but a conjugation reversing the orientation
  would change $\rho_1(\gamma)$ into a rotation of angle $-\theta$ instead
  of $\theta$, hence $\rho_1$ and $\rho_2$ are conjugate by an element of
  $\psl$.
\end{preuvede}

\subsection{Almost faithful morphisms}\label{SectionPresqueInj}

The connectedness of $\overline{m_g^u}$ (theorem \ref{cc}) strongly relies
on a property of surface groups related to the fact that these groups are
``limit groups''
(see \textsl{e.g.} \cite{Sela01,Guirardel03,ChampetierGuirardel}).

Let us fix some notations first.
Fix a standard presentation of the fundamental group
\[\pi_1\Sigma_g=\left\langle a_1,\ldots,b_g \,|\, [a_1,b_1]\cdots[a_g,b_g]=1
\right\rangle\]
of the surface $\Sigma_g$.
The set $S=\{a_1,\ldots,b_g\}$ generates $\pi_1\Sigma_g$, and we denote by
$\dF$ the free group on this generating set; and denote by
$\pi\colon\dF\rightarrow\pi_1\Sigma_g$ the canonical surjection.
Denote by $B_n\subset\pi_1\Sigma_g$ the ball of centre $1$ and radius $n$
for the Cayley metric associated to the generating set $S$. Finally, we denote
by $\dF_k$ a free group of rank $k$.

A group $\Gamma$ is said to be {\em residually free} if for all
$\gamma\in\Gamma\smallsetminus\{1\}$, there exists a morphism
$\varphi\colon\Gamma\rightarrow\dF_2$ such that $\varphi(\gamma)\neq 1$.
We say that $\Gamma$ is {\em fully residually free} if for every finite subset
$\{\gamma_1,\ldots,\gamma_n\}\subset\Gamma\smallsetminus\{1\}$, there exists
a morphism $\varphi\colon\Gamma\rightarrow\dF_2$ such that for all
$i\in\{1,\ldots,n\}$, $\varphi(\gamma_i)\neq 1$.
We will use the following celebrated result:

\begin{theo}[G. Baumslag \cite{Baumslag62}]\label{lemtechMCG1}
  For all $g\geq 2$, the group $\pi_1\Sigma_g$ is fully residually free.
  In other words, for every $n\geq 0$, there exists a morphism
  $\varphi_n\colon\pi_1\Sigma_g\rightarrow\dF_2$ such that
  $\ker(\varphi_n)\cap B_n=\{1\}$.
\end{theo}

Heuristically, the morphisms $\varphi_n$ are ``more and more injective''.
In the language of \cite{Sela01,ChampetierGuirardel}, the group $\pi_1\Sigma_g$
is a ``limit group'' of the group $\dF_2$. In fact, we will need a statement
a little more general: we will need to make two explicitly given morphisms
``more and more injective'', by composing them with automorphisms of
$\pi_1\Sigma_g$.
For all $g\geq 3$, we denote
$e_g\colon\pi_1\Sigma_g\rightarrow\pi_1\Sigma_{g-1}$ the morphism consisting
of collapsing the last handle.
More precisely, given the two standard presentations
\[\pi_1\Sigma_g=\left\langle a_1,\ldots,b_g\left|\,\prod_{i=1}^g[a_i,b_i]=1
\right.\right\rangle\]
and
\[\pi_1\Sigma_{g-1}=\left\langle a_1,\ldots,b_{g-1}\left|\,
\prod_{i=1}^{g-1}[a_i,b_i]=1\right.\right\rangle,\]
the map $e_g$ is defined by $e_g(\gamma)=\gamma$ for
$\gamma=a_1,b_1,\ldots,a_{g-1},b_{g-1}$ and $e_g(\gamma)=1$ for
$\gamma=a_g,b_g$.
Now let $g\geq 2$;
we will consider a Fuchsian group $G_g$ and a morphism
$p_g\colon\pi_1\Sigma_g\rightarrow G_g$, depending on the parity of $g$.

If $g$ is even, $g=2g'$, we fix a Fuchsian group $G_g$ of signature
$(g';2)$. It is the fundamental group of the hyperbolic orbifold of genus
$g'$ and with one conic singularity of angle $\pi$ ($=\frac{2\pi}{2}$).
This orbifold possesses a covering $f_g$ of degree $2$, branched over the conic
singularity,
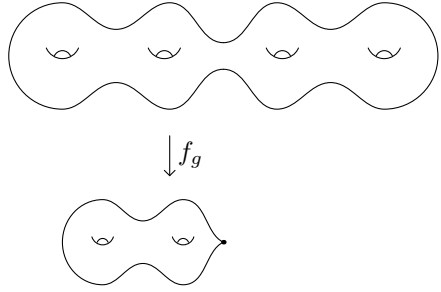
\[\epsfbox{Figures/FiguresExtraTwist.1}\]
and we set $p_g={f_g}_*\colon\pi_1\Sigma_g\rightarrow G_g$.
Recall that the Fuchsian group $G_g$, of signature $(g';2)$, has the following
presentation:
\[G_g=\left\langle\alpha_1,\ldots,\beta_{g'}\left|
\left([\alpha_1,\beta_1]\cdots[\alpha_{g'},\beta_{g'}]
\right)^2=1 \right.\right\rangle.\]
Then we define the map
$p_g\colon\pi_1\Sigma_g\rightarrow G_g$
by letting
$p_g(a_i)=p_g(a_{g'+i})=\alpha_i$ and
$p_g(b_i)=p_g(b_{g'+i})=\beta_i$ for all $i$ between $1$ and $g'$.

In the case when $g$ is odd, $g=2g'+1$, we fix a Fuchsian group
$G_g$ of signature $(g';2,2,2)$. Recall that such a group has the following
presentation:
\[G_g=\left\langle q_1,q_2,\alpha_1,\ldots,\beta_{g'}\left|\,q_1^2=q_2^2=1,
\left(q_1 q_2[\alpha_1,\beta_1]\cdots[\alpha_{g'},\beta_{g'}]
\right)^2=1 \right.\right\rangle.\]
We denote by $p_g\colon\pi_1\Sigma_g\rightarrow G_g$ the morphism defined by
$p_g(a_1)=q_1^{-1}$, $p_g(b_1)=q_2^{-1}$,
and $p_g(a_i)=(q_1q_2)^{-1}\alpha_{i-1}(q_1q_2)$,
$p_g(b_i)=(q_1q_2)^{-1}\beta_{i-1}(q_1q_2)$, for all $i$ between $2$
and $g'+1$, and $p_g(a_i)=\alpha_{i-g'-1}$ and $p_g(b_{i-g'-1})=\beta_i$
when $g'+2\leq i\leq g$. It is a little more delicate, in the case when
$g$ is odd, to draw $p_g$ as a covering.

In the two cases of parity of $g$,
it has been proved (see \cite{noninjrep}, proposition 4.5) that the discrete
representation $p_g$ (as an element of $R_g$) has Euler class $2g-3$; this is
why we consider it.

\begin{lem}\label{lemtechMCG2}
  Let $g\geq 4$. Then for all $n\geq 0$, there exists an element
  $\gamma_n\in Aut(\pi_1\Sigma_g)$
  such that the kernel of the morphism
  $p_g\circ\gamma_n\colon\pi_1\Sigma_g\rightarrow G_g$
  does not contain any non-trivial element of length less than $n$.
\end{lem}
  
In other words, we can ``mix'' the curves representing $a_i$ and
$b_i$ by Dehn twists, so that the map
$p_g$ does not kill any word of given length, that is, we can change this
map into maps which are ``arbitrarily injective''.

Similarly:

\begin{lem}\label{lemtechMCG3}
  Let $g\geq 3$. Then for all $n\geq 0$, there exists an element
  $\gamma_n\in Aut(\pi_1\Sigma_g)$
  such that the kernel of the map
  $e_g\circ\gamma_n\colon\pi_1\Sigma_g\rightarrow G_g$
  does not contain any non-trivial elements of length less than $n$.
\end{lem}

V. Guirardel pointed out to me the following proof, which is due to Z. Sela.

\begin{prop}\label{Vincent}
  Let $\Sigma$ be a surface, possibly with boundary, of Euler characteristic
  less or equal to $-1$. Let $\varphi\colon\pi_1\Sigma\rightarrow\dF$ be a
  morphism of non abelian image in a free group $\dF$, whose restrictions to
  the fundamental groups of the boundary components are injective. Then for all
  finite subset $P\subset\pi_1\Sigma\smallsetminus\{1\}$, there exists a
  diffeomorphism $\gamma_P$ of $\Sigma$, preserving pointwise the boundary
  components, and such that $Ker(\varphi\circ{\gamma_P}_*)\cap P=\emptyset$.
\end{prop}

\begin{cor}\label{enfinlibre}
  Let $g\geq 2$, and let $\varphi\colon\pi_1\Sigma_g\rightarrow\dF$ be a
  morphism of non abelian image. Then for all $n$, there exists
  $\gamma_n\in Aut(\pi_1\Sigma_g)$ such that
  $Ker(\varphi\circ\gamma_n)\cap B_n=\{1\}$.
\end{cor}

The lemmas \ref{lemtechMCG2} and \ref{lemtechMCG3} follow:

\begin{preuvede} of lemma \ref{lemtechMCG2}.

The map $\varphi_g\colon G_g\rightarrow\dF_2$ defined by
$\varphi_g(a_1)=x$, $\varphi_g(a_2)=y$ and $\varphi_g(u)=1$ for all the other
generators $u$ of $G_g$, with $\dF_2=\langle x,y\rangle$, is a morphism of non
abelian image, and
$\varphi_g\circ p_g\colon\pi_1\Sigma_g\rightarrow\dF_2$ is therefore a
morphism satisfying the hypotheses of corollary \ref{enfinlibre}.
Thus, we can conjugate it by automorphisms of $\pi_1\Sigma_g$ in order to make
it ``arbitrarily injective''.
\end{preuvede}

\begin{preuvede} of lemma \ref{lemtechMCG3}.

  If $g\geq 3$, the map
  $\varphi_g\colon \pi_1\Sigma_{g-1}\rightarrow\dF_2$ defined by
  $\varphi_g(a_1)=x$, $\varphi_g(a_2)=y$ and $\varphi_g(u)=1$ for all the other
  generators $\pi_1\Sigma_{g-1}$ is still a morphism of non abelian image.
\end{preuvede}

The proof of the proposition relies on the following two lemmas.

\begin{lem}[Z. Sela \cite{Sela01}, lemma 5.13]
  Let $\Sigma$ be a surface, possibly with boundary, of Euler characteristic
  less or equal to $-1$. Let $\varphi\colon\pi_1\Sigma\rightarrow\dF$ be a
  morphism of non abelian image in a free group $\dF$, whose restrictions to
  the fundamental groups of the boundary components are injective.
  Then there exists a family of disjoint closed simple curves
  $c_1$, \ldots, $c_p$ in $\Sigma$, which cut $\Sigma$ into pairs of pants,
  and such that the restriction of $\varphi$ to the fundamental group of each
  pair of pant is injective.
\end{lem}

\begin{lem}[G. Baumslag \cite{Baumslag62}, proposition 1]
  Let $\dF$ be a free group and let $a_1$, \ldots, $a_n$, $c\in\dF$ be such
  that $c$ does not commute with any of the $a_i$'s. Then for all
  $k_0,\ldots,k_n$ large enough, the element
  $c^{k_0}a_1 c^{k_1}a_2\cdots c^{k_{n-1}}a_n c^{k_n}$ is non trivial in $\dF$.
\end{lem}

\begin{preuvede} of proposition \ref{Vincent}.

  Denote by $\chi(\Sigma)$ and $g(\Sigma)$ the Euler characteristic and the
  genus of $\Sigma$. We shall work by induction on
  $(-\chi(\Sigma),g(\Sigma))$, following the lexicographic order. The number
  $g(\Sigma)$ being between $0$ and $1-\frac{\chi(\Sigma)}{2}$,
  the element $(-\chi(\Sigma),g(\Sigma))$ indeed describes a set which is in
  bijection preserving the order with $\dN$.

  If $\chi(\Sigma)=-1$, then $\varphi$ is injective (see \textsl{e.g.}
  \cite{ChampetierGuirardel}, proposition 3.1). Hence, suppose that
  the proposition is true for all $\Sigma'$ such that
  $(-\chi(\Sigma'),g(\Sigma'))<(-\chi(\Sigma),g(\Sigma))$ (for the
  lexicographic order) and consider curves $c_1,\ldots,c_p$ as in Z. Sela's
  lemma.

  Suppose first that $c_1$ is a separating curve: denote
  $\Sigma=\Sigma_1\displaystyle{\mathop{\cup}_{c_1}}\Sigma_2$. Put the base
  point near $c_1$, on the side of $\Sigma_1$. We have
  $\pi_1\Sigma=\pi_1\Sigma_1\displaystyle{\mathop{*}_{\alpha}}\pi_1\Sigma_2$,
  where $\alpha$ is represented by the curve $c_1$, deformed so that it passes
  through the base point. Fix a finite subset $P\subset\pi_1\Sigma$.
  For every $m\in P$, choose a writing
  $m=a_1\alpha^{k_1}b_1\alpha^{l_1}\cdots a_n\alpha^{k_n}b_n\alpha^{l_n}$,
  with $a_i\in\pi_1\Sigma_1$ and $b_i\in\pi_1\Sigma_2$, and such that
  $a_i$, $b_i$ do not commute with $\alpha$ (except maybe $a_1$ or $b_n$,
  in which case we do not write them in $m$). Denote by $P_1$ the subset of
  $\pi_1\Sigma_1$ defined by the elements $\alpha a_i\alpha^{-1}a_i^{-1}$
  and denote by $P_2$ the subset of $\pi_1\Sigma_2$ defined by the elements
  $\alpha b_i\alpha^{-1}b_i^{-1}$. By induction hypothesis, there exists a
  diffeomorphism $\gamma_1$ of $\Sigma_1$ fixing the boundary of $\Sigma_1$
  (as well as the boundary of the curve $c_1$), and a diffeomorphism
  $\gamma_2$ of $\Sigma_2$ fixing the boundary of $\Sigma_2$ such that for all
  $u\in P_1$ we have $\varphi\circ{\gamma_1}_*(u)\neq 1$ and such that for all
  $u\in P_2$ we have $\varphi\circ{\gamma_2}_*(u)\neq 1$. Consider then a
  diffeomorphism $\gamma_k\colon\Sigma\rightarrow\Sigma$ defined by $\gamma_1$
  on $\Sigma_1$, $\gamma_2$ on $\Sigma_2$ and by $k$ Dehn twists along $c_1$.
  Then ${\gamma_k}_*\colon\pi_1\Sigma\rightarrow\pi_1\Sigma$ is defined as
  follows: if $a_1,\ldots,a_n\in\pi_1\Sigma_1$ and
  $b_1,\ldots,b_n\in\pi_1\Sigma_2$, we have
  \[ {\gamma_k}_*\left(a_1\alpha^{k_1}b_1\alpha^{l_1}\cdots a_n\alpha^{k_n}b_n
  \alpha^{l_n}\right)\!=\!{\gamma_1}_*(a_1)\alpha^{k_1+k}{\gamma_2}_*(b_1)
  \alpha^{l_1-k}\cdots
  {\gamma_1}_*(a_n)\alpha^{k_n+k}{\gamma_2}_*(b_1)\alpha^{l_n-k}. \]
  Let $m\in P$, and consider the expression
  $m=a_1\alpha^{k_1}b_1\alpha^{l_1}\cdots a_n\alpha^{k_n}b_n\alpha^{l_n}$
  chosen before. We then have
  \small
  \[ \varphi\circ{\gamma_k}_*\left(m\right)=
  \varphi\circ{\gamma_1}_*(a_1)\varphi(\alpha)^{k_1+k}
  \varphi\circ{\gamma_2}_*(b_1)\varphi(\alpha)^{l_1-k}\cdots
  \varphi\circ{\gamma_1}_*(a_n)\varphi(\alpha)^{k_n+k}
  \varphi\circ{\gamma_2}_*(b_n)\varphi(\alpha)^{l_n-k}\, . \]
  \normalsize
  Since $\alpha a_i\alpha^{-1}a_i^{-1}\in P_1$, we have
  $\varphi\circ{\gamma_1}_*(\alpha a_i\alpha^{-1}a_i^{-1})\neq 1$, but
  ${\gamma_1}_*(\alpha)=\alpha$: hence $\varphi\circ{\gamma_1}_*(a_i)$ does
  not commute with $\varphi(\alpha_1)$. All the conditions of G. Baumslag's
  lemma are satisfied, and hence for all $k$ large enough,
  $\varphi\circ{\gamma_k}_*$ sends every non trivial element of $P$ on a non
  trivial element of $\dF$.

  Suppose finally that $c_1$ is a non separating curve. Then
  $\Sigma$, this time, is obtained by glueing two boundary components of
  a surface $\Sigma_1$, and we have $-\chi(\Sigma)=-\chi(\Sigma_1)$
  but $g(\Sigma_1)<g(\Sigma)$. We have
  $\pi_1\Sigma=\langle\pi_1\Sigma_1,t|t^{-1}\alpha t=\beta\rangle$, where
  $\alpha,\beta\in\pi_1\Sigma_1$ are represented by the boundary components of
  $\Sigma_1$ concerned by the glueing, and where $t$ is represented, in
  $\Sigma$, by a simple curve intersecting the curve $c_1$ at a single point.
  The elements $m\in\pi_1\Sigma$ can be written (non uniquely) as
  $m=t^{k_0}a_1t^{k_1}\cdots a_nt^{k_n}$ with $a_i\in\pi_1\Sigma_1$; if
  $\gamma$ is a diffeomorphism of $\Sigma_1$ fixing its boundary
  (as well as a neighbourhood of $c_1$ containing the base point) then
  $\gamma_*(m)=t^{k_0}\gamma_*(a_1)t^{k_1}\cdots\gamma_*(a_n)t^{k_n}$, and
  the image of $m$ under $k$ Dehn twists along $c_1$ is equal to
  $(\alpha^k t)^{k_0}a_1(\alpha^k t)^{k_1}\cdots a_n(\alpha^k t)^{k_n}$.
  The same argument as in the preceding case transposes here, thereby finishing
  the induction.
\end{preuvede}

\section{Compactifications, degenerations and orientation}\label{ChapitreCompact}

Before going into the study of the compactifications $\overline{m_g^u}$ and
$\overline{m_g^o}$, we will need to prove a technical fact, namely that the
connected components of the spaces $m_g^o$ and $m_g^u$ are one-ended. This
implies that for every compactification considered, the boundaries of these
connected components are connected spaces.

\subsection{The connected components of $m_g^u$ and $m_g^o$ are one-ended}

Let us begin with some recalls.

\begin{deuf}
  Let $X$ be a connected locally compact space. The suppremum of the number
  of unbounded components (\textsl{i.e.}, whose closure is non compact) of
  $X\smallsetminus K$, as $K$ describes the set of compact subsets of $X$,
  is called the {\em number of ends of $X$}.
\end{deuf}

Here, once again, we are interested in the case $\Gamma=\pi_1\Sigma_g$
and $n=2$; we denote $m_g^o=m_{\pi_1\Sigma_g}^o(2)$, with the same notations
as before. We fix the generating set
\[ S=\left\{a_1,a_1^{-1},b_1,b_1^{-1},\ldots,b_g,b_g^{-1}\right\}. \]
Recall (corollary \ref{dcontinupropre}) that the function
$d\colon m_g^o\rightarrow\dR_+$ defined as
\[ d(\rho)=\min_{x\in\dH^2}\max_{s\in S}d(x,\rho(s)x) \]
is well defined, continuous and proper. In particular $m_g^o$ is locally
compact, and, similarly, $m_g^u$ is locally compact.

We are going to prove that the connected components of $m_g^o$, as well as
those of $m_g^u$, are one-ended. Since the proof is exactly the same in both
cases, until the end of this section (and {\em only} in this section) we
will denote by $m_g$ these representation spaces, and by $m_{g,k}$ the
corresponding connected components, being vague whether we consider the
oriented or the unoriented representation space.

In \cite{Hitchin}, N. Hitchin proved that for all $k\in\{0,\ldots,2g-2\}$,
the connected component $m_{g,k}$ of $m_g$
is homeomorphic to a complex vector bundle of dimension $g-1+|k|$ on the
$(2g-2-|k|)$-th symmetric product of the surface. It follows, in particular,
that the connected component $m_{g,k}$ is one-ended, for every $k\neq 0$.
We shall give a much more elementary proof (than the one of \cite{Hitchin}) of
this result, and generalize it to the case $k=0$.

\begin{prop}\label{oneend}
  For all $k$ such that $|k|\leq 2g-2$, the space $m_{g,k}$ is one-ended.
\end{prop}

First let us fix some notations. If $S'\subset S$, we denote
\[d_{S'}(\rho)=\inf_{x\in\dH^2}\max_{s\in S'}d(x,\rho(s)x).\]
If $r>0$ and $S'\subset S$, we denote
\[K_{S'}^r=\left\{\left.[\rho]\subset m_g^o\right|
d_{S'}(\rho)\leq r\right\},\]
$K^r=K_S^r$, $U_{S'}^r=m_g^o\smallsetminus K_{S'}^r$. If $s\in S$,
$U_{\{s\}}^r$ is therefore the set of conjugacy classes of representations
$\rho$ such that $\rho(s)$ is a hyperbolic element whose translation length
is strictly greater than $r$, which is equivalent to $\Tr(\rho(s))>2\cosh r$.
We denote
$V^r=U_{\{a_1\}}^r\cap U_{\{a_2\}}^r$ et $W^r=U_{\{a_1\}}^r\cup U_{\{a_2\}}^r$.
Of course,
$V^r\subset W^r\subset m_g^o\smallsetminus K^r$.
Also, for all $A\in\psl$, we wish to choose a one-parameter subgroup passing
through $A$. If $A=\Id$, we set $A_t=\Id$, for all $t\in\dR$.
If $A\in Par\cup Hyp$, we choose $A_t$ such that $A_0=\Id$ and $A_1=A$. If
$A\in Ell$, we further require that $A_t\neq\Id$ for all $t\in(0,1]$, and that
$A_t$ is a rotation of positive angle for small $t$.
Note that for all $n\in\dZ$, we have $A_n=A^n$.

We shall actually establish the following result:

\begin{prop}\label{oneendtech}
  For every $k\in\dZ$ such that $|k|\leq 2g-2$, and all $r>0$, every two
  representations
  $\rho_1,\rho_2\in m_{g,k}\smallsetminus K^{r+8g\delta_{\dH^2}}$ can be
  joined by a path in $m_{g,k}\smallsetminus K^r$.
\end{prop}

\begin{lem}\label{lemtechoneend1}
  Let $A,B\in PSL(2,\dR)$ be such that $\Tr(B)>2$. Then for all $r>0$, there
  exist $n\in\dZ$ and $x\in\dR$ such that
  $\Tr\left(A\cdot(BA_x)^n\right)>2\cosh r$.
\end{lem}

\begin{preuve}
  We have $B\in Hyp$, so we can write
  $B=\left(\begin{array}{cc}\lambda & 0 \\ 0 & \frac{1}{\lambda}
  \end{array}\right)$ with $\lambda>1$,
  and
  $A=\left(\begin{array}{cc} a & b \\ c & d \end{array}\right)$
  in some adapted basis.
  Then
  $\Tr(A\cdot B^n)=\left|a\lambda^n+\frac{d}{\lambda^n}\right|$. This trace
  can be made arbitrarily large, as one takes
  $|n|$ big enough, except in the case when $a=d=0$,
  or, equivalently, when $A$ is an elliptic element which exchanges the fixed
  points of $B$ in $\partial_\infty \dH^2$. In that case, if $x$ is small
  enough, then $B\cdot A_x$ is still a hyperbolic element, whose fixed points
  are distinct to those of $B$ in $\partial_\infty\dH^2$
  (indeed, up to conjugating $\rho$ by a diagonal matrix, we have
  $A=\left(\begin{array}{cc}0 & -1 \\ 1 & 0\end{array}\right)$ hence
  $A_x=\left(\begin{array}{cc}\cos \frac{x}{\pi} & -\sin \frac{x}{\pi} \\
  \sin \frac{x}{\pi} & \cos \frac{x}{\pi} \end{array}\right)$ and then
  $B\cdot A_x=\left(\begin{array}{cc}\lambda \cos \frac{x}{\pi} &
  -\lambda\sin \frac{x}{\pi} \\ \frac{1}{\lambda}\sin \frac{x}{\pi} &
  \frac{1}{\lambda}\cos \frac{x}{\pi} \end{array}\right)$ does not fix $0$
  either $\infty$ when $\sin\frac{x}{\pi}\neq 0$).
\end{preuve}

\begin{lem}\label{lemtechoneend2}
  Let $A,B\in PSL(2,\dR)$ such that $[A,B]\neq 1$. Then there exist $x\in\dR$
  and $n\in\dZ$ such that $\Tr\left(B\cdot(AB_x)^n\right)>2$.
\end{lem}

\begin{preuve}
  If $B\in Hyp$, we take $n=0$. If $A\in Hyp$, then we exchange the roles of
  $A$ and $B$, we apply lemma \ref{lemtechoneend1}. If $A$ or $B$ is parabolic,
  and $A,B\not\in Hyp$, then $\Tr\left(B\cdot(AB_x)^n\right)>2$ as soon as $x$
  or $n$ is large enough, as in lemma \ref{lemtechoneend1} (it suffices to
  write the matrices).
  Since $A,B\neq\Id$, the only remaining case is when $A,B\in Ell$. For some
  small $x$, the matrix $A\cdot B_x$ is elliptic, of infinite order
  (indeed, $AB_x$ is still elliptic since traces are continuous, and one proves
  very easily (see \textsl{e.g.} \cite{noninjrep}, lemma 3.2) that this trace
  is non constant, so that the elliptic element $AB_x$ has irrational angle
  for infinitely many $x$),
  and does not commute with $B$. In some adapted basis,
  $B=\left(\begin{array}{cc} \cos \theta_1 & -\sin \theta_1 \\
  \sin \theta_1 & \cos \theta_1 \end{array}\right)$, and
  $AB_x=U^{-1}\cdot\left(\begin{array}{cc} \cos \theta_2 & -\sin \theta_2 \\
  \sin \theta_2 & \cos \theta_2 \end{array}\right)\cdot U$, where
  $U$ sends the fixed point of $B$ on that of $AB_x$. In particular, we can
  require that $U$ is hyperbolic. Since $\theta_2$ is irrational, $(AB_x)^n$
  can be arbitrarily close to $U^{-1}B^{-1}U$, and $\Tr(BU^{-1}B^{-1}U)$ is
  the trace of the commutator of a hyperbolic element with an elliptic element
  in $PSL(2,\dR)$, and we can check by easy computation that this is always
  strictly greater than $2$.
  Therefore, $\Tr\left(B\cdot(AB_x)^n\right)>2$, for some $n\in\dZ$.
\end{preuve}

\begin{lem}\label{lemtechoneend3}
  Let $\rho\in m_g\smallsetminus K^{r+8g\delta_{\dH^2}}$. Then there exist
  $s_1$, $s_2\in\{a_1,\ldots,b_g\}$ such that $d_{\{s_1,s_2\}}(\rho)>r$.
\end{lem}

\begin{preuve}
  If $\rho\in R_g$, $S'\subset S$ and $\alpha>0$ denote
  $F_{S'}^\alpha(\rho)=\left\{x\in\dH^2\left|\forall s\in S',\,
  d(x,\rho(s)x)\leq\alpha\right.\right\}$.

  Let $\alpha>0$, and suppose that
  $F_{\{s_1,s_2\}}^\alpha(\rho)\neq\emptyset$ for every pair
  $\{s_1,s_2\}\subset\{a_1,\ldots,b_g\}$. Then for every triple
  $\{s_1,s_2,s_3\}$, the {\em convex} sets $F_{\{s_1\}}^\alpha(\rho)$,
  $F_{\{s_2\}}^\alpha(\rho)$, $F_{\{s_3\}}^\alpha(\rho)$
  intersect pairwise, hence, since $\dH^2$ is $\delta_{\dH^2}$-hyperbolic,
  there exists a point $x$ at distance at most $2\delta_{\dH^2}$ of each of
  these three sets. This implies that
  $x\in F_{\{s_1,s_2,s_3\}}^{\alpha+4\delta_{\dH^2}}(\rho)$.
  More generally, if $x>0$ is such that for every
  $k$-tuple $S'$ of $\{a_1,\ldots,b_g\}$, $F_{S'}^x\neq\emptyset$, then,
  if $\{s_1,\ldots,s_{k+1}\}\subset\{a_1,\ldots,b_g\}$, the convex sets
  $F_{\{s_1,\ldots,s_{k-1}\}}^x(\rho)$,
  $F_{\{s_1,\ldots,s_{k-2},s_k\}}^x(\rho)$ and
  $F_{\{s_1,\ldots,s_{k-2},s_{k+1}\}}^x(\rho)$ intersect pairwise, and it
  follows that
  $F_{\{s_1,\ldots,s_{k+1}\}}^{x+4\delta_{\dH^2}}(\rho)\neq\emptyset$.
  Therefore, by induction on $Card(\{a_1,\ldots,b_g\})$ we get that
  $d(\rho)\leq \alpha+8g\delta_{\dH^2}$, as soon as
  $F_{\{s_1,s_2\}}^\alpha(\rho)\neq\emptyset$ for every pair
  $\{s_1,s_2\}\subset \{a_1,\ldots,b_g\}$.
  
  Since the inequality $d(\rho)>r+8g\delta_{\dH^2}$ is strict, there exists
  $\varepsilon>0$ such that $d(\rho)>r+8g\delta_{\dH^2}+\varepsilon$, which
  implies that there exists a pair $\{s_1,s_2\}\subset\{a_1,\ldots,b_g\}$ such
  that $F_{\{s_1,s_2\}}^{r+\varepsilon}(\rho)=\emptyset$; this implies that
  $d_{\{s_1,s_2\}}(\rho)\geq r+\varepsilon$ hence
  $d_{\{s_1,s_2\}}(\rho)>r$.
\end{preuve}

\begin{preuvede} of proposition \ref{oneendtech}.

  {\bf Step 1: }{\em Let $\rho_0$, $\rho_1\in V^r\cap e^{-1}(k)$. Then there
  exists a path $\rho_t\in W^r$, joining $\rho_0$ and $\rho_1$.}

  By lemma 10.1 of W. Goldman \cite{Goldman88}, for every $k\in\dZ$ such that
  $|k|\leq 2g-2$, the set of representations $\rho$ such that
  $[\rho(a_i),\rho(b_i)]\neq\Id$ is path-connected and dense in $m_{g,k}$.
  We can thus perturb $\rho_0$ and $\rho_1$, and find a path
  $\rho_t\in m_{g,k}$ joining $\rho_0$ to $\rho_1$ and such that for all
  $t\in[0,1]$ and all $i\in\{1,\ldots,g\}$,
  $[\rho_t(a_i),\rho_t(b_i)]\neq\Id$.

  Let $t\in[0,1]$ and $i\in\{1,2\}$. Denote $A(t)=\rho_t(a_i)$ and
  $B(t)=\rho_t(b_i)$. Then by lemma \ref{lemtechoneend2}, there exist
  $n_i(t)\in\dZ$ and $x_i(t)\in\dR$ such that
  $\Tr\left(B(t)\left(A(t)B(t)_{x_i(t)}\right)^{n_i(t)}\right)>2$,
  and, of course, we have
  \[ \left[A(t)B(t)_{x_i(t)},B(t)\left(A(t)B(t)_{x_i(t)}\right)^{n_i(t)}\right]
  =[A(t),B(t)]. \]
  Hence, by lemma \ref{lemtechoneend1}, there exist $n_i'(t)\in\dZ$ and
  $x_i'(t)\in\dR$ such that
  \begin{equation}\label{eqnend2}
    \Tr\left(A(t)B(t)_{x_i(t)}\left(B(t)(A(t)B(t)_{x_i(t)})^{n_i(t)}
    (A(t)B(t)_{x_i(t)})_{x_i'(t)}\right)^{n_i'(t)}\right)>2\cosh r,
  \end{equation}
  and we have again
  \begin{equation}\label{eqnend}
    \left[AB_{x_i(t)}\cdot\left(B(AB_{x_i(t)})^{n_i(t)}\cdot
    (AB_{x_i(t)})_{x_i'(t)}
    \right)^{n_i'(t)},B(AB_{x_i(t)})^{n_i(t)}\right]=[A,B].
  \end{equation}
  The equality (\ref{eqnend2}) being strict, for all $\tau\in[0,1]$ there
  exists an interval $(\tau-\delta,\tau+\delta)$ such that for all
  $t\in(\tau-\delta,\tau+\delta)\,\cap\,[0,1]$ we have
  \[ \Tr\left(A(t)B(t)_{x_i(\tau)}\left(B(t)(A(t)B(t)_{x_i(\tau)})^{n_i(\tau)}
    (A(t)B(t)_{x_i(\tau)})_{x_i'(\tau)}\right)^{n_i'(\tau)}\right)>2\cosh r. \]
  The compact set $[0,1]$ is covered by finitely many such intervals, hence
  there exists a subdivision
  $0=t_0<t_1<\cdots<t_k=1$, and elements $n_i^j,{n_i^j}'\in\dZ$,
  $x_i^j,{x_i^j}'\in\dR$, such that for every $i\in\{1,2\}$ and
  $j\in\{0,\ldots,k-1\}$ such that for all $t\in[t_j,t_{j+1}]$ we have
  \begin{equation}\label{eqnend3}
    \Tr\left(A(t)B(t)_{x_i^j}\left(B(t)(A(t)B(t)_{x_i^j})^{n_i^j}
    (A(t)B(t)_{x_i^j})_{{x_i^j}'}\right)^{{n_i^j}'}\right)>2\cosh r,
  \end{equation}
  with
  $x_i^0={x_i^0}'=n_i^0={n_i^0}'=0$.
  We also set $x_i^k={x_i^k}'=n_i^k={n_i^k}'=0$.
  Fix $j\in\{0,\ldots,k-1\}$.
  For all $t\in[t_j,t_{j+1}]$ and $i\in\{1,2\}$, we set
  \[ f_t^j(a_i)=A_i(t)B_i(t)_{x_i^j}\cdot\left(B_i(t)
  (A_i(t)B_i(t)_{x_i^j})^{n_i^j}
  \cdot(A_i(t)B_i(t)_{x_i^j})_{{x_i^j}'}\right)^{{n_i^j}'} \]
  and
  \[ f_t^j(b_i)=B_i(t)(A_i(t)B_i(t)_{x_i^j})^{n_i^j}, \]
  with $A_i(t)=\rho_t(a_i)$ and $B_i(t)=\rho_t(b_i)$,
  and we put $f_t^j(a_i)=\rho_t(a_i)$, $f_t^j(b_i)=\rho_t(b_i)$ for
  $i\leq 3$.
  The identity (\ref{eqnend}) guarantees that we are indeed defining a
  representation $f_t^j\in R_g$, for all $j\in\{0,\ldots,k-1\}$ and
  $t\in[t_j,t_{j+1}]$.
  Note that we have $f_0^0=\rho_0$.
  Denote also $f_1^k=\rho_1$.

  Now we can construct a path joining $\rho_0$ to $\rho_1$ in $W^r$, as
  follows:
  \begin{itemize}
    \item The path $f_t^j$ enables to travel from $f_{t_j}^j$ to
      $f_{t_{j+1}}^j$ without leaving $W^r$, by (\ref{eqnend3}).
    \item In order to go from $f_{t_{j+1}}^j$ to $f_{t_{j+1}}^{j+1}$, we first
      use Dehn twists in the first handle, and then we use Dehn twists in the
      second handle.
      More precisely, for every $j$ between
      $0$ and $k-1$, denote
      $y_i(t)=tx_i^j+(1-t)x_i^{j+1}$,
      $y_i'(t)=t{x_i^j}'+(1-t){x_i^{j+1}}'$,
      $m_i(t)=tn_i^j+(1-t)n_i^{j+1}$,
      $m_i'(t)=t{n_i^j}'+(1-t){n_i^{j+1}}'$,
      et $A_i=A_i(t_{j+1})$, $B_i=B_i(t_{j+1})$. We then set
      \[ g_t^j(a_1)=A_1{B_1}_{y_1(t)}\cdot\left(B_1(A_1
      {B_1}_{y_1(t)})_{m_1(t)}
      \cdot(A_1{B_1}_{y_1(t)})_{y_1'(t)}\right)_{m_1'(t)}, \]
      \[ g_t^j(b_1)=B_1(A_1{B_1}_{y_1(t)})_{m_1(t)} \]
      and $g_t^j(x)=f_{t_{j+1}}^j(x)$ for all $t\in[0,1]$ and all
      $x\in\{a_2,b_2,\ldots,a_g,b_g\}$, and we set
      \[ h_t^j(a_2)=A_2{B_2}_{y_2(t)}\cdot\left(B_2(A_2
      {B_2}_{y_2(t)})_{m_2(t)}\cdot(A_2{B_2}_{y_2(t)})_{y_2'(t)}
      \right)_{m_2'(t)}, \]
      \[ h_t^j(b_2)=B_2(A_2{B_2}_{y_2(t)})_{m_2(t)} \]
      and $h_t^j(x)=f_{t_{j+1}}^{j+1}(x)$ for all $t\in[0,1]$ and all
      $x\in\{a_1,b_1,a_3,b_3,\ldots,a_g,b_g\}$.

      We then have $g_0^j=f_{t_{j+1}}^j$, $g_1^j=h_0^j$ and
      $h_1^j=f_{t_{j+1}}^{j+1}$, and the equality (\ref{eqnend3})
      still guarantees that these paths do not leave $W^r$.
  \end{itemize}

  \medskip
  
  {\bf Step 2:} {\em Let
  $\rho\in m_{g,k}\smallsetminus K^{r+8g\delta_{\dH^2}}$. Then there exists
  a path $\rho_t$ taking values in
  $m_{g,k}\smallsetminus K^r$, such that $\rho_0=\rho$
  and $\rho_1\in V^r$}.

  By lemma \ref{lemtechoneend3}, there exist $s_1,s_2\in\{a_1,\ldots,b_g\}$
  such that $F_{\{s_1,s_2\}}^r(\rho)\neq\emptyset$.
  If $g\geq 3$, then there exists $i\in\{1,\ldots,g\}$ such that
  $\{a_i,b_i\}\cap\{s_1,s_2\}=\emptyset$. In that case, as in step 1, we can
  use Dehn twists in the handle $i$, without entering $K^r$ since we do not
  touch $\rho(s_1)$, $\rho(s_2)$. This completes the proof of proposition
  \ref{oneendtech}, in the case $g\geq 3$.

  In the case when $g=2$ and $\{s_1,s_2\}=\{a_1,b_1\}$ or
  $\{s_1,s_2\}=\{a_2,b_2\}$, we do the same as in the preceding cases.
  Now suppose for instance that $s_1=a_1$ et $s_2=a_2$ (the other cases are
  dealt with similarly). If $\rho(a_1)$ or $\rho(a_2)$ is hyperbolic (say, for
  instance, $\rho(a_1)$), then as in lemma \ref{lemtechoneend1}, Dehn twists
  enable to get $\Tr(\rho_t(b_1))>2\cosh r$ (almost) without touching
  $\rho(a_1)$ and $\rho(a_2)$. If $\rho(a_1)$ or $\rho(a_2)$ is parabolic
  (say, $\rho(a_1)$), then
  one proves easily
  (see \textsl{e.g.} lemma 3.2)
  that $\Tr(\rho(a_1)\rho(b_1)_t)$ is non constant
  at $t=0$, and in that way we can get to the case when $\rho(a_1)$ is
  hyperbolic.

  The last case, in which $\rho(a_1),\rho(a_2)\in Ell$, is still going to
  require some work. Up to conjugating $\rho$, we have
  $\rho(a_1)=\left(\begin{array}{cc}\cos\theta & -\sin\theta \\ \sin\theta
  & \cos\theta\end{array}\right)$, and
  $\rho(b_1)=\left(\begin{array}{cc}a & b \\ c & d\end{array}\right)$.
  There exists $\alpha\in\dR$ such that the vectors
  $\left(\begin{array}{c}\cos\alpha \\ \sin\alpha\end{array}\right)$ and
  $\left(\begin{array}{c}a+d \\ b-c\end{array}\right)$ of $\dR^2$ are
  colinear. And $(a+d)^2+(b-c)^2=4+(a-d)^2+(b+c)^2>4$, since
  $\rho(b_1)$ does not commute with $\rho(a_1)$ (indeed, if $a=d$ and
  $b=-c$ then $a^2+b^2=1$ hence $\rho(b_1)$ is under the form
  $\left(\begin{array}{cc}\cos x & -\sin x \\
  \sin x & \cos x\end{array}\right)$).
  We then have
  $\Tr\left(\left(\begin{array}{cc}a & b \\ c & d\end{array}\right)\cdot
  \left(\begin{array}{cc}\cos\alpha & -\sin\alpha \\ \sin\alpha & \cos\alpha
  \end{array}\right)\right)>2$, so that, for some $t\in\dR$,
  $\rho(b_1)\rho(a_1)_t\in Hyp$. Note also that the axis of that hyperbolic
  element passes through the fixed point of $\rho(a_1)$
  (indeed, $\rho(b_1)\rho(a_1)_t$ is of the form
  $\left(\begin{array}{cc}x & y \\ y & z\end{array}\right)$. By conjugating
  this element by the rotation
  $\left(\begin{array}{cc}0 & -1 \\ 1 & 0\end{array}\right)$ we get
  $\left(\rho(b_1)\rho(a_1)_t\right)^{-1}$, so the element
  $\left(\begin{array}{cc}0 & -1 \\ 1 & 0\end{array}\right)$
  exchanges the fixed points of $\rho(b_1)\rho(a_1)_t$, which implies that the
  axis of $\rho(b_1)\rho(a_1)_t$ passes through the fixed point of
  $\left(\begin{array}{cc}0 & -1 \\ 1 & 0\end{array}\right)$).
  This is not extremely important for the proof, but explains figure
  \ref{figureend}.

  Denote $M=[\rho(a_1),\rho(b_1)]$. It is the commutator of two non-commuting
  elliptic elements in $\psl$. An easy computation then yields $\Tr(M)>2$,
  hence $M\in Hyp$. The number $d_{\{a_1,a_2\}}(\rho)$ depends only on the
  angles of the rotations $\rho(a_1)$ and $\rho(a_2)$, and of the distance
  between their respective centres $x_1$, $x_2\in\dH^2$. Regardless of the
  repartition of these fixed points and the axis of $M$, there exists
  $\varepsilon\in\{-1,1\}$ such that for all $t>0$,
  $d(x_1,M_{\varepsilon t}\cdot x_2)>d(x_1,x_2)$.
  Denote by $p\in\partial\dH^2$ the attractive point of $M_\varepsilon$.
  By conjugating $\rho(a_2)$ and $\rho(b_2)$ by $M_{\varepsilon t}$
  (equivalently, by making a Dehn twists along the simple closed curve
  representing $[a_1,b_1]$), we get a representation $\rho_t$ such that
  $d_{\{a_1,a_2\}}(\rho_t)$ is arbitrarily large; the fixed point of
  $\rho_t(a_2)$ being sent to $p$.

  Now suppose that $\rho(b_1)$ is hyperbolic, with its axis passing through the
  fixed point of $\rho(a_1)$. Then for all $t$ such that $|t|$ is small enough,
  the element $\rho(a_1)\rho(b_1)_t$ is elliptic; its fixed point can be found
  as suggests the following picture:
  \begin{figure}[h]
    \[\epsfbox{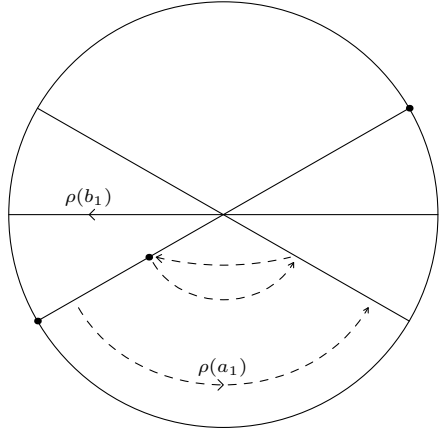}\]
    \caption{The element $\rho(a_1)\cdot\rho(b_1)_t$.}
    \label{figureend}
  \end{figure}

  There exist two real numbers $t_1$ and $t_2$ such that
  $\rho(a_1)\rho(b_1)_{t_i}$ are parabolic elements, whose unique fixed point
  is one of the points drawn on the above picture. At least one of these points
  does not depend on $p$; denote by $T=t_i$ a corresponding element.
  Then for all $t\in[0,T]$, $p$ is not a fixed point of
  $\rho(a_1)\rho(b_1)_t$. Hence, if we have sent the fixed point of
  $\rho(a_2)$ far enough, this defines a path $\rho_t'$ such that
  $d_{\{a_1,a_2\}}(\rho_t')$ stays large, and for all $t$ big enough,
  $\rho_t'(a_1)\in Par$, hence we recover the preceding situation.
\end{preuvede}

Similarly, we can prove that every loop can be pushed out of every compact set.
In other words, the fundamental group of the space $e^{-1}(k)$ is entirely
carried by this only end.

\bigskip

This fact will be useful here for the following reason:

\begin{cor}
  Let $k\in\{2-2g,\ldots,2g-2\}$ and denote $X=m_{g,k}$.
  Let $(\overline{X},i)$ be a compactification of $X$.
  Then $\partial\overline{X}$ is connected.
\end{cor}

\begin{preuve}
  The space $X$ has at least one end, hence is non compact: therefore,
  $\partial\overline{X}$ is non empty. Suppose that
  $\partial\overline{X}=A\cup B$, where $A$ and $B$ are disjoint, open and
  closed subsets of $\partial\overline{X}$.
  The boundary $\partial\overline{X}$ being closed, $A$ and $B$ are closed
  subsets of $\overline{X}$.
  Since $\overline{X}$ is compact Hausdorff, it is normal. Hence, there exist
  two open subsets $U$ and $V$ of $\overline{X}$ such that $A\subset U$,
  $B\subset V$, and $U\cap V=\emptyset$. The open set $U\cup V$ contains
  $\partial\overline{X}$ hence it is the complement of a compact subset $K$ of
  $X$. This compact set $K$ is contained in $K^r$, for $r$ large enough. And
  $X\smallsetminus K^r$ possesses a unique unbounded connected component,
  denote it by $W$. By proposition \ref{oneendtech}, $X\smallsetminus W$ is
  contained in $K^{r+8g\delta_{\dH^2}}$, thus
  $K'=X\smallsetminus W$ is a compact subset of $X$.
  Now $W$ is dense in $W'=\overline{X}\smallsetminus K'$, which is therefore
  connected. Denote $U'=W'\cap U$ and $V'=W'\cap V$. Then
  $U'$ contains $A$, $V'$ contains $B$, $U'$ and $V'$ are open,
  $U'\cap V'=\emptyset$, and $W'=U'\cup V'$ is connected. This implies that
  $U'=\emptyset$ or $V'=\emptyset$, hence $A=\emptyset$ or $B=\emptyset$.
\end{preuve}

\subsection{Oriented compactification}

In most of the sequel we will consider the case $n=2$, and we consider the
compactification of the space $m_\Gamma^u$ of actions on the hyperbolic plane
$\dH^2$. We will prove, at least in the case when $\Gamma$ is a surface group,
that this space is degenerated; and it seems more natural to study a
compactification of the oriented version $m_\Gamma^o$.

As we said before, the ideal points of the compactification $\overline{m_g^u}$
of M. Bestvina and F. Paulin are (equivariant isometry classes of) actions of
$\pi_1\Sigma_g$ by isometries on $\dR$-trees. We shall prove that it is
possible to equip these $\dR$-trees with an orientation, so that we can define
an Euler class on these trees. This will enable to define a compactification
of $m_g^o$, in which the Euler class extends continuously to the boundary.

\subsubsection{Fat $\dR$-trees}

Let $X$ be a hyperbolic space in the sense of Gromov. A {\em germ of rays}, in
$X$, is an equivalence class of rays, for the following equivalence relation:
we say that two rays are equivalent if they coincide on some non trivial
initial segment.

In the rest of this section, $T$ will be an $\dR$-tree non reduced to a point.
At every point $x\in T$, denote by $\mathcal{G}(x)$ the set of germs of rays
issued of $x$. An {\em orientation of $T$} is the data, for all $x\in T$, of a
total cyclic order $or(x)$ in $\mathcal{G}(x)$.

\begin{deuf}\label{arbrelepais}{\ }

  \begin{itemize}
    \item An $\dR$-tree equipped with an orientation is called a
      {\em fat $\dR$-tree}.
    \item Let $(T,or)$ and $(T',or')$ be two fat $\dR$-trees and let
      $h\in Isom(T,T')$ be an isometry.
      Of course, $h$ defines, at every point $x\in T$, a bijection
      $\mathcal{G}h_x\colon\mathcal{G}(x)\rightarrow\mathcal{G}(h(x))$.
      We say that {\em $h$ preserves the orientation} if for every
      $x\in T$, the following diagram commutes:
      \[ \xymatrix{\mathcal{G}(x)^3\ar[r]^{or(x)}\ar[d]_{\left(
      \mathcal{G}h_x\right)^3} & \{-1,0,1\} \, .\\ \mathcal{G}
      (h(x))^3\ar[ru]_{or'(h(x))}} \]
      The set of isometries of $T$ which preserve the orientation $or$ form a
      subgroup of $Isom(T)$, which will denote $Isom^{or}(T)$.
    \item We say that an action $\rho\colon\Gamma\rightarrow Isom(X)$ of a
      group $\Gamma$ by isometries {\em preserves the orientation $or$} if it
      takes values into $Isom^{or}(T)$.
  \end{itemize}
\end{deuf}

Since we are interested in defining the Euler class of actions on $\dR$-trees
preserving the orientation, we need to consider more particularly the boundary
of the tree. We say that a total cyclic order $o$ on $\partial_\infty T$ is
{\em coherent} if for every $x\in T$ and every non degenerate triple
$\{c_1,c_2,c_3\}$ of rays issued from $x$, the element $o(r_1,r_2,r_3)$ does
not depend on the chosen representants $r_1$, $r_2$, $r_3$ of $c_1$, $c_2$,
$c_3$.

For instance, in the following configuration
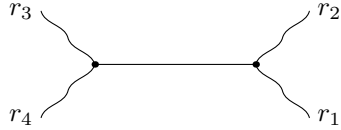
\[\epsfbox{Figures/FiguresExtraComp.1}\]
a total cyclic order $o$ on the boundary $\{r_1,r_2,r_3,r_4\}$ is coherent if
and only if it verifies $o(r_1,r_2,r_3)=o(r_1,r_2,r_4)$ et
$o(r_1,r_3,r_4)=o(r_2,r_3,r_4)$.

\bigskip

Of course, there is a relation between orientations on an $\dR$-tree and
coherent total cyclic orders on its boundary. Let us begin with some notations.

Naturally, we say that an $\dR$-tree $T$ is the
{\em convex hull of its boundary} if it is the closure of the union of the
geodesics which are isometric to the whole
$\dR$. In other words, if
\[ T=\overline{\mathop{\bigcup_{a,b\in\partial_\infty T}}_{a\neq b}(a,b)}. \]

We also say that two non degenerate tripods are equivalent if their
intersection is a non degenerate tripod; we denote by $Trip(T)$ the set of
equivalence classes.
Denote
\[O_T'
=\{-1,1\}^{Trip(T)},\]
and equip it with the product topology (the topology of pointwise convergence).
By Tikhonov's theorem, it is a compact space. Denote by $O_T$ the (closed) set
of functions $o$ which satisfy the following conditions: for every tripod
$Trip(a_0,a_1,a_2,a_3)$ (where $a_0$ is the central point),
we have
\[o(Trip(a_0,a_1,a_2,a_3))=o(Trip(a_0,a_2,a_3,a_1))=
-o(Trip(a_0,a_1,a_3,a_2)).\]
Non degenerate triples of germs of rays issued from a same point are naturally
identified with classes of tripods, so that an orientation of $T$ is naturally
identified with an element of $O_T$.
More precisely, denote by $OR_T\subset O_T$ the set consisting of functions
satisfying the following condition: for all $a\in T$ and all tripods
$Trip(a,b,c,d)$ and $Trip(a,b,d,e)$ of central point $a$,
\[or(Trip(a,b,c,d))=or(Trip(a,b,d,e))=1\Rightarrow or(Trip(a,b,c,e))=1.\]
Then $OR_T$ identifies to the set of orientations of $T$; remark that $OR_T$ is
a closed set.

Now let $E_T$ be the set
\[E_T
=\{-1,0,1\}^{(\partial T)^3},\]
equipped with the product topology. Let $E_T'$ be the (closed) subset of $E_T$,
consisting of elements $o\in E_T$ satisfying:
\[ \forall (x,y,z)\in (\partial T)^3,\ o(x,y,z)=0\Leftrightarrow
Card\{x,y,z\}\leq 2. \]
Denote by $ORD_T$ the set of coherent cyclic orders on $\partial_\infty T$.

\medskip

For every non degenerate triple $(x,y,z)\in(\partial_\infty T)^3$,
denote by $t(x,y,z)\in Trip(T)$ the class of the tripod
$(x,y)\,\cup\,(x,z)\,\cup\,(y,z)\,$. This function $t$ defines, by
precomposition, a function $f\colon O_T\rightarrow E_T$.

\begin{prop}
  Suppose that $T$ is the convex hull of its boundary. Then the restriction
  of $f$ to $OR_T$ is a homeomorphism on its image $ORD_T$.
\end{prop}

\begin{preuve}
  Denote
  \[T'=\mathop{\bigcup_{a,b\in\partial_\infty T}}_{a\neq b}(a,b).\]
  We easily check that $T'$ is connected (if $x,y\in T'$, then $x\in \,(a,b)$
  and $y\in\,(c,d)$: then $(a,c)\,\cap\,(a,b)$ is non empty since $T$ is a
  tree, thus there exists a path in $T'$ joining $x$ to a point of $(a,c)\,$,
  and $(a,c)\,\cap\,(c,d)\,\neq\emptyset$ guarantees that there is a path in
  $T'$ joining $x$ to a point of $(c,d)\,$, and then to $y$),
  hence $T'$ is a subtree of $T$, dense in $T'$. Therefore, every class of
  tripods possesses a representant in $T'$; one checks also that every segment
  in $T'$ is a subset of some line $(a,b)\,$. The function $t$ is thus
  surjective, so that $f$ is injective.
  
  The continuity of $f$ is straightforward, and hence $f$ is a homeomorphism
  from $OR_T$ (which is compact) on its image (which is Hausdorff); now we
  need to check that $f(OR_T)=ORD_T$.

  Let $or$ be an orientation of $T$, and $o=f(or)$. It is immediate that the
  function $o\colon\partial_\infty T^3\rightarrow\{-1,0,1\}$ satisfies the
  conditions $1$ and $2$ of definition \ref{defordrecyclique}.
  Let us check that condition $3$ is also satisfied: suppose that
  $o(r_1,r_2,r_3)=o(r_1,r_3,r_4)=1$. For every non degenerate triple
  $(a,b,c)\in\partial_\infty T$, $\,(a,b)\,\cap\,(a,c)\,\cap\,(b,c)$
  is a singleton; denote it by $\{P_{abc}\}$. If
  $P_{r_1r_3r_4}\in\,(r_1,P_{r_1r_2r_3})$, then
  
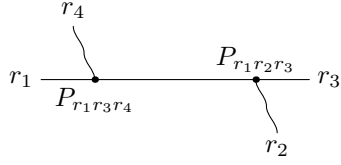
\[ \epsfbox{Figures/FiguresExtraComp.2} \]
  $P_{r_1r_2r_4}=P_{r_1r_3r_4}$ and the germs of
  $[P_{r_1r_2r_4},r_2)\,$ and $[P_{r_1r_2r_4},r_3)\,$ are identical,
  so that $o(r_1,r_2,r_4)=1$. If $P_{r_1r_3r_4}\in\,(P_{r_1r_2r_3},r_3)\,$
  the argument is the same. If $P_{r_1r_3r_4}=P_{r_1r_2r_3}$, then the germs of
  the rays $[P_{r_1r_2r_3},r_2)\,$ and $[P_{r_1r_2r_3},r_4)\,$ are distinct
  since we have $o(r_1,r_3,r_2)\neq o(r_1,r_3,r_4)$: hence
  $[P_{r_1r_2r_3},r_1)\,$, \ldots, $[P_{r_1r_2r_3},r_4)\,$ define four germs of
  rays issued from $P_{r_1r_2r_3}$, which are all distinct, hence
  $o(r_1,r_2,r_4)=1$, since $or\left(P_{r_1r_2r_3}\right)$ satisfies condition
  $3$ of definition \ref{defordrecyclique}. And by construction, the order $o$
  is coherent.
  We have just proved that $f(OR_T)\subset ORD_T$.
  But by definition, a total cyclic order is coherent if and only if it is in
  the image of $f$, and if $f(or)$ is a coherent total cyclic order we check
  easily that $or$ is indeed an orientation of $T$.
\end{preuve}

Now recall that an action of a group $\Gamma$ by isometries on an $\dR$-tree
$T$ is called {\em minimal} if $T$ possesses no subtree
$T'\subset T$, invariant under the action of $\Gamma$, distinct from
$\emptyset$ and $T$.

In the sequel, $\Gamma$ is a finitely generated group and we consider minimal
actions of $\Gamma$ on $\dR$-trees. In that case, $T$ is
the union of the translation axes of the hyperbolic elements in the image of
$\Gamma$ (see \textsl{e.g.} \cite{MorganShalen,Paulin89}).
In particular, such trees are the convex hulls of their boundaries.

\bigskip

In the sequel, we shall need to consider the {\em set} of classes of minimal
actions of $\Gamma$ on fat $\dR$-trees, preserving the orientation, up to
equivariant isometry preserving the order. General arguments of cardinality
enable to do that,
but the following proposition enables us to see this as an explicit set.

If $(T,or)$ is a fat $\dR$-tree and if $u,v,w\in T$ are not aligned, then they
define a class of tripods and we denote by $or(u,v,w)\in\{-1,1\}$ the image of
this tripod by $or$. If $u,v,w$ are aligned, then we write $or(u,v,w)=0$.

\begin{prop}
  Let $(T,or)$ be a fat $\dR$-tree,
  let $x_0\in T$, and let $\rho\colon\Gamma\rightarrow Isom^{or}(T)$ be a
  minimal action of a finitely generated group $\Gamma$, preserving the
  orientation. Then the $\dR$-tree $T$, the orientation $or$ and the action
  $\rho$ are entirely determined by the functions
  $f\colon\Gamma^2\rightarrow\dR$ et $g\colon\Gamma^3\rightarrow\{-1,0,1\}$
  defined by
  \[f(\gamma_1,\gamma_2)=d_T(\gamma_1 x_0,\gamma_2 x_0)\hspace{0.4cm}\text{et}
  \hspace{0.4cm}
  g(\gamma_1,\gamma_2,\gamma_3)=o(\gamma_1 x_0,\gamma_2 x_0,\gamma_3 x_0).\]
  More precisely, if $\rho\colon\Gamma\rightarrow Isom^{or}(T)$ and
  $\rho'\colon\Gamma\rightarrow Isom^{or'}(T')$ define the same functions $f$
  and $g$, for some choices of base points in $T$ and $T'$, then there exists
  an equivariant isometry $\varphi\colon T\rightarrow T'$ preserving the order
  and the base point.
\end{prop}

\begin{preuve}
  This is well-known (see \textsl{e.g.} \cite{Paulin89}) for minimal actions of
  a finitely generated groups on $\dR$-trees; here we simply need to add the
  orientation.

  Let us begin with the following remark. Let $T$ and $T'$ be two fat
  $\dR$-trees and $x_1,\ldots,x_n\in T$, $x_1',\ldots,x_n'\in T'$ such that for
  for all $i,j$, $d(x_i,x_j)=d(x_i',x_j')$ and all $i,j,k$,
  $or(x_i,x_j,x_k)=or'(x_i',x_j',x_k')$. Denote $K=\{x_1,\ldots,x_n\}$.
  Then the function
  $\varphi_K\colon\{x_1,\ldots,x_n\}\rightarrow\{x_1',\ldots,x_n'\}$ defined by
  $\varphi_K(x_i)=x_i'$  extends uniquely to an orientation-preserving isometry
  \[\varphi_{Hull(K)}\colon Hull(\{x_1,\ldots,x_n\})\rightarrow
  Hull(\{x_1',\ldots,x_n'\})\]
  (the subsets $Hull(\{x_1,\ldots,x_n\})$ and $Hull(\{x_1',\ldots,x_n'\})$, as
  subtrees of $T$ and $T'$, are oriented trees, by the restrictions of $or$
  and $or'$).

  We prove this by induction. If $n=1$, there is not very much to do.
  Denote again $K=\{x_1,\ldots,x_n\}$ and suppose that
  $\varphi_{Hull(K)}$ is an orientation-preserving isometry between
  $Hull(\{x_1,\ldots,x_n\})$ and $Hull(\{x_1',\ldots,x_n'\})$. Denote by
  $y_{n+1}$ the projection of $x_{n+1}$ on $Hull(K)$. The following relation,
  true in every $\dR$-tree,
  \[ \epsfbox{Figures/FiguresExtraComp.3} \]
  enables to find $y_{n+1}$ in the tree $Hull(K)$: the real numbers
  $d(x_1,x_{n+1})$, \ldots, $d(x_n,x_{n+1})$ determine a unique point
  $y_{n+1}\in Hull(K)$; and similarly they determine a unique point
  $y_{n+1}'\in Hull(\{x_1',\ldots,x_n'\})$, and we have
  $\varphi_{Hull(K)}(y_{n+1})=y_{n+1}'$.
  If $x_{n+1}=y_{n+1}$ then
  $Hull(\{x_1,\ldots,x_n\})=Hull(\{x_1,\ldots,x_{n+1}\})$ and
  $Hull(\{x_1',\ldots,x_n'\})=Hull(\{x_1',\ldots,x_{n+1}'\})$ and the
  induction is proved. Otherwise, $Hull(\{x_1,\ldots,x_{n+1}\})$ is obtained
  by glueing at the point $y_{n+1}$ the tree $Hull(\{x_1,\ldots,x_n\})$ and
  the segment $[x_{n+1},y_{n+1}]$, whose length is determined by the real
  numbers $d(x_1,x_{n+1})$, \ldots, $d(x_n,x_{n+1})$; and similarly for
  $Hull(\{x_1',\ldots,x_{n+1}'\})$ in $T'$. In that way, the isometry
  $\varphi_{Hull(K)}$ extends to a unique isometry
  $\varphi_{Hull(\{x_1,\ldots,x_{n+1}\})}$. Since $\varphi_{Hull(K)}$
  preserves the orientation, we need only check that
  $\varphi_{Hull(\{x_1,\ldots,x_{n+1}\})}$ preserves the orientation at the
  vertex $y_{n+1}$. But this follows from the fact that all the classes of the
  tripods of centre $y_{n+1}$ have a representant of the type
  $Hull(\{x_i,x_j,x_{n+1}\})$, where $Card(\{i,j,n+1\})=3$.


  Now suppose that $(\rho,T)$ and $(\rho',T')$ define the same functions
  $f$ and $g$; denote by $x_0$ and $x_0'$ the base points.
  Let $P$ be a finite subset of $\Gamma$. We then have a unique isometry
  $\varphi_{Hull(P\cdot x_0)}$ between $Hull(P\cdot x_0)$ and
  $Hull(P\cdot x_0')$, preserving the orientation, such that
  $\varphi_{Hull(P\cdot x_0)}(\gamma\cdot x_0)=\gamma x_0'$.
  In particular, for every finite subset $Q$ of $P$ the restriction of
  $\varphi_{Hull(P\cdot x_0)}$ to $Hull(Q\cdot x_0)$ equals
  $\varphi_{Hull(Q\cdot x_0)}$, so that we can construct an isometry
  $\varphi\colon\bigcup_{P\subset\Gamma}Hull(P\cdot x_0)\rightarrow
  \bigcup_{P\subset\Gamma}Hull(P\cdot x_0')$, such that for all
  $\gamma\in\Gamma$ we have $\varphi(\gamma\cdot x_0)=\gamma\cdot x_0'$; and
  and we also deduce that for every
  $\gamma\in\Gamma$ and every finite subset $P$ of $\Gamma$ we have
  \[  \forall  y\in Hull(P\cdot x_0),\ \ \ 
    \rho'(\gamma)\cdot\varphi_{Hull(P\cdot x_0)}=
    \varphi_{Hull(\gamma P\cdot x_0)}(\rho(\gamma)\cdot y)  \]
  which ensures that the isometry
  $\varphi$ is equivariant for the actions $\rho$, $\rho'$ on the trees
  $\bigcup_{P\subset\Gamma}Hull(P\cdot x_0)$ and
  $\bigcup_{P\subset\Gamma}Hull(P\cdot x_0')$.
  Every tripod in $\bigcup_{P\subset\Gamma}Hull(P\cdot x_0)$ is in
  $Hull(P\cdot x_0)$ for $P$ big enough, hence $\varphi$ preserves the
  orientation.
  Finally, the actions $\rho$ and $\rho'$ being minimal, we have
  $\bigcup_{P\subset\Gamma}Hull(P\cdot x_0)=T$ and
  $\bigcup_{P\subset\Gamma}Hull(P\cdot x_0')=T'$; whence
  $\varphi\colon T\rightarrow T'$ is an orientation-preserving equivariant
  isometry, such that $\varphi(x_0)=x_0'$.
\end{preuve}

Of course, $f$ and $g$ depend on $x_0$, anyway it follows from this proposition
that the classes of actions of $\Gamma$ on fat $\dR$-trees not reduced to a
point, up to orientation-preserving equivariant isometry, form a set, which we
denote $\mathcal{T}''(\Gamma)$.
We denote by $\mathcal{T}^o(\Gamma)$ its subset formed by those
$(\rho,T)\in\mathcal{T}''(\Gamma)$
such that
$\displaystyle{\min_{x_0\in T}\max_{\gamma\in S}}\,d_T(x_0,\gamma\cdot x_0)=1$
and such that whenever $\rho$ possesses at least one global fixed point in
$\partial_\infty T$, the tree $T$ is isometric to $\dR$.

Our aim now is to define a topology on the set
$m_{\Gamma}^o(2)\cup\mathcal{T}^o(\Gamma)$.

\subsubsection{Rigidity of the order}

We first need to give some technical lemmas indicating that the orders given by
triples of points in fat $\dR$-trees, as well as in $\dH^2$, are stable under
small perturbations of the tree or of the plane.
In all this section, $X$ will be a fat $\dR$-tree or the hyperbolic plane
$\dH^2$, equipped with its orientation (and hence, with a total cyclic order
on its boundary), and with a metric $\dfrac{d_{\dH^2}}{d}$ proportional to its
usual metric.
Its best hyperbolic constant is then $\delta(X)=\dfrac{\delta_{\dH^2}}{d}$.

\begin{lem}\label{Toricelli}
  Let $x_1, x_2, x_3\in X$. Then there exists a unique $x_0\in X$ which
  minimizes the function $x\mapsto d(x,x_1)+d(x,x_2)+d(x,x_3)$. Moreover, the
  function $X^3\rightarrow X$ defined by $(x_1,x_2,x_3)\mapsto x_0$ is
  continuous.
\end{lem}

\begin{preuve}
  If $X$ is an $\dR$-tree, then we check easily that the unique point $m\in X$
  such that $[x_1,x_2]\cap[x_1,x_3]=[x_1,m]$ is the point $x_0$ wanted.
  If $X$ is the hyperbolic plane $\dH^2$ equipped with a proportional metric
  to $d_{\dH^2}$, then the function $x\mapsto d(x,x_1)+d(x,x_2)+d(x,x_3)$ is
  convex, proper, hence achieves a minimum. And the $CAT(0)$ inequality
  implies that this function cannot be constant on any non degenerate segment,
  hence this minimum is unique. Moreover, this convex function depends
  continuously on $x_1$, $x_2$ and $x_3$, hence its unique minimum also depends
  continuously on $x_1$, $x_2$ and $x_3$.
\end{preuve}

\begin{ps}
  In the Euclidean plane $\dR^2$, the point $x_0$ is called the
  {\em Fermat point} of the triangle $\Delta(x_1,x_2,x_3)$. If this triangle
  has angles smaller than $\dfrac{2\pi}{3}$, then this point coincides with
  the {\em Torricelli point}, which, in that case, sees every edge of the
  triangle under an angle equal to $\dfrac{2\pi}{3}$
  (see \textsl{e.g.} \cite{Fresnel}).
\end{ps}

Let $A\geq 0$. We denote by $V(A)\subset X^3$ the set of $(x_1,x_2,x_3)\in X^3$
such that for every permutation $(i,j,k)$ of $(1,2,3)$, we have
$d(x_i,x_j)+d(x_j,x_k)-d(x_i,x_k)>2A$.

\begin{lem}\label{FermatDedans}
  There exists a universal constant $C_0\geq 0$ such that for all
  $(x_1,x_2,x_3)\in V(C_0 \delta(X))$, we have $x_0\not\in\{x_1,x_2,x_3\}$,
  where $x_0$ is as in the lemma \ref{Toricelli}.
\end{lem}

\begin{preuve}
  The proof in the case when $X$ is an $\dR$-tree follows from the same
  observation as in the proof of lemma \ref{Toricelli}.
  In that case, any $C_0\geq 0$ satisfies the lemma.
  Now, suppose that $X=(\dH^2,d_{\dH^2})$. First notice that if all the angles
  of a non degenerate triangle $\Delta(x_1,x_2,x_3)$ are strictly lower than
  $\dfrac{\pi}{2}$, then $x_0\not\in\{x_1,x_2,x_3\}$.
  Indeed, suppose that the angle at $x_1$ is smaller than $\dfrac{\pi}{2}$.
  Then, for some $x_3'\in[x_1,x_3]$, such that $x_3'\neq x_1$, the orthogonal
  projection $p(x_3')$ of $x_3'$ on the geodesic $(x_1,x_2)$ is strictly
  between $x_1$ and $x_2$. Then
  $d(x_1,p(x_3'))+d(x_2,p(x_3'))+d(x_3,p(x_3'))<d(x_1,x_2)+d(x_1,x_3)$,
  so that $x_0\neq x_1$. Similarly $x_0\not\in\{x_2,x_3\}$.
  Now,
  we shall prove that the condition
  $(x_1,x_2,x_3)\in V(A)$, for $A$ large enough, implies that
  all the angles of the triangle $\Delta(x_1,x_2,x_3)$ are lower than
  $\dfrac{\pi}{2}$.
  Consider Poincar\'e's disk model of $\dH^2$, and let $a$ be the centre of the
  disk. Let $u_1,u_2$ be two points of $\partial\dH^2$,
  such that the measure of the angle $\widehat{u_1 au_2}$ is equal to
  $\dfrac{\pi}{2}$. Let $A$ be the hyperbolic distance between $a$ and the
  geodesic $(u_1u_2)$ (explicit computation yields $A=\frac{1}{2}\ln 3$, see
  \textsl{e.g.} \cite{Papadopoulos}, lemme 4.1).
  Now suppose that in the triangle
  $\Delta(x_1,x_2,x_3)$,
  the angle at $x_1$ is greater or equal to $\dfrac{\pi}{2}$. Then
  $d(x_1,[x_2,x_3])\leq A$, and this implies that
  $d(x_1,x_2)+d(x_1,x_3)-d(x_2,x_3)\leq 2A$, so that
  $(x_1,x_2,x_3)\not\in V(A)$. It follows that
  $C_0=\dfrac{A}{\delta_{\dH^2}}=\frac{\ln 3}{2\ln(1+\sqrt{2})}$
  satisfies the lemma, in the case when $X=(\dH^2,d_{\dH^2})$.
  Finally, in the case when $X$ is the hyperbolic plane equipped with a
  proportional metric to $\dfrac{d_{\dH^2}}{d}$, the preceding arguments prove
  the lemma, with the same constant $C_0$.
\end{preuve}

In the sequel, $C_0$ will be the constant that has been explicited in the proof
of lemma \ref{FermatDedans}.
By choosing another equivalent definition of the $\delta$-hyperbolicity in the
sense of Gromov
and by considering, for instance, the real number
$A=\frac{1}{2}\ln 3$ as the hyperbolicity constant of
$\dH^2$, we would have got $C_0=1$. Therefore, the value of $C_0$ does not have
any deep meaning here.

Now define a set $U\subset X^6$, as follows: we say that
$(x_1,x_2,x_3,y_1,y_2,y_3)\not\in U$ if there exist
$i,j\in\{1,2,3\}$, $i\neq j$, and rays $r_i$, $r_j$ issued from
$x_i$, $x_j$ and passing through $y_i$, $y_j$ respectively, such that
$r_i$, $r_j$ represent the same point of $\partial_\infty X$, or end at the
same end point of $T$, in the case when $T$ is a tree
(heuristically, $(x_1,\ldots,y_3)\in U$ if the oriented segments
$[x_i,y_i]$ point three distinct directions). This implies, in particular,
that $x_i\neq y_i$, for all $i$.
If $X=\dH^2$, then there is a unique ray (up to parametrization)
issued from $x_i$ and passing through $y_i$; the condition
$(x_1,\ldots,y_3)\in U$ expresses the fact that the ends of these three rays
are three distinct points of $\partial\dH^2$. In the case of an $\dR$-tree, we
can also give the following equivalent condition:

\begin{lem}\label{lemarbreextrapau}
  Let $X$ be an $\dR$-tree. Then $(x_1,\ldots,y_3)\in U$ if and only if
  for every $i\in\{1,2,3\}$, the points $x_i$, $y_{i+1}$ and $y_{i+2}$
  (we are using here a cyclic notation for the indices)
  are in the same connected component of $X\smallsetminus\{y_i\}$.
\end{lem}

\begin{preuve}
  We first check that $(x_1,\ldots,y_3)\in U$ implies that
  $y_1$, $y_2$ and $y_3$ are not aligned. Suppose that $y_1$, $y_2$ and
  $y_3$ are three points pairwise distinct, and are aligned (the case when some
  of them coincide is treated similarly).
  For instance, take $y_2\in [y_1,y_3]$. If $x_1,x_3$ are in the same connected
  component of $X\smallsetminus\{y_1,y_3\}$ as $y_2$, and if $x_2$ is not in
  the same connected component as $y_1$ then there exist rays $r_1$, $r_2$
  issued from $x_1$, $x_2$ and passing through $y_1$, $y_2$ respectively, such
  that $r_2$ passes through $y_1$, and such that $r_1$, $r_2$ are identical
  after their passage at $y_1$. The point $x_2$ cannot be in several connected
  components at the same time so in that case we have
  $(x_1,\ldots,y_3)\not\in U$. Suppose then that $x_1$ and $y_2$ are in two
  distinct components of $X\smallsetminus\{y_1\}$.
  Then in order to have $(x_1,\ldots,y_3)\in U$, the point $x_3$ cannot be in
  the same connected component of $X\smallsetminus \{y_3\}$ as $y_1$. As in the
  preceding case, we cannot place the point $x_2$, once again, so that
  $(x_1,\ldots,y_3)$ is in $U$.

  In order to have $(x_1,\ldots,y_3)\in U$ it is necessary that
  $y_1$, $y_2$, $y_3$ form a non degenerate tripod.
  If, for instance, $x_1$ and $y_2$, $y_3$ are in two distinct connected
  components of $X\smallsetminus\{y_1\}$, then: either $x_2$ is in the same
  component as $X\smallsetminus\{y_2\}$ as $y_1$, $y_3$ and in that case we
  construct rays $r_1$, $r_2$ both passing though $y_2$ and identical after
  their passage at $y_2$; either $x_2$ is in another connected component of
  $X\smallsetminus\{y_2\}$ as $y_1$, $y_3$ and then we construct rays $r_1$,
  $r_2$ which are identical after their passage at $y_3$. The only remaining
  case is when for every $i$, $x_i$ is in the same connected component of
  $X\smallsetminus\{y_i\}$ as $y_1,y_2$ and in that case, regardless to the
  three rays $r_i$ issued from $x_i$ and passing through $y_i$, the three rays
  $r_1$, $r_2$ and $r_3$ leave the tripod at the respective vertices $y_1$,
  $y_2$, $y_3$ and cannot go towards a same end (or end point) of
  $\partial_\infty X$.
\end{preuve}

\begin{lem}
  If $(x_1,\ldots,y_3)\in U$ and if $r_i,r_i'$ are rays issued from $x_i$
  and passing through $y_i$, then $o(r_1,r_2,r_3)=o(r_1',r_2',r_3')$.
\end{lem}

We write $o(x_1,\ldots,y_3)=o(r_1,r_2,r_3)$ in that case.

\begin{preuve}
  We have $x_i\neq y_i$ so in the case when $X=\dH^2$, there exists a unique
  ray $r_i$ issued from $x_i$ and passing through $y_i$.
  Now suppose that $X$ is a fat $\dR$-tree. By lemma \ref{lemarbreextrapau},
  $y_1$, $y_2$ and $y_3$ are not aligned;
  let then $y_0$ be such that $[y_1,y_2]\cap[y_1,y_3]=[y_1,y_0]$. Then
  $y_1,y_2,y_3$ define three distinct germs of rays issued from $y_0$.
  Still by lemma \ref{lemarbreextrapau}, the condition $(x_1,\ldots,y_3)\in U$
  implies that for every $i$, $x_i$ lies in the connected component of
  $X\smallsetminus\{y_i\}$ containing $y_0$.
  In particular, every ray $r$ issued from $x_i$ and passing through $y_i$
  defines a unique ray issued from $y_0$ and passing through $y_i$: it is the
  ray joining $y_0$ to $y_i$, and which then continues as the ray $r$.
  Therefore, the equality $o(r_1,r_2,r_3)=o(r_1',r_2',r_3')$ follows from the
  coherence condition on the order $o$ of the fat $\dR$-tree we are
  considering.
\end{preuve}

\begin{lem}\label{continu}
  The function $o\colon U\rightarrow\{-1,1\}$ thereby defined is continuous.
\end{lem}

\begin{preuve}
  If $X=\dH^2$, it is immediate that the function
  $X^2\smallsetminus\Delta\rightarrow\partial X$ (where $\Delta$ is the
  diagonal) sending $(x,y)$ on the end of the ray issued from $x$ and passing
  through $y$ is continuous (where $\partial X=\dS^1$ is equipped with the
  usual topology), and hence $o\colon U\rightarrow\{-1,1\}$ is simply the
  composition of two continuous functions.
  In the case when $X$ is a fat $\dR$-tree, the proof is similar to that of the
  preceding lemma.
\end{preuve}

\begin{lem}\label{OrdreConvexe}
  Let $x_1,\ldots,y_3$ such that for every $i\in\{1,2,3\}$, we have
  $\sum_j d(x_i,y_j)<\sum_j d(y_i,y_j)$.
  Then $(x_1,\ldots,y_3)\in U$.
\end{lem}

\begin{preuve}
  First suppose that $X$ is an $\dR$-tree.
  By lemma \ref{lemarbreextrapau}, if $(x_1,\ldots,y_3)\not\in U$, then for
  some $i$, either $x_i$ is not in the same component of
  $X\smallsetminus\{y_i\}$ as $y_{i+1}$ and $y_{i+2}$, and in that case
  $\sum_j d(x_i,y_j)=\sum_j d(y_i,y_j)+3d(x_i,y_i)>\sum_j d(y_i,y_j)$,
  either $y_i$ lies in the segment $[y_{i+1},y_{i+2}]$ and then $y_i$ realizes
  the minimum, in $X$, of the function
  $x\mapsto d(x,y_1)+d(x,y_2)+d(x,y_3)$, so that
  $\sum_j d(x_i,y_j)\geq\sum_j d(y_i,y_j)$.

  Now consider the case when $X$ is the hyperbolic plane, and suppose that
  $x_i\neq y_i$ for every $i$, and that $(x_1,\ldots,y_3)\not\in U$.
  Suppose, for instance, that the rays $[x_1,y_1)$ and $[x_2,y_2)$ have the
  same end. Suppose also, by contradiction, that
  $\sum_j d(x_i,y_j)<\sum_j d(y_i,y_j)$.
  Then $d(x_1,y_1)+d(x_1,y_2)+d(x_1,y_3)<d(y_1,y_2)+d(y_1,y_3)$ hence
  $d(x_1,y_2)<d(y_1,y_2)$, and, similarly, $d(x_2,y_1)<d(y_1,y_2)$.
  Since these inequalities are strict, we can move slightly the point $x_2$
  so that the rays $[x_1,y_1)$ and $[x_2,y_2)$ intersect at some point
  $a\in\dH^2$, without losing these inequalities.
  By putting $a$ at the centre of Poincar\'e's disk model of $\dH^2$, we see
  that the condition $d(x_1,y_2)<d(y_1,y_2)$ implies $d(a,y_1)<d(a,y_2)$.
  Indeed, $y_2$ is closer to $x_1$ than to $y_1$, hence $y_2$ lies in the
  hyperbolic half-plane determined by the mediator of the segment $[x_1,y_1]$
  which does not contain $a$, and this half-plane does not meet the ball
  $B(a,d(a,y_1))$.
  Similarly, we have $d(a,y_2)<d(a,y_1)$, hence a contradiction.
\end{preuve}

In particular, if $(x_1,x_2,x_3)\in V(C_0\delta(X))$ and if $x_0$ realizes the
minimum of the function $x\mapsto d(x,x_1)+d(x,x_2)+d(x,x_3)$, then
we have $(x_0,x_0,x_0,x_1,x_2,x_3)\in U$. We set
$o(x_1,x_2,x_3)=o(x_0,\ldots,x_3)$ in that case.

\begin{ps}\label{lemtechniqueordrehyperbolique}
  Let $y_1,y_2,y_3\in X$.
  The hypotheses of lemma \ref{OrdreConvexe} being convex conditions on
  $x_1$, $x_2$ and $x_3$, they define convex subsets of $X$.
  In particular, by lemma \ref{continu}, if $(y_1,y_2,y_3)\in V(C_0 \delta(X))$
  and if for every $i\in\{1,2,3\}$, $\sum_j d(x_i,y_j)<\sum_j d(y_i,y_j)$,
  then we have $o(x_1,\ldots,y_3)=o(y_1,y_2,y_3)$.
\end{ps}

\begin{ps}\label{ineq}
  Suppose that the spaces $X$ and $X'$ are the hyperbolic plane
  (equipped with a proportional distance to the usual distance)
  or an $\dR$-tree, and let $x_1,x_2,x_3\in X$, $x_1',x_2',x_3'\in X'$.
  Suppose that $(x_1,x_2,x_3)\in V(C_0(\delta(X)+a)+b)$,
  $|\delta(X)-\delta(X')|<\varepsilon_1$, and
  $|d(x_i,x_j)-d(x_i',x_j')|<2\varepsilon_2$ for all $i,j\in\{1,2,3\}$.
  Then we have
  $(x_1',x_2',x_3')\in V(C_0(\delta(X')+a-\varepsilon_1)+b-3\varepsilon_2)$.
\end{ps}

Let $(\rho,X,o), (\rho',X',o')\in m_{\Gamma}^o(2)\cup\mathcal{T}^o(\Gamma)$,
and let
$\varepsilon>0$, $K=\{x_1,\ldots,x_p\}\subset X$, and let $P$ be any finite
subset of $\Gamma$. Suppose that $(\rho',X')\in U'_{K,\varepsilon,P}(\rho,X)$.
This means (see section \ref{sectionTopGrom}) that there exists a collection
$K'=\{x_1',\ldots,x_p'\}\subset X'$ such that for every $g,h\in P$, and
every $i,j\in\{1,\ldots,p\}$ we have
\[\left|d(\rho(g)\cdot x_i,\rho(h)\cdot x_j)-d'(\rho'(g)\cdot x_i',
\rho'(h)\cdot x_j')\right|<\varepsilon
\text{ et }\left|\delta(X)-\delta(X')\right|<\varepsilon.\]

\begin{deuf}\label{defpresord}
  If $o(x_i,x_j,x_k)=o(x_i',x_j',x_k')$ for every $i,j,k\in \{1,\ldots,p\}$
  such that $(x_i,x_j,x_k)\in V(C_0(\delta(X)+\varepsilon)+3\varepsilon)$,
  we say that $K$ and $K'$ {\em come in the same order}.
\end{deuf}

Note that if
$(x_i,x_j,x_k)\in V(C_0(\delta(X)+\varepsilon)+3\varepsilon)$,
then $(x_i',x_j',x_k')\in V(C_0\delta(X'))$ and hence $o(x_i',x_j',x_k')$
is well-defined, so that definition \ref{defpresord} makes sense.

\subsubsection{Oriented equivariant Gromov topology}

Let $(\rho,X)\in m_{\Gamma}^o(2)\cup\mathcal{T}^o(\Gamma)$, let
$\varepsilon>0$, $K=\{x_1,\ldots,x_p\}\subset X$ and let $P$ be a finite subset
of $\Gamma$. We denote by $U_{K,\varepsilon,P}''(\rho,X)$ the set consisting of
the $(\rho',X')\in m_{\Gamma}^o(2)\cup\mathcal{T}^o(\Gamma)$ such that there
exists a collection $K'=\{x_1',\ldots,x_p'\}\subset X'$ such that for every
$g,h\in P$, and every $i,j\in\{1,\ldots,p\}$ we have
\[ \left|d(\rho(g)\cdot x_i,\rho(h)\cdot x_j)-d'(\rho'(g)\cdot x_i',
\rho'(h)\cdot x_j')\right|<\varepsilon
\text{ et }\left|\delta(X)-\delta(X')\right|<\varepsilon, \]
and such that $K$ and $K'$ come in the same order.

\begin{prop}\label{BaseTopo}The sets $U_{K,\varepsilon,P}''(\rho,X)$ form the
  basis of open sets of some topology; we call it the
  {\em oriented equivariant Gromov topology}.
\end{prop}

\begin{preuve}
  Of course, we always have $(\rho,X)\in U_{K,\varepsilon,P}''(\rho,X)$.
  Hence, we need only check that if
  $(\rho,X)\in U_{K_1,\varepsilon_1,P_1}''(\rho_1,X_1)\cap
  U_{K_2,\varepsilon_2,P_2}''(\rho_2,X_2)$
  then, for some $K$, $\mu>0$ and $P$ we have
  $U_{K,\mu,P}''(\rho,X)\subset U_{K_1,\varepsilon_1,P_1}''
  (\rho_1,X_1)\cap U_{K_2,\varepsilon_2,P_2}''(\rho_2,X_2)$.
  It is just a technical verification.
  Denote $K_1=\{a_1',\ldots,a_{n_1}'\}$ and
  $K_2=\{b_1',\ldots,b_{n_2}'\}$. Then there exists a finite collection
  $K=\{a_1,\ldots,a_{n_1},b_1,\ldots,b_{n_2}\}\subset X$ such that:
  \[\forall i,j\leq n_1,\forall\gamma_1,\gamma_2\in P_1,
  |d_X(\rho(\gamma_1)\cdot a_i,\rho(\gamma_2)\cdot a_j)-
  d_{X_1}(\rho_1(\gamma_1)\cdot a_i',\rho_1(\gamma_2)\cdot a_j')|
  <\varepsilon_1,\]
  $|\delta(X)-\delta(X_1)|<\varepsilon_1$,
  and such that $o_{X_1}(a_i',a_j',a_k')=o_X(a_i,a_j,a_k)$ for every $i,j,k$
  such that
  $(a_i',a_j',a_k')\in V(C_0(\delta(X_1)+\varepsilon_1)+3\varepsilon_1)$;
  similarly $o_{X_1}(b_i',b_j',b_k')=o_X(b_i,b_j,b_k)$ for every $i,j,k$ such
  that $(b_i',b_j',b_k')\in V(C_0(\delta(X_1)+\varepsilon_1)+3\varepsilon_1)$.
  Since all these {\em finitely many} inequalities are strict, they can all be
  refined simultaneously by some $\mu>0$.
  In other words,
  \begin{equation}\label{eqn1}
    |d_X(\rho(\gamma_1)\cdot a_i,\rho(\gamma_2)\cdot a_j)-
    d_{X_1}(\rho_1(\gamma_1)\cdot a_i',\rho_1(\gamma_2)\cdot a_j')|
    <\varepsilon_1-\mu,
  \end{equation}
  \begin{equation}\label{eqn2}
    |\delta(X)-\delta(X_1)|<\varepsilon_1-\mu,
  \end{equation}
  and the condition
  $o_{X_1}(a_i',a_j',a_k')=o_X(a_i,a_j,a_k)$ is verified for all
  \begin{equation}\label{eqn3}
    (a_i',a_j',a_k')\in
    V(C_0(\delta(X_1)+\varepsilon_1+\mu)+3\varepsilon_1+3\mu).
  \end{equation}
  Now, take $(\rho'',X'')\in U_{K,\mu,P_1\cup P_2}''(\rho,X)$.
  Then in particular, for every $i,j\leq n_1$ and every
  $\gamma_1,\gamma_2\in P_1$,
  \begin{equation}\label{eqn4}
    |d_{X''}(\rho''(\gamma_1)\cdot a_i'',\rho''(\gamma_2)\cdot a_j'')-
    d_X(\rho(\gamma_1)\cdot a_i,\rho(\gamma_2)\cdot a_j)|<\mu,
  \end{equation}
  $|\delta(X)-\delta(X'')|<\mu$, and for every $i,j,k$ such that
  $(a_i,a_j,a_k)\in V(C_0(\delta(X)+\mu)+3\mu)$,
  $o_X(a_i,a_j,a_k)=o_{X''}(a_i'',a_j'',a_k'')$.
  Now, the conditions (\ref{eqn1}) and (\ref{eqn4}) imply that for every
  $i,j\leq n_1$ and every $\gamma_1,\gamma_2\in P_1$,
  \[ |d_{X''}(\rho''(\gamma_1)\cdot a_i'',\rho''(\gamma_2)\cdot a_j'')-
  d_{X_1}(\rho_1(\gamma_1)\cdot a_i',\rho_1(\gamma_2)\cdot a_j')|
  <\varepsilon_1, \]
  and, by remark \ref{ineq}, the conditions (\ref{eqn1}), (\ref{eqn2})
  and (\ref{eqn3}) imply that
  \[(a_i,a_j,a_k)\in V(C_0(\delta(X)+2\mu)+\dfrac{3}{2}\varepsilon_1+
  \dfrac{9}{2}\mu)\subset V(C_0(\delta(X)+\mu)+3\mu),\]
  so that
  $o_{X''}(a_i'',a_j'',a_k'')=o_{X}(a_i,a_j,a_k)=o_{X_1}(a_i',a_j',a_k')$.
  In that way we get
  $(\rho'',X'')\in U_{K_1,\varepsilon_1,P_1}''(\rho_1,X_1)$,
  and, similarly,
  $(\rho'',X'')\in U_{K_2,\varepsilon_2,P_2}''(\rho_2,X_2)$.
\end{preuve}

\begin{prop}\label{memtop}
  The oriented equivariant Gromov topology coincides with the usual topology
  on $m_{\Gamma}^o(2)$.
\end{prop}

\begin{preuve}
  It is the same proof as the one of F. Paulin's proposition 6.2 in
  \cite{Paulin88}, with minor modifications. In that proof
  (recalled here, as proposition \ref{propPaulin}), the only difference is that
  the isometry $\phi$ of \ref{lemmememetop}, is now an
  orientation-preserving isometry.
\end{preuve}

Of course, we denote by $\overline{m_{\Gamma}^o(2)}$ the closure of
$m_{\Gamma}^o(2)$ in the space $m_{\Gamma}^o(2)\cup\mathcal{T}^o(\Gamma)$,
equipped with the oriented equivariant Gromov topology.

We also write $\overline{m_\Gamma^u}=\overline{m_{\Gamma}^u(2)}$ and
$\overline{m_\Gamma^o}=\overline{m_{\Gamma}^o(2)}$.

\subsubsection{The space $\overline{m_\Gamma^o}$ is compact}

Denote by $\pi:\overline{m_\Gamma^o}\rightarrow\overline{m_\Gamma^u}$ the
natural function consisting in forgetting the orientation.

\begin{prop}\label{fibcomp}
  The map $\pi$ is continuous, and its fibres are compact Hausdorff.
\end{prop}

\begin{preuve}
  First, the continuity of $\pi$ follows directly from the definition of these
  two topologies.

  Now, let $(\rho,T)\in\overline{m_\Gamma^u}$. If $T$ is $\dH^2$, then the
  fibre $\pi^{-1}(\rho,T)$ has cardinal $2$ in the Hausdorff space
  $m_\Gamma^o$ (by theorem \ref{mg}), hence it is compact Hausdorff.
  Suppose now that $T$ is an $\dR$-tree.
  By definition, the {\em set} $\pi^{-1}(\rho,T)$ is a subset of $O_T$. By
  definition of the sets $U_{K,\varepsilon,P}''(\rho,T)$, the induced topology
  on $\pi^{-1}(\rho,T)$ in $O_T$ coincides with the oriented equivariant Gromov
  topology (in each of these topologies, open sets are defined by equalities of
  $o$ on finite subsets of $Trip(T)$).

  We shall notice now that $\overline{m_\Gamma^o}$ is Hausdorff. Indeed,
  if $(\rho,X,o)$ and $(\rho',X',o')$ are distinct and are not separated by
  open sets,
  then $\delta(X)=\delta(X')$. The open set $m_\Gamma^o$ being Hausdorff, this
  means that $X$ and $X'$ are fat $\dR$-trees. Since the map
  $\pi\colon\overline{m_\Gamma^o}\rightarrow\overline{m_\Gamma^u}$ is
  continuous, this implies that these two spaces differ only by the
  orientation: but by definition of the oriented equivariant Gromov topology,
  there are two open sets separating $(\rho,X,o)$ and $(\rho',X',o')$.

  Consequently, we need only prove that $\pi^{-1}(\rho,T)$ is closed in $O_T$.
  Let $E_T''$ be the set of orders
  $o\colon(\partial T)^3\rightarrow\{-1,0,1\}$ which satisfy the hypotheses of
  definition \ref{defordrecyclique}, and which are invariant by $\rho$. These
  are closed conditions, hence $E_T''$ is closed in $f(O_T)$. By definition,
  $E_T''\subset\mathcal{T}^o(\Gamma)$, and
  $\overline{m_\Gamma^o}\cap\mathcal{T}^o(\Gamma)$ is closed in
  $\mathcal{T}^o(\Gamma)$, hence $\pi^{-1}(\rho,T)$ is closed in $E_T''$.
  It follows that $\pi^{-1}(\rho,T)$ is closed in a compact Hausdorff space,
  hence it is compact.
\end{preuve}

\begin{theo}
  The space $\overline{m_\Gamma^o}$, equipped with the function
  $m_\Gamma^o\hookrightarrow\overline{m_\Gamma^o}$, is a natural
  compactification of $m_\Gamma^o$.
\end{theo}

Here again, after \cite{Paulin04}, by ``natural'', we mean that the action of
$Out(\Gamma)$ on $m_\Gamma^o$ extends continuously to an action of
$Out(\Gamma)$ on $\overline{m_\Gamma^o}$.

\begin{preuve}
  Since $m_\Gamma^o$ is open and dense in $\overline{m_\Gamma^o}$ and since, by
  definition of the oriented equivariant Gromov topology, the action of
  $Out(\Gamma)$ on $\overline{m_\Gamma^o}$ is continuous, it suffices to prove
  that the space $\overline{m_\Gamma^o}$ is compact Hausdorff. We have already
  seen that $\overline{m_\Gamma^o}$ is Hausdorff, and by the definition of
  compacity in terms of ultrafilters (see \textsl{e.g.} \cite{Bourbaki},
  page 59), we need only prove that every ultrafilter in $\overline{m_g^o}$
  converges.

  Let $\omega$ be an ultrafilter in $\overline{m_\Gamma^o}$.
  Then the image of the ultrafilter $\pi(\omega)$ is an ultrafilter in the
  compact space $\overline{m_\Gamma^u}$ (see \textsl{e.g.} \cite{Bourbaki},
  proposition 10, page 41) hence it converges to some action
  $(\rho_\infty,X_\infty)\in\overline{m_\Gamma^u}$. If $X_\infty=\dH^2$, then
  it follows from proposition \ref{memtop} that $X$ is equipped with an
  orientation, compatible with the convergence of the ultrafilter.
  We need only prove that if $X_\infty$ is an $\dR$-tree
  (denote $(\rho,T)=(\rho_\infty,X_\infty)$ in that case)
  then there exists a coherent order
  $o\colon (\partial_\infty T)^3\rightarrow\{-1,0,1\}$
  which verifies the hypotheses of definition \ref{defordrecyclique},
  which is invariant under the action of $\Gamma$, and such that $(\rho,T)$,
  equipped with this order, is indeed the limit, in $\overline{m_\Gamma^o}$,
  of the ultrafilter $\omega$.

  Consider an increasing sequence $T_k\subset T$ of finite, closed
  subtrees of $T$, such that $\displaystyle{\bigcup_k}\,T_k=T$.
  Suppose for simplicity that $T_1$ is a singleton $\{x_0\}$.
  For all $k\geq 1$, $N\geq 0$, denote by $F_{k,N}$ the finite collection of
  elements of $T_k$
  consisting of all the end points of $T_k$, as well as all the points of $T_k$
  whose distance to $x_0$ is a multiple of $\dfrac{1}{2^N}$.
  Since the ultrafilter $\pi(\omega)$ converges to $(\rho,T)$,
  for all $k,N$ and $\varepsilon>0$, and for every finite subset $P$ of
  $\Gamma$, for all $M\in\omega$, there exists $(\rho_M,X_M,o_M)\in M$ such
  that $(\rho_M,X_M)\in U_{F_{k,N},\varepsilon,P}'(\rho,T)$, \textsl{i.e.}
  there exists a $P$-equivariant $\varepsilon$-approximation between
  $F_{k,N}$ and some finite collection $K_M$ in $X_M$, for some
  $(\rho_M,X_M,o_M)\in M$.
  For all $M\in\omega$,
  denote by $O_{k,N,\varepsilon,P,M}$ the set of functions $o\in O_T$ such that
  there exists such an approximation, and such that for every
  $x_1,x_2,x_3\in F_{k,N}$ with
  $(x_1,x_2,x_3)\in V\left((C_0+3\varepsilon)+\dfrac{1}{2^N}\right)$,
  and for all corresponding $x_1',x_2',x_3'\in K_M$, we have
  $o(Trip(x_1,x_2,x_3))=o_M(x_1',x_2',x_3')$
  (note that this is indeed well defined, thanks to remark
  \ref{lemtechniqueordrehyperbolique}).
  
  We cut the end of the proof into the two following lemmas:
  \begin{lem}\label{un}
    For every $k,N,\varepsilon$, $P$ finite subset of $\Gamma$ and every
    $M\in\omega$, the set
    $O_{k,N,\varepsilon,P,M}$ is
    closed, and non empty. Moreover, if $k>k'$, $N>N'$,
    $\varepsilon<\varepsilon'$, $P'\subset P$ and $M\subset M'$, then
    $O_{k,N,\varepsilon,P,M}\subset O_{k',N',\varepsilon',P',M'}$.
  \end{lem}
  By compacity of $E_T'$, it follows that
  $\displaystyle{\bigcap_{k,N,\varepsilon,P,M}}
  O_{k,N,\varepsilon,P,M}\neq\emptyset$.
  \begin{lem}\label{deux}
    Let
    $o\in\displaystyle{\bigcap_{k,N,\varepsilon,P,M}}
    O_{k,N,\varepsilon,P,M}$.
    Then $o$ satisfies the conditions of definition \ref{defordrecyclique}.
    Hence, $o$ defines a coherent total cyclic order. Moreover, it is invariant
    under the action of $\Gamma$, and the element $(\rho,T)$, equipped with the
    orientation $o$, is the limit, in $\overline{m_\Gamma^o}$, of
    the ultrafilter $\omega$.
  \end{lem}
\end{preuve}

\begin{preuvede} of lemma \ref{un}.

  It follows from the definition that
  $O_{k,N,\varepsilon,P,M}\subset O_{k',N',\varepsilon',P',M'}$
  if $k>k'$, $\varepsilon<\varepsilon'$, $M\subset M'$ and $P'\subset P$.

  The hypotheses concern only $Trip(T_k)$, which is a {\em finite} subset of
  $Trip(T)$, and hence $O_{k,N,\varepsilon,P,M}$ is closed.
  We need to prove that it is also non empty. Up to restricting ourselves to a
  subset of it, we may suppose that $\varepsilon$ is small enough and that $N$
  is large enough so that for every triple $u_1,u_2,u_3$ of non aligned branch
  points of $T_k$, we have
  $(u_1,u_2,u_3)\in V\left((C_0+3)\varepsilon+\dfrac{1}{2^N}\right)$.
  Choose a (non-oriented) $\varepsilon$-approximation between $E_{k,N}$ and
  $K_M\subset X_M$ (such an approximation exists, since
  $\pi(\omega)$ converges to $(\rho,T)$ in $\overline{m_g^u}$). For each tripod
  $Trip(a_0,a_1,a_2,a_3)\in Trip(T_k)$ (where $a_0$ is the center of the
  tripod), the only obstruction to defining
  $o(Trip(a_0,a_1,a_2,a_3))\in\{-1,1\}$ according to $K_M$, would be that there
  exist two representants $(a_0,a_1,a_2,a_3)$ and $(a_0,b_1,b_2,b_3)$ of this
  same class of tripods, such that the triples $(a_1,a_2,a_3)$ and
  $(b_1',b_2',b_3')$ are in
  $V\left((C_0+3)\varepsilon+\dfrac{1}{2^N}\right)$, and such that
  $o_{X_M}(a_1',a_2',a_3')=1=-o_{X_M}(b_1',b_2',b_3')$. Hence we need to prove
  that this is impossible.
  The condition
  $(u_1,u_2,u_3)\in V\left((C_0+3)\varepsilon+\dfrac{1}{2^N}\right)$ for all
  the non aligned branched points of $T_k$, implies that it is possible to go
  from the triple $(a_1,a_2,a_3)$ to the triple $(b_1,b_2,b_3)$ by a sequence
  of moves consisting of replacing $a_1$, $a_2$ or $a_3$ by one of its close
  neighbours in $T_k$, without leaving
  $V\left((C_0+3)\varepsilon+\dfrac{1}{2^N}\right)$.
  In particular, we may suppose that $(b_1,b_2,b_3)=(b_1,a_2,a_3)$ and
  $d_T(a_1,b_1)\leq\dfrac{1}{2^N}$. This implies that
  $d_{X_M}(a_1',b_1')<\varepsilon+\dfrac{1}{2^N}$. Take
  $a'(t)\in[a_1',b_1']$, with $a'(0)=a_1'$ and $a'(1)=b_1'$. Then, for all
  $t\in[0,1]$, and every permutation $(i,j,k)$ of $(1,2,3)$, we have
  $d(a_i',a_j')+d(a_j',a_k')-d(a_i',a_k')>
  2\left(\left((C_0+3)\varepsilon+\dfrac{1}{2^N}\right)-3\varepsilon
  -\dfrac{1}{2^N}\right)$,
  so that for all $t\in[0,1]$, we have
  $(a'(t),a_2',a_3')\in V(C_0\delta(X))$.
  It follows from the continuity of the order
  (and more precisely, from lemmas \ref{Toricelli} and \ref{continu})
  that $o_{X_M}(a_1',a_2',a_3')=o_{X_M}(b_1',a_2',a_3')$.
\end{preuvede}

\begin{preuvede} of lemma \ref{deux}.

  It follows from the definitions of
  the sets $O_{k,N,\varepsilon,P,M}$ that the order
  $o:(\partial_\infty T)^3\rightarrow\{-1,1\}$ indeed defines an orientation,
  \textsl{i.e.}, it is coherent, and satisfies the first two conditions of
  definition \ref{defordrecyclique}.
  The third condition, as well as the invariance of $o$ under the action of
  $\Gamma$, are proved by considering a big enough subtree $T_k$ of $T$
  containing the desired branched points, and by deriving the properties of
  $o$ from the corresponding properties for $o_M$, which are supposed to be
  true since $(\rho_M,X_M,o_M)\in\overline{m_\Gamma^o}$. As an example we prove
  that $o$ satisfies the third condition of definition \ref{defordrecyclique};
  the proof of the invariance of $o$ under $\Gamma$ is similar.
  Let $x_1,x_2,x_3,x_4\in\partial_\infty T$ be such that
  $o(x_1,x_2,x_3)=o(x_1,x_3,x_4)=1$; we want to prove that $o(x_1,x_2,x_4)=1$.
  Let $k$ and $N$ be sufficiently large so that $T_k$ contains five points
  $a_0,a_1,a_2,a_3,a_4\in F_{k,N}$ such that every triple $i,j,k$ of distinct
  elements of $\{1,2,3,4\}$, we have
  $(a_0,a_0,a_0,a_i,a_j,a_k)\in U$ in $T_k$, and such that for all
  $i\in\{1,2,3,4\}$, the ray issued from $a_0$ and passing through $x_i$ also
  passes through $a_i$, and such that $a_0\in Conv(a_1,a_2,a_3,a_4)$
  (such points indeed exist in $T$).
  Then, for all $\varepsilon$-approximation between $F_{k,N}$ and
  $K_M\subset X_M$, with $\varepsilon$ small enough, we have
  $(a_0',a_0',a_0',a_i',a_j',a_k')\in U$ in $X_M$, and then the equality
  $o(x_1,x_2,x_4)=1$ indeed follows from the fact that $o_n$
  satisfies the third condition of definition \ref{defordrecyclique}.

  Now, for all $M\in\omega$, we have
  $U_{F_{k,N},\varepsilon,P}''(\rho,T,o)\cap M\neq\emptyset$
  (for all $\varepsilon'$ such that
  $(C_0+3)\varepsilon<(C_0+3)\varepsilon'+\dfrac{1}{N}$). This means, by
  definition, that the point $(\rho,T,o)\in\overline{m_g^o}$ is adherent
  to the filter $\omega$, and since $\omega$ is an ultrafilter this implies
  that $\omega$ converges to $(\rho,T,o)$.
\end{preuvede}

\begin{cor}\label{onto}
  The map $\pi\colon\overline{m_\Gamma^o}\rightarrow\overline{m_\Gamma^u}$ is
  onto.
\end{cor}

\begin{preuve}
  Of course, $m_\Gamma^o\rightarrow m_\Gamma^u$ is onto. Now, let
  $T\in\partial\overline{m_\Gamma^u}$, and $\rho_n\in m_\Gamma^o$ be such that
  $\pi(\rho_n)$ converges to $T$. Then $\rho$ possesses a subsequence
  converging to some $\rho_\infty$, and by continuity of $\pi$ we have
  $\pi(\rho_\infty)=T$.
\end{preuve}



\begin{ps}
  It is possible to write this proof without using ultrafilters
  (see \cite{TheseMaxime}), and to prove the sequential compactness first
  (by considering a sequence instead of an ultrafilter on $\overline{m_g^o}$),
  and then to prove the compactness of $\overline{m_g^o}$ by using
  elementary general topology.
\end{ps}

\subsubsection{The space $\overline{m_g^o}$ possesses $4g-3$ connected components}

We are now going to focus on the case when $\Gamma=\pi_1\Sigma_g$.
We denote
$\mathcal{T}_g=\mathcal{T}(\pi_1\Sigma_g)$, and
$\mathcal{T}_g^o=\mathcal{T}(\pi_1\Sigma_g)^o$.
If $(\rho,T)\in \mathcal{T}_g^o$, the set $\partial_\infty T$ is equipped with
a total cyclic order, preserved by the action of $\pi_1\Sigma$, and hence it
possesses an Euler class, as defined in section
\ref{SectionEulerClass}.
Notice that if $T$ is a line, then it follows from the definition of the Euler
class that $e(\rho,T)=0$.

\begin{theo}\label{EulerContinu}
  The Euler class $e:m_g^o\cup\mathcal{T}_g^o\rightarrow\dZ$
  is a continuous function.
\end{theo}

This proof will use the technical statements established in section
\ref{SectionEulerClass}, which imply that we need only finitely many
information about the order in order to compute the Euler class of a
representation. We will be using here the notations introduced in that
section.

\begin{preuve}
  First, the set $m_g^o=\{(\rho,X)|\delta(X)\neq 0\}$ is open
  in $m_g^o\cup\mathcal{T}_g^o$, and it follows from the
  formula (\ref{AlgMil}), in section \ref{SectionAlgMil},
  that $e$ is continuous on $m_g^o$.

  Now take an element $(\rho_T,T)\in\overline{m_g^o}\smallsetminus m_g^o$.
  We shall prove that there exists a neighbourhood of $T$, in the sense of the
  oriented equivariant Gromov topology, consisting only in representations of
  the same Euler class as $T$. First suppose that $T$ is not reduced to a line,
  so that it has at least three ends (in that case, $T$ possesses infinitely
  many ends). Take $x,y\in\partial_\infty T$, such that
  $x\not\in P_{ref}\cdot y$. We may suppose that $Card(P_{ref}\cdot y)\geq 2$;
  otherwise $\pi_1\Sigma_g$ would fix every end of $T$, and, by minimality,
  $T$ would be a point or a line.
  
  Every triple $\{a,b,c\}$ of pairwise distinct elements of
  $P_{ref}\cdot\{x,y\}$ determines a unique class of tripods in $T$; we will
  denote by $P_{abc}$ the centre of this tripod.
  Let $K_1$ be the convex hull
  \[K_1=Hull\left\{P_{abc}\,|\, (a,b,c)\in\left(P_{ref}\cdot\{x,y\}\right)^3,
  Card\{a,b,c\}=3\right\}.\]
  Put $d_T=\max\{d(p,\gamma p)\,|\, p\in K_1, \gamma\in P_{ref}\}$ and let
  $d_{K_1}$ be the diameter of $K_1$. Also, for every non degenerate triple
  $\{a,b,c\}\subset P_{ref}\cdot\{x,y\}$, let
  $P_{abc}^a\subset \,(a,b)\,\cap\, (a,c)$ be such that
  $d(P_{abc}^a,P_{abc})=L$, with $L>9d_T+3d_{K_1}$, and let $K_2$ be the
  convex hull
  \[ K_2=Hull\left\{P_{abc}^a| (a,b,c)\in\left(P_{ref}\cdot\{x,y\}\right)^3,
  Card\{a,b,c\}=3\right\}. \]
  Finally, fix a point $p_0\in K_1$.

  
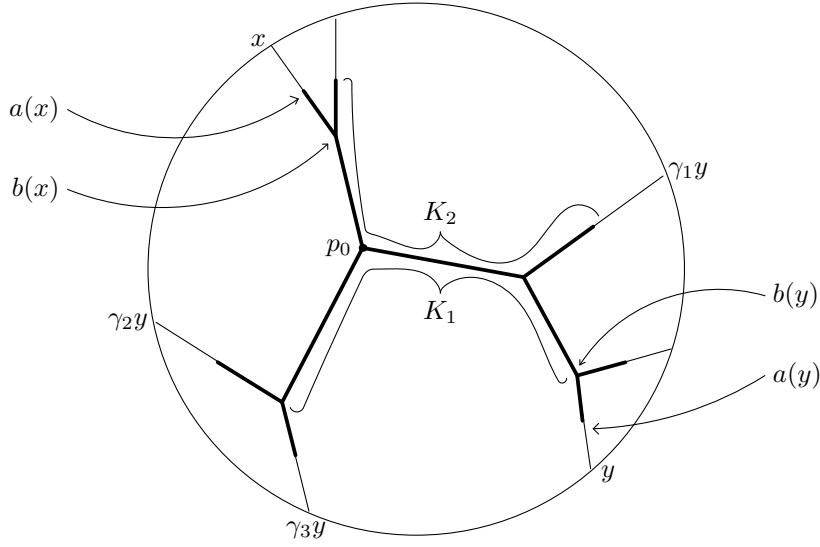
\[ \epsfbox{Figures/FiguresCompactification.2} \]

  Let $a(x)=\displaystyle{\bigcap_{p\in K_2}}[p,x)\,$.
  This is the ``closest point to $x$'' in $K_2$. And let $b(x)$ be the
  projection of $a(x)$ on $K_1$. We define $a(y)$ and $b(y)$ similarly.
  We put $K=\left\{p_0,a(x),b(x),a(y),b(y)\right\}$ and we consider
  $(\rho,X)\in U_{K,\frac{d_T}{C_0+3},P_{ref}}''(\rho_T,T)$.
  We shall prove that $(\rho_T,T)$ and $(\rho,X)$ have the same Euler class,
  by applying proposition \ref{nbrefini2}.
  In the space $X$, denote by
  $p_0',a'(y),b'(y),a'(x),b'(x)$ the corresponding points.
  Denote by $x'\in\partial_\infty X$ the end of some ray
  $[b'(x),a'(x))$
  and by $y'\in\partial_\infty X$ the end of some ray $[b'(y),a'(y))$
  (chosen arbitrarily, in the case of an $\dR$-tree).

  Let $\gamma_1,\gamma_2,\gamma_3\in P_{ref}$ be such that
  $o(\gamma_1 x,\gamma_2 x,\gamma_3 y)=1$. We want to prove that
  $o(\gamma_1 x',\gamma_2 x',\gamma_3 y')=1$.
  For this, we shall prove the three following equalities:
  \begin{equation}\label{eqn5}
    o(\gamma_1 x,\gamma_2 x,\gamma_3 y)=
    o(\gamma_1 a(x),\gamma_2 a(x),\gamma_3 a(y)),
  \end{equation}
  \begin{equation}\label{eqn6}
    o(\gamma_1 x',\gamma_2 x',\gamma_3 y')=
    o(\gamma_1 a'(x),\gamma_2 a'(x),\gamma_3 a'(y)),
  \end{equation}
  \begin{equation}\label{eqn7}
    o(\gamma_1 a(x),\gamma_2 a(x),\gamma_3 a(y))=
    o(\gamma_1 a'(x),\gamma_2 a'(x),\gamma_3 a'(y)).
  \end{equation}
  We can check that for every $\gamma\in P_{ref}$,
  $L-d_T\leq d(\gamma a(y),K_1)\leq L+d_T$, and the centre of the tripod
  determined by $\gamma_1 a(x)$, $\gamma_2 a(x)$ and $\gamma_3 a(y)$
  lies in $K_1$, hence
  \[(\gamma_1 a(x),\gamma_2 a(x),\gamma_3 a(y))\in V(8d_T+3d_{K_1})
  \subset V(d_T),\]
  so that all the terms of the equations (\ref{eqn5}), (\ref{eqn6}) and
  (\ref{eqn7}) are well-defined, and equation (\ref{eqn7}) holds.
  Put $p_1=\gamma_1a(x)$, $p_1'=\gamma_1a'(x)$, \ldots, $p_3=\gamma_3a(y)$,
  $p_3'=\gamma_3a'(y)$.
  Then $d_{X}(p_0',p_i')<d_{K_1}+L+d_T+\frac{d_T}{C_0+3}$, and
  $d_{X}(p_i',p_j')> 2L-2d_T-\frac{d_T}{C_0+3}$, so that for every
  $i\in\{1,2,3\}$, $\sum_j d_X(p_0',p_j')<\sum_j d_X(p_i',p_j')$.
  These inequalities are even finer in $T$, and by lemma \ref{OrdreConvexe} and
  remark \ref{lemtechniqueordrehyperbolique}, the equalities (\ref{eqn5})
  and (\ref{eqn6}) hold.
  Similarly, if $\gamma_1,\gamma_2,\gamma_3\in P_{ref}$ are such that
  $o(\gamma_1 x,\gamma_2 y,\gamma_3 y)\neq 0$ then
  $o(\gamma_1 x',\gamma_2 y',\gamma_3 y')=o(\gamma_1 x,\gamma_2 y,\gamma_3 y)$,
  so that the conditions of proposition \ref{nbrefini2} are satisfied. This
  finishes the proof, in the case when $T$ is not a line.

  \smallskip

  Now, suppose that $T$ is a line. We want to prove that there exists a
  neighbourhood of $T$ consisting only in representations of Euler class zero.
  For simplicity we will prove the following:
  \begin{lem}\label{continu2}
    There exists a neighbourhood $V'$ of $T$ in which every fat $\dR$-tree has
    Euler class zero.
  \end{lem}
  This lemma implies the theorem, for the following reason.
  For all $k\neq 0$ such that $|k|\leq 2g-2$, denote by
  $L\subset\mathcal{T}_g^o$ the set of actions on lines, and denote by
  $F_k\subset L$ the set of actions on lines $(\rho,T)$ such
  that every neighbourhood of $(\rho,T)$ contains representations of Euler
  class $k$.
  By lemma \ref{continu2}, we have
  $F_k=L\cap\overline{m_{g,k}^o}$. Hence, $F_k$ is a closed
  subset of $\partial m_{g,k}^o$ (indeed, $L$ is a closed set, as, by
  definition of the topology, it is an open condition to contain a non
  degenerate tripod).
  Now, let $(\rho,T)\in F_k$. By lemma \ref{continu2}, $(\rho,T)$ has a
  neighbourhood $V\subset\mathcal{T}_g^o\cup m_g^o$ in which every fat
  $\dR$-tree has Euler class $0$. Put $V'=V\cap\partial {m_{g,k}^o}$. It
  is an open subset of $\partial m_{g,k}^o$. If there was a tree
  $(\rho',T')\in V'$ not reduced to a line, then by the preceding argument,
  $(\rho',T')$ has a neighbourhood consisting of actions (on hyperbolic
  planes or on trees) of Euler class $0$, which is a contradiction since
  $(\rho',T')\in\overline{m_{g,k}^o}$. Hence, $V'$ consists of actions on
  lines, \textsl{i.e.}, $F_k$ is open in $\partial m_{g,k}^o$. By
  proposition \ref{oneend}, the space $m_{g,k}^o$ is one-ended, hence
  $\partial m_{g,k}^o$ is connected. And we can prove easily that
  $\partial m_{g,k}^o$ contains actions not reduced to a line, hence
  $F_k\neq \partial m_{g,k}^o$. Whence, $F_k=\emptyset$, for all $k\neq 0$.
  In other words, every action on a line has a neighbourhood consisting of
  actions of Euler class $0$; this finishes the proof of theorem
  \ref{EulerContinu}.
\end{preuve}

\begin{preuvede} of lemma \ref{continu2}.

  Up to add some elements, we shall suppose here that the set
  $P_{ref}$ contains $S$ and that it is symmetric.
  Let $(\rho,T)$ be a line, such that
  $\displaystyle{\min_{x_0\in T}
  \max_{\gamma\in S}}\,d(x_0,\gamma\cdot x_0)=1$.
  Consider some point $x_0\in T$ as in the equality above.
  Let $d_1$ be the greatest distance between $x_0$ and $\gamma x_0$, for every
  $\gamma\in P_{ref}$.
  Consider points $x_1$, $x_2$, $y_1$, $y_2$ of $T$ such that $x_i$, $y_i$ are
  on the same side of $x_0$, such that $d(x_0,y_i)=d_1+4$ and
  $d(x_0,x_i)=2 d_1+6$.
  Let $K$ be the finite subset of $T$ consisting of $x_1$, $x_2$, $y_1$, $y_2$
  and $P_{ref}\cdot x_0$. Let
  $(\rho',X')\in U_{K,\frac{1}{6},P_{ref}}''(\rho,T)$, we shall prove that
  $e(\rho',T')=0$.
  If $T'$ is a line, then there is nothing to do.
  Otherwise, for every point $p_i'\in K'$ approximating $K$, denote by
  $p_i''$ its projection on the segment $[x_1',x_2']$. This defines a new
  approximation, which realizes $(\rho',T')$ as an element of
  $U_{K,1,P_{ref}^{n'}}''(\rho,T)$, and such that $K''$ is contained in a
  segment.
  Let $r$ be the end of a ray issued from $x_0''$ and which leaves the segment
  $[x_1',x_2']$ at a point $p_0$ at distance at most $2$ of $x_0''$
  (such a ray does exist, since
  $\displaystyle{\max_{\gamma\in S}}\,d(x_0'',\gamma\cdot x_0'')<2$).
  For every $\gamma\in P_{ref}$, $d(\gamma\cdot p_0,x_0'')<d_1+3$, and hence
  the segment $[x_1',x_2']\cap Hull(P_{ref}\cdot p_0)$ is contained
  (strictly, on each side),
  in the segment $[y_1'',y_2'']$. Similarly, for every $\gamma\in P_{ref}$,
  $[y_1'',y_2'']\subset [\gamma\cdot x_1',\gamma\cdot x_2']$.
  For $i=1,2$, let $U_i$ be the set of ends of rays issued from $x_0''$ and
  passing through $x_i'$, and let $U_2'$ be the set of ends of rays issued from
  $x_0''$ and passing through $y_i''$. Then $P_{ref}$ sends $U_1\cup U_2$ on a
  subset of $U_1'\cup U_2'$.
  Then it follows from the coherence condition on the order on $T'$ that for
  all $x\in U_1'$, $y\in U_2'$, $o(r,x,y)\in\{-1,1\}$ is constant.
  Suppose for instance that $o(r,x,y)=1$ for all $x\in U_1'$, $y\in U_2'$.
  Thus, for every $\gamma\in P_{ref}$, the situation is the following.
  \begin{itemize}
    \item If $\gamma$ sends $U_1$ on a subset of $U_1'$, and $U_2$ on a subset
      of $U_2'$ (or, equivalently, if $\rho(\gamma)$ preserves the orientation
      of $T$), then $o(\gamma\cdot r,x,y)=1$ for all $x\in U_1'$, $y\in U_2'$.
      Denote by $A$ the set of these ends $\gamma\cdot r$.
    \item If $\gamma$ sends $U_1$ on a subset of $U_2'$, and $U_2$ on a subset
      of $U_1'$ (or, equivalently, if $\rho(\gamma)$ reverses the orientation
      of $T$), then $o(\gamma\cdot r,x,y)=-1$ for all $x\in U_1'$, $y\in U_2'$.
      Denote by $B$ the set of such ends $\gamma\cdot r$.
  \end{itemize}
  Now equip the set $\{a,u_1,b,u_2\}$ with the cyclic order in which we wrote
  them here. Denote by $h$ the order-preserving bijection which exchanges $a$
  and $b$ and exchanges $u_1$ and $u_2$. Then we can consider the action
  $\pi_1\Sigma_g\rightarrow\{1,h\}$ on this ordered set, defined as follows: if
  $\gamma\in \pi_1\Sigma_g$ preserves the orientation of the line $T$ then we
  send it to $1$, otherwise we send it to $h$. Of course, this action has Euler
  class zero, and now it follows from proposition \ref{nbrefini2}, which
  applies here, that $e(\rho',T')=0$.
\end{preuvede}

\subsection{Degeneracy of the non oriented compactification}

\subsubsection{Non-orientable $\dR$-trees}

Now we are going to prove that the existence of an orientation, on an
$\dR$-tree, preserved by the action of the group, is indeed a restrictive
condition. More precisely:

\begin{prop}
  Let $g\geq 3$. Then the inclusion
  $\partial\overline{m_g^u(2)}\subset\partial\overline{m_g^u(3)}$ is strict.
\end{prop}

\begin{preuve}
  Of course, every isometric embedding of $\dH^2$ into $\dH^3$ gives rise to an
  embedding $m_g^u(2)\subset m_g^u(3)$ and it follows that
  $\partial\overline{m_g^u(2)}\subset\partial\overline{m_g^u(3)}$.
  In order to prove that this inclusion is strict, we shall prove that there
  exists an element $(T,\rho_\infty)\in\partial\overline{m_g^u(3)}$ such that
  no orientation on $T$ is preserved by $\rho_\infty$. Since the map
  $\pi\colon\overline{m_g^o(2)}\rightarrow\overline{m_g^u(2)}$ is onto, this
  implies that $(T,\rho_\infty)\not\in\partial\overline{m_g^u(2)}$.

  Consider the realization of $\pi_1\Sigma_2$ as a Fuchsian group acting on
  $\dH^2$ with a fundamental domain as symmetric as possible, that is, a
  regular octagon, $\rho_0\colon\pi_1\Sigma_2\rightarrow PSL(2,\dR)$.
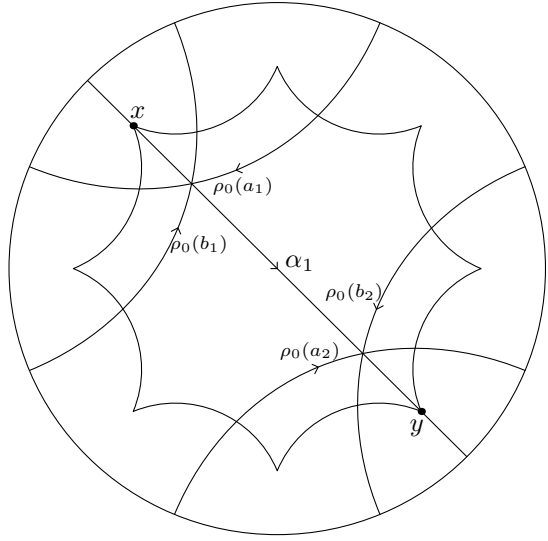
  $$\epsfbox{Figures/FiguresExtraComp.4}$$
  Denote by $\rho_0(a_1),\ldots,\rho_0(b_2)\in PSL(2,\dR)$ the hyperbolic
  isometries suggested in the above picture. The element
  $\alpha_1=\rho_0(a_1^{-1}b_1^{-1}a_2b_2)$ is hyperbolic, of axis
  $(x,y)$ represented above (indeed, $\alpha_1\cdot x=y$, and if
  $\overrightarrow{u}$ is a unit tangent vector at $x$ pointing towards $y$,
  the angles of $\overrightarrow{u}$ and of its successive images with the
  edges of the octagon enable to check that the image of $\overrightarrow{u}$
  is again a vector whose direction is the one of the axis $(x,y)$, pointing
  in the opposite direction of $x$).
  Note that $\alpha_1$ is represented by a non separating simple closed curve
  on the surface $\Sigma_2$
  (in order to check that this curve is non separating, we can see for instance
  that the isometry $\rho_0(a_1)$ identifies two points of the boundary of the
  fundamental domain represented by the above picture; these two points are in
  the two different sides of the axis of $\alpha_1$; in other words: the axis
  of $\rho_0(a_1)$ determines a simple closed curve on the surface, which
  intersects exactly once, transversely, the axis of $\alpha_1$, hence this
  axis defines a non separating curve).
  Thus, we can complete it into $(\alpha_1,\beta_1,\alpha_2,\beta_2)$
  represented by a system of curves on $\Sigma_2$, with
  $\alpha_1,\beta_1,\alpha_2,\beta_2\in\rho_0(\pi_1\Sigma_2)$.
  Define then $\rho_n\colon\pi_1\Sigma_2\rightarrow PSL(2,\dR)$ by the
  formulas $\rho_n(a_i)=\alpha_i$, $\rho_n(b_2)=\beta_2$,
  $\rho_n(b_1)=\beta_1\alpha_1^n$, for all $n\geq 1$.
  Then $\rho_n$ is faithful and discrete, and, as subgroups of $PSL(2,\dR)$,
  we have $Im(\rho_n)=Im(\rho_0)$. Since $\rho_0$ is purely hyperbolic
  (\textsl{i.e.}, every element of $\pi_1\Sigma_2\smallsetminus\{1\}$ is sent
  to a hyperbolic element), the elements $\alpha_1,\beta_1\in Hyp$ do not
  share any fixed points on $\partial\dH^2$. Hence
  $\Tr(\beta_1\alpha_1^n)\rightarrow+\infty$ as $n\rightarrow+\infty$.
  Now denote by $S\in Isom(\dH^2)$ the inversion with respect to the axis
  $(x,y)$ (that is, the reflection whose axis is the one of $\alpha_1$).

  We now define $h_n\colon\pi_1\Sigma_g\rightarrow Isom^+(\dH^3)$,
  for every $g\geq 3$, as follows. Consider an isometric embedding
  $i\colon\dH^2\hookrightarrow\dH^3$; this determines an injection
  $PSL(2,\dR)\hookrightarrow Isom^+(\dH^3)$. Every reflection in $\dH^2$ can
  then be realized as a rotation in $\dH^3$, and we denote again by
  $\alpha_1$, $\beta_1$, $\alpha_2$, $\beta_2$ and $S$ the corresponding
  elements of $Isom^+(\dH^3)$. We put $h_n(a_i)=\rho_n(a_i)$ and
  $h_n(b_i)=\rho_n(b_i)$ for $i=1,2$, and we put
  $h_n(a_i)=h_n(b_i)=S$ for $3\leq i\leq g$.

  We then have $h_n\in R_g(3)$, and $h_n(b_1)$ is a hyperbolic element whose
  translation distance tends to $+\infty$ as $n\rightarrow+\infty$, so that
  $\displaystyle{\mathop{\lim}_{n\rightarrow+\infty}}d(h_n)=+\infty$; hence
  there exists an accumulation point
  $(T,h_\infty)\in\partial\overline{m_g^u(3)}$ of this sequence of
  representations. We claim that no orientation on $T$ is preserved by
  $\rho_\infty$. In order to prove this, it suffices to find a non degenerate
  tripod $Trip(a,b,c,d)\subset T$, with central point $a$, and an element
  $\gamma\in\pi_1\Sigma_g$, such that $\gamma(a)=a$,
  $\gamma(b)=b$, $\gamma(c)=d$ and $\gamma(d)=c$.

  Note that if $A,B\in Hyp$ do not have any common fixed points in
  $\partial\dH^2$, then the attractive fixed point of $AB^n$, as
  $n\rightarrow+\infty$, converges to the one of $B$, whereas the repulsive
  fixed point of $AB^n$ converges to the preimage, by $A$, of the one of $B$.
  Hence, the axis of $\rho_n(b_1)$ converges to some fixed geodesic line in
  $\dH^2$. Since the images $\rho_n(\gamma)$ of the other generators
  $\gamma\in\{a_1,a_2,b_2,\ldots,a_g,b_g\}$ are fixed, there exists a sequence
  $(x_0^n)_n\in\left(\dH^2\right)^\dN$ converging to a point
  $x_\infty\in\dH^2$, such that for all $n$, $x_0^n\in\min(\rho_n)$ and
  $i(x_0^n)\in\min(h_n)$.
  Then $d(i(x_0^n),h_n(a_3)\cdot i(x_0^n))$ is bounded.
  Moreover, $h_n(b_1^{-1})=\alpha_1^{-n}\beta_1^{-1}$.
  Since the axis of the symmetry $S=h_n(a_3)$ is the axis of the translation
  $\alpha_1$, we have
  \[ d(\alpha_1^{-n}\beta_1^{-1}i(x_0^n),S.\alpha_1^{-n}\beta_1^{-1}i(x_0^n))
  =d(\beta_1^{-1}i(x_0^n),S\beta_1^{-1}i(x_0^n)), \]
  and this number is bounded, hence the distance
  $d(h_n(b_1^{-1}).i(x_0^n),h_n(a_3b_1^{-1}).i(x_0^n))$ is bounded.
  It follows from the construction of the $\dR$-tree $T$ given by M. Bestvina
  (\cite{Bestvina88}) that there exist two distinct points of $T$,
  $x_0^\infty$ and $h_\infty(b_1^{-1})\cdot x_0^\infty$, which are fixed by
  $h_\infty(a_3)$. Thus, the segment
  $\left[x_0^\infty,\rho_\infty(b_1^{-1})\cdot x_0^\infty\right]$ is globally
  fixed by $h_\infty(a_3)$. Since $\left(h_\infty(a_3)\right)^2=id_T$, in order
  to prove that there is a tripod
  $Trip(a,b,c,d)$ such that $h_\infty(a_3)(a)=a$, $h_\infty(a_3)(b)=b$,
  $h_\infty(a_3)(c)=d$ and $h_\infty(a_3)(d)=c$, it suffices to prove that
  $h_\infty(a_3)\neq id_T$. This follows for instance from the fact that
  $h_n(a_2b_1a_2^{-1})$ is a hyperbolic element whose fixed points in
  $\partial\dH^2$ are distinct from those of
  $\alpha_1$ and of $h_n(b_1)$, hence for $n$ large enough the distance between
  $h_n(a_2b_1a_2^{-1})\cdot i(x_n)$ and the axis of the symmetry $h_n(a_1)$
  is of the order of $d(h_n)$: still by the construction of M. Bestvina,
  this yields a point of $T$ which is not fixed by $h_\infty(a_3)$.
\end{preuve}

\subsubsection{The space $\overline{m_g^u}$ has at most $3$ connected components}

Now we shall exhibit another example of degeneracy.
Fix an injective representation
$\rho\colon\pi_1\Sigma_{g-1}\rightarrow PSL(2,\dR)$.
The circle $\dS^1$ being non numerable, there exists a point
$a_0\in\dS^1=\partial\dH^2$ such that for all
$\gamma\in\pi_1\Sigma_{g-1}$, $\rho(\gamma)a_0=a_0\Leftrightarrow \gamma=1$.
Choose a point $x_0\in\dH^2$. Denote by $A_n$ the hyperbolic element whose axis
passes through $x_0$, with attractive point $a_0$, and translation length $n$.
We define a representation
$\rho_n'\colon\pi_1\Sigma_g\rightarrow PSL(2,\dR)$ by letting
$\rho_n'(a_i)=\rho(a_i)$, $\rho_n'(b_i)=\rho(b_i)$ for $i\leq g-1$, and
$\rho_n'(a_g)=1$, $\rho_n'(b_g)=A_n$.

\begin{prop}\label{degg}
  The sequence $\left(\rho_n'\right)_n\in \overline{m_g^u}^\dN$ converges to an
  action $\rho_\infty$ on an $\dR$-tree $T$, which does not depend on $\rho$.
\end{prop}

In \cite{dBK}, J. DeBlois and R. Kent proved that every connected components
of $R_{g-1}$ contain injective representations.
It follows that this $\dR$-tree is a common point to all the
$\partial m_{g,k}^u$ in $\overline{m_g^u}$, for every $k\in\{0,\ldots,2g-4\}$,
since, of course, the representation $\rho_n'$ is also of
Euler class $k$ (this follows immediately from Milnor's algorithm, see formula
\ref{AlgMil}), for all $n\in\dN$.
This proves the following result:

\begin{cor}\label{precc}
  Let $g\geq 2$. Then the space $\overline{m_g^u(2)}$ has at most $3$
  connected components. More precisely, the components of Euler class between
  $0$ and $2g-4$ all meet at their boundary.
\end{cor}

Moreover, every injective representation in $R_{g-1}$ of Euler class $k$, with
$|k|\leq 2g-5$, is non elementary (indeed, elementary subgroups of $\psl$ are
virtually abelian, hence they do not contain isomorphic copies of
$\pi_1\Sigma_{g-1}$) and non discrete (by W. Goldman's corollary C
\cite{Goldman88}, faithful and discrete representations have Euler class
$2g-4$ or $4-2g$). Hence, by theorem \ref{ordrecaracterise}, every conjugacy
class of injective representations gives rise to a distinct order. These
conjugacy classes, by the theorem of DeBlois and Kent \cite{dBK}, have the
same cardinality as $\dR$. Hence,
\begin{cor}
  The surjective map $\overline{m_g^o}\twoheadrightarrow\overline{m_g^u}$
  has a fibre which has the cardinality of $\dR$.
\end{cor}

\medskip

It is known (see \cite{Goldman88}) that the space $m_g^u(3)$ has two connected
components, one containing all the representations of even Euler class in
$m_g^u(2)$, and the other containing those of odd Euler class.
The following result follows:

\begin{cor}
  Let $g\geq 3$. Then the space $\overline{m_{g}^u(3)}$ is connected.
\end{cor}

\begin{preuvede} of proposition \ref{degg}.

  We shall first define an $\dR$-tree $T$ equipped with an action
  $\rho_\infty$ of $\pi_1\Sigma_g$; we will then prove that
  $(\dH^2,\rho_n')\rightarrow (T,\rho_\infty)$ in $\overline{m_g^u}$.

  Let $\mathcal{G}$ be the Cayley graph of the free product
  $G=\pi_1\Sigma_{g-1}*\dZ$ for the generating set
  $\{a_1,b_1,\ldots,a_{g-1},b_{g-1},z\}$ where $z$ is one of the two generators
  of the free factor $\dZ$ of $G$. Note that the action of $G$ on $\mathcal{G}$
  is free, and that if $e$ is the edge between $1$ and $z$, then $e$ separates
  $\mathcal{G}$, by definition of a free product. Let $T$ be the dual graph, in
  $\mathcal{G}$, to the set $\mathcal{E}$ of middles of edges of
  $\mathcal{G}$: its vertices are the connected components of
  $\mathcal{G}\smallsetminus\mathcal{E}$, with an edge between two such
  connected components if their closures meet. We denote by
  $\overline{1}$, $\overline{z}$ the vertices of $T$ which are the images of
  the vertices $1$ and $z$ of $\mathcal{G}$, respectively. The action of $G$ on
  $\mathcal{G}$, which is free and transitive on the vertices, induces an
  action of $G$ on $T$ which is transitive on the edges. The graph $T$ is thus
  a tree, since each of its edges separates it, and the action of $G$ on $T$ is
  therefore minimal, with trivial edge stabilizers.
  Hence, it is a {\em geometric} action of the group $G$ on the $\dR$-tree $T$,
  in the sense of G. Levitt and F. Paulin (see \cite{LevittPaulin}):
  $G$ is the fundamental group of the bouquet of $\Sigma_{g-1}$ with a circle;
  denote this complex by $S$. Then the action of $G$ on $T$ is the action dual
  to the ``lamination'' defined by any point $x$ of the circle of $S$.
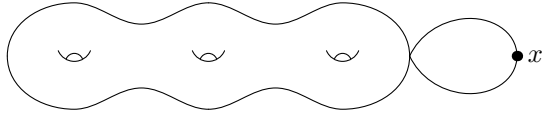
  \[ \epsfbox{Figures/FiguresCompactification.4} \]
  The surjective morphism $e\colon\pi_1\Sigma_g\rightarrow G$ defined by
  $e(a_i)=a_i$ and $e(b_i)=b_i$ if $i\leq g-1$, and $e(a_g)=1$,
  $e(b_g)=z$ induces a minimal action $\rho_\infty$ of $\pi_1\Sigma_g$ on $T$
  by isometries, in which the stabilizer of edges is exactly the normal
  subgroup of $\pi_1\Sigma_g$ generated by $a_g$.

  Now we shall prove that $(\dH^2,\rho_n')\rightarrow(T,\rho_\infty)$, in the
  sense of the equivariant Gromov topology.
  Let $P$ be a finite subset of $\pi_1\Sigma_g$. Let $T_P$ be the minimal
  subtree of $T$ containing the vertices
  $\overline{1},\overline{z}\in C\subset T$, as well as $P\cdot\overline{1}$
  and $P\cdot\overline{z}$. Let $\{x_1,\ldots,x_m\}$ be a finite collection in
  $T_P$ and let $\varepsilon>0$. We need to prove that there exists
  $\{x_1',\ldots,x_m'\}\subset\dH^2$ such that for all $i,j\in \{1,\ldots,m\}$
  and all $\gamma_1,\gamma_2\in P$,
  \[ \left|
  d\left(\rho_n'(\gamma_1)\cdot x_i',\rho_n'(\gamma_2)\cdot x_j'\right)-
  d\left(\rho_\infty(\gamma_1)\cdot x_i,\rho_\infty(\gamma_2)\cdot x_j\right)
  \right|<\varepsilon. \]

  For every $\gamma\in P$, the element $e(\gamma)\in G$ can be uniquely written
  as $e(\gamma)=u_1z^{n_1}\cdots u_k z^{n_k}$,
  with $u_1,\ldots,u_k\in\pi_1\Sigma_{g-1}$, $u_2$, \ldots, $u_k\neq 1$
  and $n_1$, \ldots, $n_{k-1}\neq 0$.
  Denote by $Q\subset\pi_1\Sigma_{g-1}$ the set described by all the elements
  $u_1,\ldots,u_k$ as $\gamma$ describes $P$.
  Put
  \[\alpha_P=\displaystyle{\mathop{\min_{u_1,u_2\in Q}}_{u_1\neq u_2}
  (\widehat{\rho_0(u_1) a_0,x_0,\rho_0(u_2) a_0})}.\]
  It is a non zero angle. Put also
  $\beta_P=\displaystyle{\max_{u\in Q}}d_{\dH^2}(x_0,\rho_0(u)\cdot x_0)$.

  If $P$ is large enough, which we suppose now, the tree $T_P$ is the orbit of
  the segment $[\overline{1},\overline{z}]$ under $\rho_\infty(P)$.
  For all $i\in\{1,\ldots,m\}$ there exists $t_i\in[0,1]$ and (at least one)
  $c_i\in P$ such that $x_i$ is on the segment
  $[c_i\cdot\overline{1},c_i\cdot\overline{z}]$ and such that
  $d(x_i,c_i\cdot\overline{1})=t_i=1-d(x,c_i\cdot\overline{z})$.
  Denote by $x_i'$ the element of $\dH^2$ such that
  $\frac{1}{n}d_{\dH^2}(x_i',\rho_n'(c_i)x_0)=t=1-\frac{1}{n}d_{\dH^2}
  (x_i',\rho_n'(c_i)A_nx_0)$. Such a point exists and is unique, since
  $d_{\dH^2}(\rho_n'(c_i)x_0,\rho_n'(c_i)A_nx_0)=n$.
  By construction of the points $x_i'$, in order to check that for all
  $i,j\in\{1,\ldots,m\}$ and $\gamma_1,\gamma_2\in P$ we have
  \[ \left|
  d\left(\rho_n'(\gamma_1)\cdot x_i',\rho_n'(\gamma_2)\cdot x_j'\right)-
  d\left(\rho_\infty(\gamma_1)\cdot x_i,\rho_\infty(\gamma_2)\cdot x_j\right)
  \right|<\varepsilon, \]
  it suffices to do it in the case when $x_i$ and $x_j$ are vertices of the
  tree $T$.
  Suppose for instance that $x_i=c_i\cdot\overline{1}$ and
  $x_j=c_j\cdot\overline{1}$ with $c_i,c_j\in P$. Let
  $\gamma_1,\gamma_2\in P$. Denote $c=\gamma_1 c_1(\gamma_2 c_2)^{-1}$ and
  $e(c)=u_{k+1} t^{n_k}\cdots u_2 t^{n_1} u_{1}$, with $n_1,\ldots,n_k\neq 0$
  and $u_2,\ldots,u_k\neq 1$.
  Since the stabilizer of the edge $\overline{1}$, in the group $G$,
  is $\pi_1\Sigma_{g-1}$ and since $d(\overline{1},\overline{z})=1$ in the
  {\em tree} $T$, we check that
  $d\left(\rho_\infty(\gamma_1)x_i,\rho_\infty(\gamma_2)x_j\right)=
  \displaystyle{\sum_{i=1}^k}|n_i|$.
  We are going to prove that, asymptotically, the distance
  $\frac{1}{n}d_{\dH^2}\left(\rho_n'\left(c\right)x_0,x_0\right)$
  approaches $\displaystyle{\sum_{i=1}^m}|n_k|$.
  We work by induction on $k$.
  Denote $c'=u_{k} z^{n_{k-1}}\cdots u_2 z^{n_1}$, and suppose that
  $u_{k}\neq 1$ and that
  $d_{\dH^2}(\rho_n'(c)x_0,x_0)=n\displaystyle{\sum_{i=1}^{k-1}}|n_i|+O(1)$.
  Then $d(\rho_n'(c')x_0,x_0)=n\displaystyle{\sum_{i=1}^{k-1}}|n_i|$, and hence
  $d(A_n^{n_k}\rho_n'(c')x_0,A_n^{n_k}x_0)=
  n\displaystyle{\sum_{i=1}^{k-1}}|n_i|$, and
  $d(A_n^{n_k}x_0,x_0)=n|n_k|$. Since $A_n^{n_k}\in PSL(2,\dR)$
  preserves the angles, we then have
  $\displaystyle{\widehat{A_n^{n_k}\rho_n'(c')x_0,A_n^{n_k}x_0,x_0}=
  \widehat{\rho_n'(c')x_0,x_0,A_n^{-n_k}x_0}\geq\dfrac{\alpha_{P'}}{2}}$, for
  $P'$ large enough (but depending only on $P$).
  This, combined with cosine law I, implies that
  $d(A_n^{n_k}\rho_n'(c')x_0,x_0)=n\displaystyle{\sum_{i=1}^{m}}|k_i|+O(1)$.
  Since $d(u_{k+1}x_0,x_0)\leq \beta$, it follows that
  $d(\rho_n'(c)x_0,x_0)=n\displaystyle{\sum_{i=1}^{k}}|n_i|+O(1)$, which
  completes the proof.
\end{preuvede}

\subsubsection{The space $\overline{m_g^u}$ is connected}

We will now prove theorem \ref{introcc}.

\begin{theo}\label{cc}{\ }
  \begin{itemize}
    \item
      For all $g\geq 2$, the space $\overline{m_g^u}$ has at most two connected
      components. More precisely, all the connected components, except possibly
      the one of Euler class $2g-3$, meet at their boundaries.
    \item
      For all $g\geq 4$, the space $\overline{m_g^u}$ is connected.
  \end{itemize}
\end{theo}

We still consider the generating set $S=\{a_1,\ldots,b_g\}$ of $\pi_1\Sigma_g$.
The main idea is the following.

\begin{lem}\label{cctech}
  Let $\Gamma$ be a finitely generated group
  and let $\phi\colon\Gamma\rightarrow\psl$ be a discrete, faithful
  representation of cocompact image.
  Denote by $B_n$ the closed ball of radius $n$ for the Cayley metric on
  $\pi_1\Sigma_g$ for the generating set $S$, and suppose that
  $\phi_n\colon\pi_1\Sigma_g\rightarrow\Gamma$ is a non injective morphism
  such that $Ker(\phi_n)\cap B_n=\{1\}$. Then
  $d(\phi\circ\phi_n)\rightarrow+\infty$ as $n\rightarrow+\infty$.
  Let $(T,\rho_\infty)$ be an accumulation point in $\partial\overline{m_g^u}$
  of the sequence $\left(\phi\circ\phi_n\right)_{n\in\dN}$. Then the action
  $\rho_\infty$ has small edge stabilizers.
\end{lem}

\begin{preuve}
  First, let us prove that
  $\displaystyle{\lim_{n\rightarrow+\infty}}d(\phi\circ\phi_n)=+\infty$.
  By contradiction, suppose that we can extract a subsequence such that
  $\left(d(\phi\circ\phi_{\varphi(n)})\right)_{n\in\dN}$ converges to a real
  number $d\in\dR_+$. Fix a point $x_0\in\dH^2$. Since the image $\phi(\Gamma)$
  is cocompact, there exist $g_{\varphi(n)}\in\phi(\Gamma)$ and
  \[ x_n\in\min\left(g_{\varphi(n)}\cdot(\phi\circ\phi_{\varphi(n)})\cdot
  g_{\varphi(n)}^{-1}\right) \]
  such that the distance $d(x_0,x_n)$ is unbounded; say $d(x_0,x_n)\leq k$.
  Denote by $\rho_n$ the representation
  $g_{\varphi(n)}\cdot(\phi\circ\phi_{\varphi(n)})\cdot
  g_{\varphi(n)}^{-1}$. Then for every $n\geq 0$, $\rho_n$ is discrete, and
  $Ker(\rho_n)\cap B_{\varphi(n)}=\{1\}$. For every $n\in\dN$ and every
  $\gamma\in S$, $d(\rho_n(\gamma)\cdot x_0,x_0)\leq d(\rho_n)+2k$,
  and $\lim_n d(\rho_n)=d$ hence, by proposition \ref{bornecompact}, up to
  extract it, $(\rho_n)_n$ converges to a representation $\rho\in R_g$.
  Since the set $\rho_n(\pi_1\Sigma_g)=\phi(\Gamma)\subset\psl$ is discrete
  and since for all $\gamma\in\pi_1\Sigma_g\smallsetminus\{1\}$, for $n$ large
  enough $\rho_n(\gamma)\neq\Id$, we have
  $\rho(\gamma)\neq\Id$, hence $\rho$ is injective. Hence, the representation
  $\rho$ is discrete and faithful, hence $|e(\rho)|=2g-2$, thus $\rho$ is the
  limit of representations $\rho_n$ which are non injective, in particular
  $|e(\rho_n)|\neq 2g-2$: this is in contradiction with the continuity of the
  Euler class (in fact it is not necessary to use the Euler class here, but it
  gives the shortest proof).

  Now denote $\rho_n=\phi\circ\phi_n$. It follows that there exists an
  accumulation point $(T,\rho_\infty)\in\partial\overline{m_g^u}$ of the
  sequence $\left(\rho_n\right)_{n\in\dN}$.
  We still have to prove that this action has small edge stabilizers.
  But our representation $\rho_n$, for all $n\geq 0$, is {\em discrete}. We
  can therefore apply the same argument as M. Bestvina and F. Paulin
  (\cite{Bestvina88,Paulin88}), consisting of applying Margulis' lemma.
  Here we follow the lines of the proof of theorem 6.7 of \cite{Paulin88},
  pages 78-79, and refer the reader to this text for more details.

  \noindent {\bf Margulis' lemma:} {\em There exists a constant $\mu>0$,
  depending only on $n$, such that for every discrete group $\Gamma$ of
  isometries of $\dH^n$, and for all $x\in\dH^n$, the subgroup generated by
  $\{\gamma\in\Gamma\,|\,d(x,\gamma x)<\mu\}$ is virtually abelian.}

  Let us suppose that there exists a segment $[x_\infty,y_\infty]$ of the
  limit $\dR$-tree, whose stabilizer contains a free group of rank $2$. Up to
  considering a subgroup of index $2$, we may suppose that there exists a free
  group of rank $2$, $\langle\alpha,\beta\rangle\subset\pi_1\Sigma_g$, such
  that $\alpha$ and $\beta$ fix $x_\infty$ and $y_\infty$. Hence, for all
  $\varepsilon>0$, and for $n$ large enough, there exist $x_n,y_n\in\dH^2$ such
  that $\left|\frac{1}{\ell(\rho_n)}d_{\dH^2}(x_n,y_n)-1\right|<\varepsilon$,
  $d_{\dH^2}(x_n,\rho_n(\alpha)x_n)<\varepsilon\ell(\rho_n)$,
  $d_{\dH^2}(y_n,\rho_n(\alpha)y_n)<\varepsilon\ell(\rho_n)$,
  and similarly for $\rho_n(\beta)$. Denote by $z_n$ the middle of the segment
  $[x_n,y_n]$: it can be proved (see \cite{Paulin88}) that if $\varepsilon$ is
  small enough and if $\ell(\rho_n)$ is large enough, then the elements
  $\left[\rho_n(\alpha),\rho_n(\beta)\right]$ and
  $\left[\rho_n(\alpha^2),\rho_n(\beta)\right]$ move the point $z_n$ by a
  distance less than $\mu$. By Margulis' lemma, the elements
  $\rho_n([\alpha,\beta]),\rho_n([\alpha^2,\beta])\in\psl$ generate a virtually
  abelian subgroup of $\psl$. But it is an easy exercise (see \textsl{e.g.}
  \cite{TheseMaxime}, lemme 1.1.18) to check that virtually abelian subgroups
  of $\psl$ are {\em metabelian}, \textsl{i.e.}, all the commutators commute.
  In particular, denote for instance
  $w(\alpha,\beta)=\left[[[\alpha,\beta],[\alpha^2,\beta]],
  [[\alpha,\beta]^2,[\alpha^2,\beta]]\right]$:
  then we have $\rho_n(w(\alpha,\beta))=1$ for all $n$ large enough. If $n$ is
  large enough, and larger than the length of the word $w(\alpha,\beta)$ in the
  generators $a_i,b_i$, this implies that $w(\alpha,\beta)=1$ in
  $\pi_1\Sigma_g$ (since $Ker(\rho_n)\cap B_n=\{1\}$), which contradicts the
  fact that $\langle\alpha,\beta\rangle$ is free. Hence, the action
  $\rho_\infty$ indeed has small edge stabilizers.
\end{preuve}

\begin{preuvede} of theorem \ref{cc}.
  \begin{itemize}
    \item
      Fix a cocompact Fuchsian group $\Gamma$, a subgroup $\dF_2\subset\Gamma$
      isomorphic to a free group of rank $2$, and
      $\phi\colon\Gamma\rightarrow\psl$ the tautologic representation
      (the inclusion).
      Consider the morphism $\phi_n\colon\pi_1\Sigma_g\rightarrow\dF_2$ given
      by theorem \ref{lemtechMCG1}, and let $\rho_n=\phi\circ\phi_n$. Then
      $\rho_n$ factors through the free group, which has a trivial $H^2$, hence
      by remark \ref{pscohom}, the representation $\rho_n$ has Euler class
      zero.
      By lemma \ref{cctech}, $\left(\rho_n\right)_{n\in\dN}$ possesses
      an accumulation point $(T,\rho_\infty)\in\partial\overline{m_g^u}$,
      which has small edge stabilizers.
      By R. Skora's theorem \cite{Skora}, this implies that this limit is also
      at the boundary of the Teichm\"uller space; in other words,
      this action on an $\dR$-tree is also the limit of representations of
      Euler class $2g-2$. Hence, the closures of the connected components of
      $m_g^u$ of Euler classes $0$ and $2g-2$ meet. By corollary \ref{precc},
      we already knew that the connected components of Euler classes
      $0$, $1$, \ldots, $2g-4$ meet at their boundaries; this concludes the
      proof of the first point.
    \item
      By proposition 4.5 of \cite{noninjrep}, the map
      $p_g\colon\pi_1\Sigma_g\rightarrow\psl$ that we defined in section
      \ref{SectionPresqueInj} is discrete and has Euler class $2g-3$.
      We have seen (lemma \ref{lemtechMCG2}) that for all $n\geq 0$, there
      exists $\phi_n\in Aut(\pi_1\Sigma_g)$ such that
      $Ker(p_g\circ\phi_n)\cap B_n=\{1\}$. The representation
      $p_g\circ\phi_n$ is still discrete, and we have
      $|e(p_g\circ\phi_n)|=2g-3$, hence, as before, the lemma
      \ref{cctech} ensures that the connected components of $m_g^u$ of Euler
      classes $2g-3$ and $2g-2$ meet at their boundaries.
  \end{itemize}
\end{preuvede}

\begin{ps}\label{toutlemondetoucheteich}
  For all $g\geq 2$, consider the action $(\rho_\infty,T)$ exhibited in
  proposition \ref{degg}.
  The representation $\rho_\infty\colon\pi_1\Sigma_g\rightarrow Isom(T)$
  factors through the group $\pi_1\Sigma_{g-1}*\dZ$, which acts on $T$ with
  trivial arc stabilizers. Since there exists a morphism
  $\pi_1\Sigma_{g-1}*\dZ\twoheadrightarrow\dF_2$ with non abelian image
  to some free group of rank $2$, as in section \ref{SectionPresqueInj} we can
  prove that there exist automorphisms $\phi_n$ of $\pi_1\Sigma_g$ such that
  $Ker(\rho_\infty\circ\phi_n)\cap B_n=\{1\}$. Following the proof of
  lemma 5.7 of \cite{Paulin88} (see also remark (1), page 73 in
  \cite{Paulin88}), we can prove that up to extract it, the sequence
  $(\rho_\infty\circ\phi_n,T)$ converges to an action on an $\dR$-tree
  with small stabilizers. This proves that for all $g\geq 2$ and all
  $k\in\{0,\ldots,2g-4\}$,
  $\overline{m_{g,k}^u}\cap\overline{m_{g,2g-2}^u}\neq\emptyset$.
  
  Hence, for all $g\geq 4$ and all $k\in\{0,\ldots,2g-3\}$, we have
  $\overline{m_{g,k}^u}\cap\overline{m_{g,2g-2}^u}\neq\emptyset$.
\end{ps}

\subsubsection{Dynamics}
Finally, here we complete the proof of theorem \ref{PaBoIntro},
which implies that the compactification $\overline{m_g^u}$ is
{\em extremely} degenerated, pledging that $\overline{m_g^o}$ is a more
natural object to consider and to try to generalize to higher dimensions
and to the compactification of other representation spaces.

\begin{prop}\label{PaBo}
  Let $g\geq 4$ and $k\in\{0,\ldots,2g-3\}$.
  Then the boundary of the Teichm\"uller space embeds in
  $\partial m_{g,k}^u\subset\overline{m_g^u}$ as a
  closed, nowhere dense subset.
\end{prop}

\begin{preuve}
  Fix $g\geq 4$ and $k\in\{0,\ldots,2g-3\}$.
  Put $F_k=\partial m_{g,2g-2}^u\cap\partial m_{g,k}^u$.
  By remark \ref{toutlemondetoucheteich}, we have $F_k\neq\emptyset$.
  Since $m_{g,2g-2}^u$ and $m_{g,k}^u$ are invariant under the (natural)
  action of $Out(\pi_1\Sigma_g)$, it follows that $F_k$ is invariant, too,
  under this action. It is well-known
  (see \cite{FLP}, expos\'e 6, theorem VII.2, p. 117 ; see also \cite{Masur})
  that the action of $Out(\pi_1\Sigma_g)$ on $\partial m_{g,2g-2}^u$ is
  minimal, that is, every closed subset of $\partial m_{g,2g-2}^u$, invariant
  under $Out(\pi_1\Sigma_g)$, is either empty of is $m_{g,2g-2}^u$ itself.
  Since $F_k$ is a closed subset, this implies that
  $\partial m_{g,2g-2}^u=F_k\subset\partial m_{g,k}^u$.

  Now denote by $G_k$ the boundary of $\partial m_{g,2g-2}^u$ in the space
  $\partial m_{g,k}^u$.
  We can easily produce elements in $\partial m_{g,k}^u$ which do not have
  small stabilizers (if $k\leq 2g-4$ then the tree $(\rho_\infty,T_\infty)$
  of proposition \ref{degg} is an example; if $k=2g-3$ then we can compose
  the map $p_g$ with Dehn twists along the last handle: this does not touch
  the kernel of the map $p_g$, hence this yields a sequence of actions
  with a fixed non trivial kernel in $\pi_1\Sigma_g$, converging (up to
  extract it) to an action on an $\dR$-tree, with this non trivial kernel).
  Hence,
  $\partial m_{g,k}^u\neq \partial m_{g,2g-2}^u$. By proposition \ref{oneend},
  the space $m_{g,k}^u$ is connected: it follows that $G_k\neq\emptyset$.
  Since $\partial m_{g,2g-2}^u$ is closed, we have
  $G_k\subset\partial m_{g,2g-2}^u$, and $G_k$ is again invariant under
  the action of the mapping class group. Hence, $G_k=\partial m_{g,2g-2}^u$.
\end{preuve}

\bibliographystyle{plain}

\bibliography{BiblioThese,BibMain}

\end{document}